\documentclass[10pt]{article}

\usepackage{amsmath,amssymb,amsthm,amsfonts}
\usepackage{mathtools}
\usepackage{graphicx,graphics}
\usepackage{latexsym,bm}
\usepackage{enumitem}
\usepackage{esint}
\usepackage{microtype}
\usepackage{mathrsfs}
\usepackage[hidelinks,hypertexnames=false]{hyperref}
\ifdefined\RenewDocumentCommand\else
  \usepackage{xparse}
\fi
\ifdefined\NewCommandCopy\else
  \usepackage{letltxmacro}
\fi
\numberwithin{equation}{section}
\allowdisplaybreaks[4]

\newtheoremstyle{manuscriptplain}
  {\topsep}{\topsep}{\itshape}{}{\normalfont\bfseries}{.}{.5em}
  {\thmname{#1}\thmnumber{ \textup{#2}}\thmnote{ \textnormal{(#3)}}}
\theoremstyle{manuscriptplain}
\setlist[enumerate]{font=\normalfont}
\AtBeginDocument{%
  \ifdefined\NewCommandCopy
    \NewCommandCopy\manuscriptref\ref
  \else
    \LetLtxMacro\manuscriptref\ref
  \fi
  \RenewDocumentCommand{\ref}{s m}{%
    \IfBooleanTF{#1}
      {\textup{\manuscriptref*{#2}}}
      {\textup{\manuscriptref{#2}}}%
  }%
}

\makeatletter
\newcommand{\proofstep}[2]{%
  \if@inlabel\leavevmode\par\nobreak\fi
  \par\addvspace{\smallskipamount}\noindent
  {\normalfont\bfseries Step~#1.}\ {\normalfont #2}%
  \par\nobreak\noindent\ignorespaces
}
\makeatother

\newtheorem{lemma}{Lemma}[section]
\newtheorem{definition}{Definition}
\newtheorem{proposition}{Proposition}[section]
\newtheorem{corollary}{Corollary}[section]
\newtheorem{remark}{Remark}[section]
\newtheorem{theorem}{Theorem}[section]
\newcommand{\ve}{\varepsilon}
\newcommand{\epstr}{\epsilon_{\mathrm{tr}}}
\newcommand{\R}{{\mathbb R}}
\newcommand{\C}{{\mathbb C}}

\providecommand{\g}{\mathbf{g}}
\renewcommand{\L}{\mathcal{L}}
\providecommand{\ul}{{\rm ul}}

\begin{document}

\date{}

\title{Boundary singularities for the Lane--Emden equation}

\author{Weiwei Ao, Sho Katayama, Juncheng Wei}
\maketitle

\begin{abstract}
For every $n\ge3$ and $(n+1)/(n-1)<p<(n+2)/(n-2)$, we construct
positive solutions of $\Delta u+u^p=0$ in the half-space with zero Dirichlet
data away from the origin, the homogeneous singular profile there, and
Poisson-kernel decay at infinity, resolving a question of
del Pino--Musso--Pacard. Under a weighted non-degeneracy condition,
these half-space solutions serve as building blocks for positive very
weak solutions in smooth bounded domains, singular
along prescribed compact smooth boundary submanifolds of normal dimension
$n$. We verify non-degeneracy in the initial stable interval and at
and just beyond its threshold when this lies below $(n+2)/(n-2)$.
For isolated boundary singularities, the construction applies throughout
the strict subcritical interval outside a locally finite exceptional set.
\end{abstract}

\setcounter{tocdepth}{2}
\tableofcontents
\clearpage

\section{Introduction and main results}\label{s1}

We consider positive solutions of the Lane--Emden--Fowler problem
\begin{equation}\label{e1}
\left\{\begin{array}{ll}
\Delta u+u^p=0 & \text{ in }\Omega,\\
u>0 & \text{ in }\Omega,\\
u=0 & \text{ on }\partial\Omega,
\end{array}\right.
\end{equation}
where $\Omega\subset \R^N$ is a bounded smooth domain, $N\geq 2$, and $p>1$. The solutions are smooth in $\Omega$ and vanish on $\partial\Omega$ away from a prescribed boundary singular set, where they may blow up. They are understood in the very weak sense described below.

The role of singularities is well understood in the interior case. If $u\in L^p(\Omega)$ is a positive distributional solution of \eqref{e1}, then singularities are removable for $1<p<\frac{N}{N-2}$, with $N/(N-2)=+\infty$ when $N=2$. Above this threshold, the local singular models were analyzed in \cite{gs}, and solutions with prescribed interior singular sets were constructed in \cite{cl,mp}. At the threshold, Pacard \cite{p} constructed singular solutions using a model with a logarithmic correction. In dimension $n\ge3$, the model equation is
\begin{equation}\label{e2}
\left\{
\begin{array}{ll}
\Delta u+u^p=0 & \text{in }\R^n\setminus\{0\},\\
u>0, & \displaystyle\lim_{|x|\to0}u(x)=+\infty.
\end{array}
\right.
\end{equation}
For
$\frac{n}{n-2}<p<\frac{n+2}{n-2}$,
the singular radial theory consists of the homogeneous solution
$u_{\mathrm{hom}}=c_{p,n}|x|^{-2/(p-1)}$
and a one-parameter family of solutions which are singular like $u_{\mathrm{hom}}$ at the origin and have fundamental-solution decay at infinity; see \cite{gs}. These fast-decaying radial solutions are used in the interior gluing construction of Mazzeo--Pacard \cite[Section~\textup{4.1}]{mp}.
For nonlocal singular models, related ODE methods for asymptotic
analysis and radial non-degeneracy modulo scaling are developed in
\cite{acdfgw-ode}.

For boundary singularities, the critical exponent and the notion of solution are different. Following del Pino--Musso--Pacard \cite[Section~\textup{1}]{dmp}, a positive function $u\in C^\infty(\Omega)$ is a very weak solution of \eqref{e1} if
\begin{equation}
u\in L^1(\Omega),\qquad \operatorname{dist}(x,\partial\Omega)u^p\in L^1(\Omega),
\end{equation}
and
\begin{equation}
\int_\Omega (u\Delta v+u^pv)\,dx=0,
\end{equation}
for every $v\in C^2(\overline\Omega)$ with $v=0$ on $\partial\Omega$.
As recalled in \cite[Section~\textup{1}]{dmp}, in the range
$1<p<\frac{N+1}{N-1}$,
every positive very weak solution belongs to $H^1_0(\Omega)$ and solves
the equation in the variational sense. It follows from elliptic regularity
that $u\in C^2(\overline\Omega)$.
Thus $\frac{N+1}{N-1}$ is the boundary analogue of the interior critical exponent $\frac{N}{N-2}$ when $N\ge3$. Boundary singular solutions may therefore occur only when
$p\ge\frac{N+1}{N-1}$.

This threshold also appears in weighted smoothing estimates for fast
diffusion: Choi--Jiang--Wang--Xiong \cite{jwx} identify
$p=(N+1)/(N-1)$ as the sharp exponent threshold for the natural
weighted $L^p$--$L^\infty$ smoothing effect, with boundary singular
elliptic solutions entering the sharpness argument. As noted in
\cite[Remark~\textup{3.6}]{jwx}, the exponent range of their singular
counterexamples is tied to the existence of boundary singular elliptic
solutions in bounded domains. For the elliptic source equation,
isolated strong boundary singularities and their asymptotic
classification have been studied in \cite{bpv1,lb}; see also
\cite{bdn,mr,s} for related very weak regularity and a priori results.

Boundary singularities also arise in other settings. Hor\'ak, McKenna,
and Reichel \cite{hmr} constructed unbounded very weak solutions at
conical corners, and Gkikas \cite{gk} constructed solutions along
prescribed edges of wedge-like domains, with fast decay at infinity
in the unbounded case. There the singular set is tied to the boundary
geometry. Problems with prescribed boundary measures and generalized
traces are treated in \cite{bvhnv,bvv}, with singular-potential
extensions in \cite{bmn,gkn}; their solvability is formulated through
measure traces and capacity conditions. Here we prescribe a singular
set on a smooth boundary and construct solutions with the strong
nonlinear singular profile and zero boundary value elsewhere.

For smooth boundaries, del Pino--Musso--Pacard \cite{dmp} constructed
such solutions at the boundary critical exponent and for exponents
slightly above it. Their construction uses a fast-decaying half-space
solution in the variables normal to the singular set. If the boundary singular set has
dimension $k$ in $\R^N$, these variables have dimension $n=N-k$;
for isolated singularities, $k=0$ and $n=N$. Throughout our gluing
construction, $n\ge3$.

In this normal dimension the critical exponent is $(n+1)/(n-1)$.
In analogy with the interior problem~\eqref{e2}, we consider, for
$p>(n+1)/(n-1)$, the half-space problem
\begin{equation}\label{e3}
\left\{\begin{array}{ll}
\Delta u+u^p=0 & \text{ in }\R^n_+\setminus\{0\},\\
u>0 & \text{ in }\R^n_+,\\
u=0 & \text{ on }\partial\R^n_+\setminus\{0\}.
\end{array}\right.
\end{equation}
We interpret the zero boundary value at the origin in the local very weak
sense: a cell belongs to $L^1_{\rm loc}(\overline{\R^n_+})$, its source
satisfies $x_nu^p\in L^1_{\rm loc}(\overline{\R^n_+})$, and
\begin{equation}\label{eq:halfspace-very-weak}
 \int_{\R^n_+}(u\Delta\zeta+u^p\zeta)\,dx=0,
 \qquad
 \zeta\in C_c^2(\overline{\R^n_+}),\quad
 \zeta|_{\partial\R^n_+}=0.
\end{equation}
For the singular profiles below, $a_p=2/(p-1)<n-1$ ensures the local
integrability. The identity follows by cutoff approximation at the origin:
since $|\zeta|\le Cx_n$, $|u|\le C|x|^{-a_p}$ and, by scaled elliptic
estimates, $|\nabla u|\le C|x|^{-a_p-1}$, the boundary terms on
$\{|x|=\rho,\ x_n>0\}$ are $O(\rho^{n-1-a_p})\to0$.
For $\frac{n+1}{n-1}<p<\frac{n+2}{n-2}$, this problem admits a homogeneous singular solution
\begin{equation}\label{e4}
u_\infty(x)=|x|^{-\frac{2}{p-1}}\phi_p(\theta),\qquad \theta=\frac{x}{|x|},
\end{equation}
where $\phi_p$ solves
\begin{equation}\label{e5}
\left\{\begin{array}{ll}
\displaystyle
\left(\Delta_{S^{n-1}}+n-1\right)\phi_p
-\frac{p+1}{p-1}\left(n-\frac{p+1}{p-1}\right)\phi_p
+\phi_p^p=0 & \text{ in }S^{n-1}_+,\\
\phi_p>0 & \text{ in }S^{n-1}_+,\\
\phi_p=0 & \text{ on }\partial S^{n-1}_+,
\end{array}\right.
\end{equation}
with
\[
S^{n-1}_+=\{\theta=(\theta_1,\ldots,\theta_n)\in S^{n-1}:\theta_n>0\}.
\]
By \cite[Theorem~\textup{1.1}]{bpv1}, the spherical problem \eqref{e5} has a unique positive solution in the range relevant here, and \cite[Theorem~\textup{1.2}]{bpv1} shows that this profile governs every non-removable isolated boundary singularity in the strict subcritical regime.

The homogeneous solution $u_\infty$ describes the desired singularity,
but rescaling does not change it. For gluing, we need a solution with
the same behavior near the origin and a faster decay at infinity.
Its rescalings then retain the singular profile while becoming small
away from the singular point, so that they can be cut off to form an
approximate solution. The required half-space cell is a solution $u_1$
of \eqref{e3} satisfying
\begin{equation}\label{e6}
u_1(x)=\left\{\begin{array}{ll}
|x|^{-\frac{2}{p-1}}\phi_p(\theta)(1+o(1)) & \text{ as } |x|\to0,\\
|x|^{-(n-1)}\theta_n(1+o(1)) & \text{ as } |x|\to\infty.
\end{array}\right.
\end{equation}
Del Pino--Musso--Pacard \cite{dmp} proved the existence of such a solution when $p>\frac{n+1}{n-1}$ is sufficiently close to the lower endpoint and used it as the basic cell in their gluing construction. In \cite[Section~\textup{1.2}, Open problem~1]{dmp}, they asked whether this
cell exists throughout the strict subcritical interval. Our first theorem answers this question affirmatively for $n\ge3$, with the same
singular profile and normalized Poisson-kernel decay as
\cite[Proposition~\textup{1.1}]{dmp}.

\begin{theorem}\label{thm1}
Let $n\ge3$ and
\[
\frac{n+1}{n-1}<p<\frac{n+2}{n-2}.
\]
Then problem~\eqref{e3} admits a solution $u_1$ satisfying
\[
 u_1(x)=
 \begin{cases}
 |x|^{-\frac{2}{p-1}}\phi_p(x/|x|)(1+o(1)),& |x|\to0,\\[1mm]
 |x|^{-n}x_n(1+o(1)),& |x|\to\infty,
 \end{cases}
\]
and
\[
0<u_1(x)\le c_p |x|^{-\frac{2}{p-1}}\frac{x_n}{|x|}
\qquad (x\in\R^n_+),
\]
for some $c_p>0$.
\end{theorem}

Related half-space solutions with nonzero singular boundary data
were constructed in \cite{fik}. The second author \cite{kat-cone} studied
existence and multiplicity for the Lane--Emden equation on cones,
including the half-space, with inhomogeneous Dirichlet boundary data.
These boundary conditions differ from \eqref{e3}.

At the Sobolev critical exponent, Xiong \cite{x} proved quantitative asymptotics
and symmetry, while Wang--Xiong \cite{wx} constructed a global
continuum of positive Delaunay-type log-periodic half-space solutions,
addressing \cite[Section~\textup{1.2}, Open problem~2]{dmp}.

The fast-decaying solution allows us to extend the boundary singularity
construction of \cite{dmp}. The remaining question is whether the error
of the approximate solution can be corrected. This requires estimates
for the linearized equation, and hence a non-degeneracy condition on
the chosen cell. For isolated singularities, Assumption H\textup{(i)} excludes
nonzero solutions of
\[
L_1\Phi:=\Delta_X\Phi+pu_1^{p-1}\Phi=0
\]
with zero boundary value and the prescribed bounds at zero and infinity.
These bounds, stated in Section~\ref{sec:assumption-H}, specify the class in which
the kernel is required to vanish.

We first state the result for isolated boundary singularities.
Here and below, convergence at $y\in\partial\Omega$ is non-tangential,
that is, it holds within every cone with vertex $y$, axis the inward
normal, and aperture strictly smaller than $\pi/2$.
\begin{theorem}\label{thm2}
Let $N\ge3$, let $\Omega\subset\R^N$ be a bounded smooth domain, and let
\[
\frac{N+1}{N-1}<p<\frac{N+2}{N-2}.
\]
Suppose that problem~\eqref{e3} in $\R^N_+$ admits a fast-decaying solution $u_1^{(p)}$ satisfying Assumption~H\textup{(i)} for some pair of nonindicial weights $(\mu,\nu)$ as in~\eqref{e59-twoend}. Then, for every finite set
\[
S=\{y_1,\ldots,y_\ell\}\subset\partial\Omega,
\qquad y_i\ne y_j\quad(i\ne j),
\]
problem~\eqref{e1} admits a positive very weak solution $u$ such that
\[
u\in C^2(\overline\Omega\setminus S),
\qquad
u(x)\to+\infty
\quad\text{non-tangentially as }x\to y_i,
\quad i=1,\ldots,\ell.
\]
\end{theorem}

The same construction applies when the singular set is a boundary submanifold, in analogy
with the interior results of \cite{mp,p}. In this case rescaling also
introduces tangential variables $Y\in\R^k$. The limiting linearized
equation is
\[
(\Delta_Y+L_1)\Phi=0\quad\text{in }\R^k\times\R^n_+.
\]
Assumption H\textup{(ii)} excludes nonzero solutions of this equation with zero boundary
value, the prescribed normal bounds, and polynomial growth in $Y$.
It is this condition that is needed when the singular set has positive
dimension. The precise hypotheses and admissible weights are given
in Section~\ref{sec:assumption-H}; their compatibility with the gluing weights
is proved in Remark~\ref{rem:weight-compatibility}.
\begin{theorem}\label{thm3}
Let $\Omega\subset\R^N$ be a bounded smooth domain and let
\[
0<k<N-2,\qquad n:=N-k,\qquad
\frac{n+1}{n-1}<p<\frac{n+2}{n-2}.
\]
Assume that \eqref{e3} in $\R^n_+$ admits a fast-decaying solution $u_1^{(p)}$ satisfying Assumption~H\textup{(ii)} for some pair of nonindicial weights $(\mu,\nu)$ as in~\eqref{e59-twoend}. Then, for every smooth compact embedded $k$-dimensional submanifold $S\subset\partial\Omega$, \eqref{e1} admits a positive very weak solution $u$ such that
\[
u\in C^2(\overline\Omega\setminus S),
\qquad
u(x)\to+\infty
\quad\text{non-tangentially as }x\to y,
\quad y\in S.
\]
\end{theorem}

The non-degeneracy hypotheses can be verified for the following
families of cells. Let $p^*(n)$ be the right endpoint of the connected
stable interval of the homogeneous solution issuing from
$(n+1)/(n-1)$; see Section~\ref{sec:stability}. The cells constructed in this
interval satisfy both H\textup{(i)} and H\textup{(ii)}. When the threshold lies below
$(n+2)/(n-2)$, these properties also hold at the threshold and on a
small right neighborhood. For isolated singularities, analytic
continuation verifies H\textup{(i)} outside a locally finite set of exponents
in the whole strict subcritical interval. Thus the preceding theorems
imply the following conclusions without additional non-degeneracy
assumptions.

\begin{corollary}\phantomsection\label{cor:gluing-applications}
\begin{enumerate}[label=(\roman*)]
\item\label{cor:stable-unconditional}\label{cor:near-threshold-unconditional}
Let $\Omega\subset\R^N$ be a bounded smooth domain, let $0\le k<N-2$,
and set $n=N-k$. If $k=0$, let $S\subset\partial\Omega$ be finite;
if $k>0$, let $S\subset\partial\Omega$ be a smooth compact embedded
$k$-dimensional submanifold. Suppose that
\[
\frac{n+1}{n-1}<p<p^*(n),
\]
or, when $p^*(n)<(n+2)/(n-2)$, that
\[
p^*(n)\le p<p^*(n)+\eta_*(n),
\qquad 0<\eta_*(n)<\frac{n+2}{n-2}-p^*(n),
\]
where $\eta_*(n)$ is supplied by Proposition~\ref{prop:near-threshold-H}.
Then the conclusion of Theorem~\ref{thm2} for $k=0$, or of
Theorem~\ref{thm3} for $k>0$, holds.

\item\label{thm:generic-isolated}
For each $N\ge3$, there is a locally finite set
\[
\mathcal E_N\subset\left(\frac{N+1}{N-1},\frac{N+2}{N-2}\right)
\]
such that, for every $p$ in this interval outside $\mathcal E_N$, every
bounded smooth domain $\Omega\subset\R^N$, and every finite set
$S\subset\partial\Omega$, the conclusion of Theorem~\ref{thm2} holds.
Here $\mathcal E_N$ may be taken to be the set in
Proposition~\ref{prop:generic-exponents} with $n=N$.
\end{enumerate}
\end{corollary}

The exceptional set in Corollary~\ref{cor:gluing-applications}\textup{(ii)}
is defined by the absence of any cell satisfying H\textup{(i)}, so it is
independent of the parametrization of a continuation branch.
For positive-dimensional singular sets, H\textup{(ii)} remains an assumption
outside the ranges covered by Corollary~\ref{cor:gluing-applications}\textup{(i)}; the generic argument for H\textup{(i)} does
not establish H\textup{(ii)}.

We now explain the main difficulties and the strategy of the proofs.
The interior fast-decaying cell is radial, so its construction and
linearized equation can be studied by ordinary differential equations.
The boundary cell has a nonconstant angular profile and requires a
different construction. We first construct the stable cells
by a Kelvin transform, a lifting that removes the degeneracy at the
flat boundary, and blow-up of bounded minimal solutions. A limit from
the stable side gives the threshold cell. Existence over the full
subcritical interval is then proved by analytic continuation.

More precisely, for this continuation we solve the equation near the
homogeneous profile on a half-cylinder and match it to an exterior
solution. The far-field coefficient is held fixed along the branch.
Uniform far-field remainders give compactness in the exterior space,
while the bounds on the matching data and positivity exclude limits on
the other boundary components of the admissible set. The two exponent
endpoint limits therefore make the branch cover the subcritical
interval. Each selected solution can then be normalized by dilation.
This proves existence as stated in Theorem~\ref{thm1}, without
requiring uniqueness throughout the interval.

Non-degeneracy requires a separate argument. In the stable range,
energy estimates exclude the relevant kernels. At the threshold,
the two indicial roots of the first angular mode coincide, and the repeated mode must
be compared with the scaling Jacobi field. A Wronskian identity and
the far-field normalization remove this remaining kernel and allow
continuation through the threshold. Fourier transformation in the
tangential variables reduces H\textup{(ii)} to shifted normal equations and
yields the required Liouville property in the stable range and in the
stated neighborhood of the threshold.
On the global branches, analytic Fredholm theory gives generic H\textup{(i)}.
We remark that the definition of $p^*(n)$ refers only to the initial
connected stable interval; no global sign assumption beyond it enters
this argument.

Finally, the gluing construction requires a right inverse whose bound
is uniform as the concentration scale tends to zero. A rescaling
argument has three possible limits: the homogeneous singular model,
the fast-decaying cell, and the harmonic equation. The corresponding
kernel estimates exclude these limits at the higher weight. At the
lower weight, where a kernel may remain, we choose the solution
orthogonal to it in a weighted Hilbert space and estimate the resulting
coupled system. This gives the uniform right inverse used in the
contraction argument. Near the flat boundary, both the cell and the
correction have a boundary factor. Keeping this factor in the
nonlinear estimates is essential, especially when $1<p<2$, and also
allows us to prove positivity and the very weak boundary identity.

The paper will be structured as follows. Theorem~\ref{thm1} is proved
in Sections~\ref{s2}--\ref{sec:analytic-continuation}.
Section~\ref{s2} constructs the stable cells and their threshold limit,
and proves their asymptotics and rotation invariance.
Section~\ref{sec:prelim} collects estimates for fast-decaying cells
which are uniform in the exponent, the spectral properties of the
linearized profile operator, and the analytic continuation theorem.
Section~\ref{sec:analytic-continuation} continues the cells analytically
through the full subcritical range, completing the proof in
Section~\ref{sec:proof-cell-existence}. Section~\ref{sec:model-linear}
proves the non-degeneracy results used in
Corollary~\ref{cor:gluing-applications}. After the tools of
Section~\ref{sec:massive-tools}, which reduce H\textup{(ii)} to shifted
normal operators, the stable and near-threshold
cases are treated in Sections~\ref{sec:stable-nondegeneracy}
and~\ref{sec:near-threshold-H}, and the generic case in
Section~\ref{sec:generic-nondegeneracy}.
The proofs of Theorems~\ref{thm2} and~\ref{thm3} occupy
Sections~\ref{s6}--\ref{s9}. Section~\ref{s6} defines the approximate
solution, Section~\ref{s8} constructs a uniformly bounded right inverse,
and Section~\ref{s9} solves the nonlinear equation. Positivity and the
very weak identity complete both proofs in
Section~\ref{sec:gluing-proof-completion}. Corollary~\ref{cor:gluing-applications}
is then proved in Section~\ref{sec:completion-main} by applying the
non-degeneracy results of Section~\ref{sec:model-linear}.
Appendix~\ref{app:continuation-estimates} contains the auxiliary estimates
for analytic continuation, and Appendix~\ref{app:fixed-regularity}
proves fixed-parameter regularity of the canonical inverse.

\section{Half-space cells in the stable regime}\label{s2}

Fix
\[
\frac{n+1}{n-1}<p<\frac{n+2}{n-2}.
\]
The homogeneous singular solution provides an upper bound for the
construction of a fast-decaying cell. In the stable regime, bounded
minimal solutions approach this singular profile, and a blow-up yields
the desired cell. We first identify the range in which stability makes
this construction possible. The threshold is then reached by compactness;
existence beyond it is obtained in Section~\ref{sec:analytic-continuation}.

\subsection{Stability of the homogeneous solution}\label{sec:stability}

Write $x=r\theta$, $\theta\in S^{n-1}_+$, and consider the homogeneous solution
\[
u_\infty(x)=r^{-\frac{2}{p-1}}\phi_p(\theta)
\]
of
\begin{equation}
\Delta u+u^p=0\quad\text{in }\R^n_+,
\qquad
u=0\quad\text{on }\partial\R^n_+\setminus\{0\}.
\end{equation}
The solution $u_\infty$ is called stable if
\[
Q_{u_\infty}(\phi):=
\int_{\R^n_+}|\nabla\phi|^2\,dx
-p\int_{\R^n_+}u_\infty^{p-1}\phi^2\,dx\ge0,
\qquad \phi\in H_0^1(\R^n_+).
\]

The stability of $u_\infty$ is determined by its angular linearization.
Consider the eigenvalue problem
\begin{equation}\label{e7}
\left\{ \begin{array}{l}
\Delta_{S^{n-1}}\psi_i-\beta_p \psi_i+p\phi_p^{p-1}\psi_i=-\lambda_i\psi_i \text{ in }S^{n-1}_+\\
\psi_i=0\text{ on }\partial S^{n-1}_+
\end{array}
\right.,
\end{equation}
where $\beta_p=\frac{p+1}{p-1}(n-\frac{p+1}{p-1})-(n-1)=\frac{2}{p-1}(n-2-\frac{2}{p-1})$.

Let $\lambda_1(p)\le\lambda_2(p)\le\cdots$ be its eigenvalues with $L^2(S^{n-1}_+)$-orthonormal eigenfunctions $\{\psi_i\}_{i\ge1}$.

\begin{proposition}\label{prop2}
Let $u_\infty$ be given by \eqref{e4}, and set
\begin{equation}\label{e8}
\Lambda_i(p):=\frac{(n-2)^2}{4}-\beta_p+\lambda_i(p).
\end{equation}
Then $u_\infty$ is stable if and only if
\begin{equation}\label{e9}
\Lambda_1(p)=\frac{(n-2)^2}{4}-\beta_p+\lambda_1(p)\ge0.
\end{equation}
\end{proposition}

\begin{proof}
For $x=r\theta$ and $t=\log r$, set
\[
w(t,\theta):=e^{\frac{n-2}{2}t}\phi(e^t\theta).
\]
For $\phi\in C_c^\infty(\overline{\R^n_+}\setminus\{0\})$ with zero Dirichlet trace,
\begin{align}
\int_{\R^n_+}|\nabla\phi|^2\,dx
&=\int_{\R\times S^{n-1}_+}
\left(|\partial_tw|^2+|\nabla_\theta w|^2+\frac{(n-2)^2}{4}w^2\right)
\,dt\,d\theta,\label{e10}\\
\int_{\R^n_+}u_\infty^{p-1}\phi^2\,dx
&=\int_{\R\times S^{n-1}_+}\phi_p^{p-1}w^2\,dt\,d\theta.
\label{e11}
\end{align}
Expanding $w=\sum_{i\ge1}a_i(t)\psi_i(\theta)$ yields
\begin{equation}\label{e12}
Q_{u_\infty}(\phi)
=\sum_{i\ge1}\int_\R
\left(|a_i'|^2+\Lambda_i(p) a_i^2\right)dt.
\end{equation}
Since $\Lambda_i\ge\Lambda_1$, this form is nonnegative if
$\Lambda_1\ge0$. Conversely, testing with
$w_R(t,\theta)=\eta(t/R)\psi_1(\theta)$, where
$0\ne\eta\in C_c^\infty(\R)$, gives
\[
Q_{u_\infty}(\phi_R)
=R^{-1}\int_\R|\eta'|^2dt
+\Lambda_1(p)R\int_\R\eta^2dt<0
\]
for large $R$ when $\Lambda_1<0$. This proves \eqref{e9}; density
extends the conclusion to the form domain.
\end{proof}

As $p\downarrow p_0:=\frac{n+1}{n-1}$, one has $\phi_p\to0$ and
$\beta_p\to-(n-1)$. The first Dirichlet eigenvalue of
$-\Delta_{S^{n-1}}$ is $n-1$, so \eqref{e7} implies $\lambda_1(p)\to0$ and
\[
 \Lambda_1(p)\longrightarrow\frac{(n-2)^2}{4}+n-1=\frac{n^2}{4}>0.
\]
Choose $\epsilon_0>0$ so that $\Lambda_1>0$ on $(p_0,p_0+\epsilon_0)$.
By uniqueness of $\phi_p$ and elliptic compactness, $\phi_p$ depends
continuously on $p$ on compact subintervals; the Rayleigh characterization
then ensures continuity of $\lambda_1(p)$. Define $p^*(n)$ as the right
endpoint of the connected component of
\[
 \{p\in(p_0,(n+2)/(n-2)):\Lambda_1(p)>0\}
\]
which contains $(p_0,p_0+\epsilon_0)$. This definition is independent of
the sufficiently small choice of $\epsilon_0$. Thus $p^*(n)>p_0$ and
$u_{\infty,p}$ is strictly stable on $(p_0,p^*(n))$. If
$p^*(n)<(n+2)/(n-2)$, continuity implies $\Lambda_1(p^*(n))=0$.

The same angular decomposition determines the behavior of the linearized
equation near the singular point. Writing
$\phi(r,\theta)=\sum_{i=1}^\infty \phi_i(r)\psi_i(\theta)$, we obtain
\begin{equation*}
\Delta_r \phi_i(r)+\frac{\beta_p-\lambda_i}{r^2}\phi_i(r)=0.
\end{equation*}
The indicial roots are the exponents $\gamma\in\C$ for which the
homogeneous linearized operator admits a nonzero separated solution
$r^\gamma\psi_i(\theta)$. Inserting this expression into the radial
equation yields
\begin{equation}\label{e15}
 \gamma^2+(n-2)\gamma+\beta_p-\lambda_i=0,
 \qquad
 \gamma_i^\pm=-\frac{n-2}{2}\pm\sqrt{\Lambda_i(p)}.
\end{equation}
We take the nonnegative square root when $\Lambda_i\ge0$ and the
square root with positive imaginary part otherwise. The indicial spectrum
and the set of real parts of its elements are, respectively,
\begin{equation}\label{eq:indicial-spectrum-definition}
 \mathfrak I_p:=\{\gamma_i^\pm:i\ge1\},\qquad
 \mathfrak R_p:=\{\operatorname{Re}\gamma:\gamma\in\mathfrak I_p\}.
\end{equation}
A real weight is nonindicial if it lies outside $\mathfrak R_p$; an
interval is root-free if it is disjoint from $\mathfrak R_p$. At a
repeated root, the two radial solutions are $r^\gamma\psi_i$ and
$r^\gamma\log r\,\psi_i$. We call terms $r^\gamma(\log r)^j\psi_i$
Frobenius terms; the integer $j\ge0$ records a logarithmic factor.
Thus statements about avoiding roots refer to
their real parts even when the roots are nonreal.
By the quadratic formula, in the stable range the two indicial roots in the first angular mode are real and satisfy $\gamma_1^-<-\frac{n-2}{2}<\gamma_1^+$. When the stability quantity in \eqref{e9} becomes negative, the corresponding first roots form a complex conjugate pair with real part $-\frac{n-2}{2}$.

When $\Lambda_1(p)\ge0$, we compare the first real roots as follows. Let $\mu_1=\lambda_1-\beta_p$ and
$\bar\mu_1$ be the principal eigenvalues of
$-\Delta_{S^{n-1}}-p\phi_p^{p-1}$ and
$-\Delta_{S^{n-1}}-\phi_p^{p-1}$, respectively, with Dirichlet data.
By the Rayleigh principle, $\mu_1<\bar\mu_1=-\beta_p$; the latter
eigenvalue has the positive eigenfunction $\phi_p$. The smaller root
of $\gamma^2+(n-2)\gamma-\mu=0$ decreases with $\mu$.
Since $2/(p-1)>(n-2)/2$, the smaller root for $\bar\mu_1$ is
$\bar\gamma_1^-=-2/(p-1)$. Thus
\begin{equation}\label{e16}
\gamma_1^+\geq -\frac{n-2}{2}\geq \gamma_1^->\bar{\gamma}_1^-=-\frac{2}{p-1}.
\end{equation}

\subsection[Construction in the stable range]{Construction in the stable range: $\frac{n+1}{n-1}<p<p^*(n)$}\label{sec:stable-construction}
We construct the cell by applying the Kelvin transform and dividing by the
flat-boundary factor $x_n$. In these variables, a positive solution
continuous at the origin will ensure the required fast decay at infinity.
Set
\[
\kappa_p:=p(n-2)-(n+2)\in\left(-\frac{2n}{n-1},0\right),
\qquad
\sigma_p:=n-1-\frac{2}{p-1}>0,
\]
so that $\sigma_p=[p(n-1)-(n+1)]/(p-1)$.
The homogeneous solution becomes
\[
\bar v_\infty(x):=
|x|^{\frac{2}{p-1}-(n-1)}
\frac{\phi_p(x/|x|)}{x_n/|x|},
\]
where the quotient is extended evenly across $\{x_n=0\}$.
Boundary Schauder estimates for the odd extension of $\phi_p$, followed
by division by its simple normal zero, produce a $C^{2,\alpha}$ quotient
away from the origin for $0<\alpha<\min\{1,p-1\}$; positivity up to
the equator follows from the Hopf lemma.

\smallskip
\noindent\emph{Bounded minimal solutions.}
We approximate $\bar v_\infty$ by bounded solutions with smaller boundary
values. For $\lambda\in(0,1)$, consider
\begin{equation*}
\left\{\begin{array}{l}
\operatorname{div}(x_n^2\nabla\bar v)+|x|^{\kappa_p}|x_n|^{p+1}\bar v^p=0 \text{ in }B_1\\
\bar v=\lambda \bar v_\infty \text{ on }\partial B_1\\
\bar v>0 \text{ in }B_1
\end{array}
\right.
.
\end{equation*}
Lift $\bar v$ to $V(x',y)=\bar v(x',|y|)$ on
$\R^m=\R^{n-1}\times\R^3$, $m=n+2$, with $z=(x',y)$. The equation becomes
\[
-\Delta V=QV^p,\qquad Q(z)=|z|^{\kappa_p}|y|^{p-1}.
\]
Since $d_p=p(n-1)-(n+3)>-2$, we have
$Q\in L^{r_0}_{\rm loc}$ for some $r_0>m/2$. Standard
elliptic estimates ensure local compactness of bounded solution families,
including at the origin. Let $V_\infty$ be the lifted homogeneous
solution and set $g=V_\infty|_{\partial B_1^m}$,
$M_g=\max_{\partial B_1^m}g$.
For a fixed $\lambda\in(0,1)$, set $t:=\lambda^{p-1}$. A bounded
supersolution is obtained by applying to $V_\infty$ the nondecreasing
concave function
\[
\Phi_t(s):=
\begin{cases}
s,&0\le s\le M_g,\\[1mm]
\bigl(t s^{1-p}+(1-t)M_g^{1-p}\bigr)^{\frac1{1-p}},&s>M_g.
\end{cases}
\]
Then $0\le\Phi_t(s)\le s$,
\[
\Phi_t'(s)s^p\ge t\Phi_t(s)^p,
\qquad
\sup_{s\ge0}\Phi_t(s)=M_g(1-t)^{-\frac1{p-1}}<\infty.
\]
The singularity of $V_\infty$ causes no extra source at the origin:
since $V_\infty\sim |z|^{-\sigma_p}$ and
\[
n-2\sigma_p=2-n+\frac4{p-1}>0,
\]
we have $V_\infty\in H^1(B_1^m)$ and $QV_\infty^p\in L^1(B_1^m)$.
The concave chain rule, first on punctured balls and then by smooth
approximation and removal of the inner cutoff, gives
\[
-\Delta\Phi_t(V_\infty)
\ge \Phi_t'(V_\infty)QV_\infty^p
\ge tQ\Phi_t(V_\infty)^p
\]
weakly on $B_1^m$.
Since $\Phi_t(V_\infty)=g$ on $\partial B_1^m$,
$\Phi_t(V_\infty)$ is a bounded supersolution of
\[
-\Delta W=tQW^p\quad\text{in }B_1^m,
\qquad W=g\quad\text{on }\partial B_1^m.
\]
Monotone iteration between the harmonic extension of $g$ and this
supersolution constructs a bounded minimal solution $W_t$. Standard Dirichlet
estimates imply $W_t\in W^{2,r_0}(B_1^m)$, and
\[
0<W_t\le\Phi_t(V_\infty)\le V_\infty.
\]
Multiplication by $\lambda$ identifies the $t$-problem with the lifted
problem having boundary value $\lambda g$, and therefore preserves
minimality. Returning to the original variables, we obtain the bounded minimal
solution $\bar v_\lambda$ from $V_\lambda:=\lambda W_t$; explicitly,
$\bar v_\lambda(x',t)=V_\lambda(x',(t,0,0))$, satisfying
\[
0<\bar v_\lambda\le\lambda\bar v_\infty<\bar v_\infty.
\]
This comparison also proves semistability. Indeed,
set $h:=V_\infty-V_\lambda>0$. By convexity, away from the origin,
\[
-\Delta h-pQV_\lambda^{p-1}h\ge0.
\]
For
$\xi\in C_c^\infty(B_1^m\setminus\{0\})$, integration by parts yields
\[
\begin{aligned}
\int_{B_1^m}\bigl(|\nabla\xi|^2-pQV_\lambda^{p-1}\xi^2\bigr)
&=\int_{B_1^m}h^2|\nabla(\xi/h)|^2\\
&\quad+\int_{B_1^m}
\frac{-\Delta h-pQV_\lambda^{p-1}h}{h}\,\xi^2\ge0.
\end{aligned}
\]
The cutoff at the origin has energy $O(\rho^{m-2})\to0$, while
$QV_\lambda^{p-1}\in L^{r_0}$ with $r_0>m/2$ makes the quadratic form
continuous on $H_0^1(B_1^m)$. Cutoff approximation and density therefore
extend this inequality to all $H_0^1$ test functions.

\smallskip
\noindent\emph{The limit as $\lambda\uparrow1$.}
Minimality and local elliptic estimates imply
$\bar v_\lambda\uparrow\bar v_*\le\bar v_\infty$ locally on
$\overline{B_1}\setminus\{0\}$, with boundary value $\bar v_\infty$.
For $w=\bar v_\infty-\bar v_*\ge0$, testing with cutoffs at the origin yields
\[
\begin{aligned}
\int_{B_1}x_n^2|\nabla w|^2\,dx
&=\int_{B_1}|x|^{\kappa_p}|x_n|^{p+1}
(\bar v_\infty^p-\bar v_*^p)w\,dx\\
&\le p\int_{B_1}|x|^{\kappa_p}|x_n|^{p+1}
\bar v_\infty^{p-1}w^2\,dx.
\end{aligned}
\]
Scale-invariant estimates for the lifted equation imply
$w=O(|x|^{-\sigma_p})$ and $|\nabla w|=O(|x|^{-\sigma_p-1})$.
Since $n-2\sigma_p>0$, the cutoff limit is valid and $w$ belongs to
the weighted form domain. Strict stability of $u_\infty$
(Section~\ref{sec:stability}), transferred by the Kelvin transform and the lifting,
forces $w=0$. Therefore
$\bar v_\lambda\to\bar v_\infty$ locally away from the origin and
$m_\lambda:=\|\bar v_\lambda\|_\infty\to\infty$ as $\lambda\uparrow1$.

\smallskip
\noindent\emph{A nontrivial entire limit.}
Choose $\lambda_j\to 1$, set
\begin{equation}
m_j=\|\bar v_{\lambda_j}\|_{L^\infty}, \ R_j=m_j^{\frac{p-1}{p(n-1)-(n+1)}},
\end{equation}
and define
\[
\bar v_j(x):=
m_j^{-1}\bar v_{\lambda_j}(\frac{x}{R_j}),
\qquad x\in B_{R_j}.
\]
Then $R_j\to\infty$, $0<\bar v_j\le\min\{1,C|x|^{-\sigma_p}\}$,
and $\bar v_j$ solves the equation in \eqref{e18} on $B_{R_j}$.
At a maximum point $x_j$, the bound just obtained implies
\[
1=\bar v_j(x_j)\le C|x_j|^{-\sigma_p},
\qquad |x_j|\le C^{1/\sigma_p}.
\]
Thus the maximum points stay in a fixed ball.
Since the lifted right-hand sides are bounded by
$Q\in L^{r_0}_{\rm loc}$, $r_0>m/2$, local elliptic compactness ensures,
along a subsequence, $x_j\to x_*$ and locally uniform convergence
$\bar v_j\to\bar v$, including at the origin. The limit solves
\begin{equation}
\label{e18}
\operatorname{div}(x_n^2\nabla \bar v)
+|x|^{\kappa_p} |x_n|^{p+1} \bar v^p=0
\quad\text{in }\R^n,
\end{equation}
and satisfies $\bar v(x_*)=1$. The strong maximum principle for the
lifted equation implies $\bar v>0$ and, in particular, $\bar v(0)>0$.
The limit is even in $x_n$ and $\bar v\le\bar v_\infty$.

\subsection{The half-space cell and its asymptotic behavior}

Let $\bar v$ be the positive even solution of \eqref{e18} just
constructed. The scaling $\bar v(x)\mapsto\rho^{\sigma_p}\bar v(\rho x)$
allows us to arrange $\bar v(0)=1$. Define
\[
v(x):=x_n\bar v(x),\qquad
u_1(x):=|x|^{-(n-2)}v\!\left(\frac{x}{|x|^2}\right),
\qquad x\in\R^n_+.
\]
Then $u_1$ solves \eqref{e3}. By continuity of $\bar v$ at the origin,
\[
u_1(x)=\bigl(1+o(1)\bigr)|x|^{-n}x_n
\quad\text{as }|x|\to\infty.
\]

It remains to determine the behavior at the singular end. The origin
cannot be removable, since otherwise the half-space Liouville theorem
for $1<p<(n+2)/(n-2)$ would exclude this positive solution with zero
Dirichlet boundary values; see \cite{bdn,gs}. As
$(n+1)/(n-1)<p<(n+2)/(n-2)$, the classification of isolated boundary
singularities \cite[Theorem~\textup{1.2}]{bpv1} therefore yields
\begin{equation}\label{eq:classification-C0}
\lim_{r\downarrow0}r^{\frac{2}{p-1}}u_1(r\theta)=\phi_p(\theta)
\quad\text{uniformly for }\theta\in\overline{S^{n-1}_+},
\end{equation}
where $\phi_p$ is the unique positive solution of \eqref{e5}.

For the relative asymptotic, rescale
$U_r(z)=r^{2/(p-1)}u_1(rz)$ on fixed half-annuli. The uniform
convergence in \eqref{eq:classification-C0}, the equation, and the
common zero trace ensure $C^1$ convergence by boundary $W^{2,q}$
estimates with $q>n$. The Hopf lemma therefore yields
\[
 |U_r(\theta)-\phi_p(\theta)|\le o(1)\theta_n,
 \qquad c_p\theta_n\le\phi_p(\theta)\le C_p\theta_n.
\]
Dividing by $\phi_p$ proves
\[
u_1(x)=|x|^{-\frac2{p-1}}\phi_p\!\left(\frac{x}{|x|}\right)(1+o(1))
\quad\text{as }|x|\to0,
\]
uniformly up to the flat boundary in the relative sense.

\smallskip
\noindent\emph{Passage to the threshold.}
Suppose $p^*(n)<(n+2)/(n-2)$ and let $p_j\uparrow p^*(n)$.
Before normalization at the origin, the lifted profiles satisfy
$-\Delta V_j=Q_{p_j}V_j^{p_j}$, $0<V_j\le1$, and attain $1$ in a
fixed ball. Since $d_{p_j}$ stays uniformly above $-2$, the coefficients
$Q_{p_j}V_j^{p_j-1}$ are bounded in $L^r$ on fixed balls for some
fixed $r>m/2$. Harnack's inequality consequently gives
\[
 0<c\le V_j(0)\le1,\qquad
 \rho_j:=V_j(0)^{-1/\sigma_{p_j}},\qquad
 \sigma_{p_j}=\frac{d_{p_j}+2}{p_j-1}.
\]
Both $\rho_j$ and $\rho_j^{-1}$ are bounded. The scaled profiles
$\rho_j^{\sigma_{p_j}}V_j(\rho_j z)$ solve the same equations and
equal $1$ at the origin. Local elliptic compactness and restriction to $y=(x_n,0,0)$
produce a positive threshold solution $\bar v_*$ with
$\bar v_*(0)=1$. Its Kelvin transform has normalized fast decay;
the Liouville and boundary classification arguments above establish its
singular asymptotic. Non-degeneracy at the repeated root is proved
separately in Section~\ref{sec:near-threshold-H}.

For both the stable cells and their threshold limit, the two end
asymptotics also imply the pointwise bound in Theorem~\ref{thm1}.  Near the
singular end, $\phi_p(\theta)\le C_p\theta_n$ by the Hopf lemma.
At infinity, the normalized Poisson-kernel asymptotic yields
$u_1(x)\le C_p|x|^{1-n}\theta_n$, and $1-n<-\frac2{p-1}$.  On a fixed intermediate annulus, boundary regularity
and the Hopf lemma bound $u_1/\theta_n$.  Combining the three regions yields
\[
u_1(x)\le c_p|x|^{-\frac2{p-1}}\frac{x_n}{|x|}.
\]

\subsection{Tangential rotation invariance}
The cells just constructed are rotation invariant by minimality and passage
to the limit. For later use we establish the same invariance for every
positive cell with the two prescribed end asymptotics, independently of
its construction.
\begin{lemma}\label{lem:cell-rotation}
Every positive solution of \eqref{e3} with the two end asymptotics in
\eqref{e6} satisfies
\begin{equation}\label{eq:cell-rotation-invariance}
 u_1(RX',X_n)=u_1(X',X_n),\qquad R\in O(n-1).
\end{equation}
\end{lemma}
\begin{proof}
We apply the moving-plane method in an arbitrary tangential direction. For $\lambda<0$, set
$\Sigma_\lambda=\{X_1<\lambda,X_n>0\}$ and let $u_\lambda$ be the
reflection of $u_1$. The function $v=(u_1-u_\lambda)_+$ vanishes near
the reflected singular point: there $u_\lambda/X_n$ tends to infinity
by the inner profile and the Hopf boundary factor, whereas $u_1/X_n$
is bounded. Its fast decay and boundary values imply
$v\in\mathcal D^{1,2}_0(\Sigma_\lambda)$. Testing the difference equation,
with cutoff approximation if needed, yields
\[
 \int|\nabla v|^2\le p\int_{\{v>0\}}u_1^{p-1}v^2
 \le Cp\|u_1^{p-1}\|_{L^{n/2}(\{v>0\})}\int|\nabla v|^2.
\]
Since $u_1=O(r^{1-n})$ and $(n-1)(p-1)>2$, the last coefficient is
small when the plane is sufficiently far to the left.

If the plane
stopped at $\lambda_0<0$, the strong maximum principle would give
strict comparison on compact interior sets; equality is excluded by
the reflected singular point. A small tail, the neighborhood where
$v=0$, and the small measure of the remaining boundary strip then
keep the same coefficient below one for nearby planes. This allows
the comparison to continue, contradicting the stopping assumption. The plane reaches zero; moving from
the opposite direction proves equality. Repeating in all tangential
directions proves \eqref{eq:cell-rotation-invariance}.
\end{proof}

\section{Uniform estimates and the linearized profile operator}\label{sec:prelim}
\label{sec:general-model-spaces}

In this section we collect estimates for fast-decaying cells which are
uniform in the exponent, the spectral properties of the linearized
profile operator, and the analytic continuation theorem used in
Section~\ref{sec:analytic-continuation}. We first introduce the
two-ended weighted spaces used throughout the paper. As before,
$n\ge3$ is the normal dimension.

For arbitrary weights $\mu,\nu\in\R$ and $q>n$, define
\begin{equation}
\label{e59-weight-twoend}
\omega_{\mu,\nu}(r):=
\begin{cases}
r^\mu, & 0<r\le1,\\
r^\nu, & r\ge1.
\end{cases}
\end{equation}

For $\ell\in\mathbb{Z}_{\ge 0}$ and $\mu,\nu\in\R$, set
    \begin{gather*}
        W^{\ell,q}_{\mu,\nu}(\R^n_+):=\left\{v:|v|_{\ell,q;\mu,\nu}<\infty\right\},\\
        |v|_{\ell,q;\mu,\nu}:=\sup_{j\in\mathbb{Z}}\omega_{\mu,\nu}(2^j)^{-1}\|v(2^j\,\cdot\,)\|_{W^{\ell,q}\left(B_2^+(0)\setminus\overline{B_\frac{1}{2}(0)}\right)},
        \\
        W^{\ell,q}_{\mu,\nu;0}(\R^n_+):=\left\{v\in W^{\ell,q}_{\mu,\nu}(\R^n_+):v|_{\partial\R^n_+\setminus\{0\}}=0\right\},\\
        L^q_{\mu,\nu}(\R^n_+):=W^{0,q}_{\mu,\nu}(\R^n_+).
    \end{gather*}
The traces are imposed when the Sobolev order admits them. We also use
the corresponding space of pointwise controlled functions; that is,
\begin{equation}\label{eq:two-ended-pointwise-space}
 C^0_{\mu,\nu}:=\left\{v\in C^0(\overline{\R^n_+}\setminus\{0\}):
 \|v\|_{C^0_{\mu,\nu}}
 :=\sup_{X\in\R^n_+}\frac{|v(X)|}{\omega_{\mu,\nu}(|X|)}<\infty\right\}.
\end{equation}
These supremum norms impose pointwise orders at the two ends; they do
not require square integrability with respect to a weighted measure.

For the half-space problem \eqref{e3}, write
$u_p^*:=u_\infty=|x|^{-2/(p-1)}\phi_p(x/|x|)$ for the homogeneous
solution~\eqref{e4} when its dependence on $p$ is relevant, and
retain $\beta_p=\frac2{p-1}(n-2-\frac2{p-1})$ from Section~\ref{s2}.
Before normalization, the positive far-field coefficient is denoted by
$c_u$. Rescaling $u(x)\mapsto\lambda^{2/(p-1)}u(\lambda x)$ changes it
to $c_u\lambda^{2/(p-1)+1-n}$, so it can be set to $1$ as in
\eqref{e6}.

\subsection{Uniform bounds and far-field asymptotics}
The continuation in Section~\ref{sec:analytic-continuation} requires
estimates which are uniform over families of fast-decaying cells.
We first control the rescaled profile, and then use the source moment
to obtain a uniform remainder at infinity.

For $p\in\left(\frac{n+1}{n-1},\frac{n+2}{n-2}\right)$, set
    \[
    w(\theta,t):=e^{\frac{2}{p-1}t}u(e^t\theta),\quad t\in\R,\quad \theta\in S^{n-1}_+.
    \]
Then problem~\eqref{e3} is transformed into
    \begin{equation}\label{e73}\tag{P'}
        \left\{\begin{aligned}
            \partial_t^2w-\alpha_p\partial_tw-\beta_pw+\Delta_{\theta}w+w^p&=0& &\text{in}& &Q,\\
            w&>0& &\text{in}& &Q,\\
            w&=0& &\text{on}& &\partial Q,
        \end{aligned}
        \right.
    \end{equation}
where $Q:=S^{n-1}_+\times\R$ and
\[
\alpha_p:=-n+2+\frac{4}{p-1}>0.
\]
Define the energy functionals
    \begin{align}
    E_p(w;t)&:=-\frac{1}{2}\int_{S^{n-1}_+}\left(\partial_tw\right)^2\,d\sigma+F_p(w(\,\cdot\,,t)),\label{e74}\\
    \begin{split}
    F_p(\phi)&:=\frac{\beta_p}{2}\int_{S^{n-1}_+}\phi^2\,d\sigma
    +\frac{1}{2}\int_{S^{n-1}_+}|\nabla_{S^{n-1}}\phi|^2\,d\sigma\\
    &\qquad-\frac{1}{p+1}\int_{S^{n-1}_+}\phi^{p+1}\,d\sigma,
    \end{split}\label{e75}
    \end{align}
where $w\in C^2(Q)\cap C^1\left(\overline{Q}\right)$ and $\phi\in C^2\left(S^{n-1}_+\right)\cap C^1\left(\overline{S^{n-1}_+}\right)$.
Every solution $w$ of~\eqref{e73} satisfies
    \begin{equation}\label{e76}
    \partial_tE_p(w;t)=-\alpha_p\int_{S^{n-1}_+}\left(\partial_t w\right)^2\,d\sigma\le0.
    \end{equation}
    Set
    \[
    A_p(w;t):=\int_{S^{n-1}_+}w\partial_tw\,d\sigma,
    \qquad
    H_p(w;t):=\int_{S^{n-1}_+}w^{p+1}\,d\sigma.
    \]
    Multiplying \eqref{e73} by $w$ and using the definition of $E_p$ yields
    \[
    \partial_tA_p
    =2\int_{S^{n-1}_+}(\partial_tw)^2\,d\sigma
      +\alpha_pA_p+2E_p
      -\frac{p-1}{p+1}H_p.
    \]
    Combining this identity with \eqref{e76} gives
    \begin{equation}\label{e77}
    \partial_t\left[
    e^{-\alpha_pt}\left(E_p(w;t)+\frac{\alpha_p}{2}A_p(w;t)\right)
    \right]
    =-\frac{\alpha_p(p-1)}{2(p+1)}
    e^{-\alpha_pt}H_p(w;t).
    \end{equation}
    In particular, if $u$ has a boundary singularity, then
    \cite[Theorem~\textup{1.2}]{bpv1} implies
    \begin{equation}\label{e78}
        \lim_{t\to-\infty}E_p(w;t)=E_p^*,
    \end{equation}
where
    \[
    E_p^*:=F_p(\phi_p)
    =\frac{p-1}{2(p+1)}
      \int_{S^{n-1}_+}\phi_p^{p+1}\,d\sigma>0.
    \]
For each fixed $p$, estimate~\eqref{e79} is contained in
\cite[Theorem~\textup{1.5}]{bpv1}. We also control its constant uniformly on
compact $p$-intervals and determine its behavior near the lower endpoint.
\begin{proposition}\label{prop6}
Let
\[
\frac{n+1}{n-1}<p<\frac{n+2}{n-2},
\]
and let $u\in C^2(\overline{\R^n_+}\setminus\{0\})$ be a fast-decaying boundary singular solution of~\eqref{e3}. Then
\begin{equation}\label{e79}
 |u(x)|\le C_{n,p}|x|^{-\frac{p+1}{p-1}}x_n,
 \qquad x\in\R^n_+.
\end{equation}
The constant $C_{n,p}$ is locally uniform for $p<\frac{n+2}{n-2}$ and, as $p\downarrow\frac{n+1}{n-1}$,
\[
 C_{n,p}
 =O\!\left(\|\phi_p\|_{L^{p+1}}^{\frac{p+1}{2p}}\right).
\]
\end{proposition}

\begin{proof}
Set
\[
\vartheta_p:=\frac{p-1}{2(p+1)},\qquad
G_p(t):=\alpha_p\vartheta_p e^{\alpha_pt}
\int_t^\infty e^{-\alpha_ps}H_p(w;s)\,ds.
\]
The fast-decay asymptotic, \eqref{e76}--\eqref{e78}, and the energy identity imply
\begin{equation}\label{e80}
G_p(t)=E_p(w;t)+\frac{\alpha_p}{2}A_p(w;t),
\end{equation}
\begin{equation}\label{e81}
\int_{\R\times S^{n-1}_+}(\partial_tw)^2\,d\sigma\,dt
=\alpha_p^{-1}E_p^*,
\end{equation}
and $G_p(t)\to E_p^*$ as $t\to-\infty$.  On
$I_t=[t,t+\alpha_p^{-1}]$ choose $t'$ so that
\[
\|\partial_tw(\cdot,t')\|_2^2\le CE_p^*,\qquad
H_p(w;t')\le CG_p(t).
\]
Then
$|A_p(w;t')|\le C(E_p^*)^{1/2}G_p(t)^{1/(p+1)}$, and \eqref{e80} yields
\begin{equation}\label{e82}
G_p(t+\alpha_p^{-1})
\le C\left(E_p^*+\alpha_p(E_p^*)^{1/2}G_p(t)^{1/(p+1)}\right).
\end{equation}
Young's inequality and iteration establish
\begin{equation}\label{e83}
G_p(t)
\le C_p\left(E_p^*+\alpha_p^{\frac{p+1}{p}}
(E_p^*)^{\frac{p+1}{2p}}\right),\qquad t\in\R.
\end{equation}
The constants are uniform on compact subintervals of the strict subcritical
range.  Thus $w$ is uniformly bounded in $L^{p+1}$ on translated unit
cylinders, and boundary elliptic estimates imply
$0<w(\theta,t)\le C_{n,p}\theta_n$, which is \eqref{e79}.  Near the lower
endpoint the local bootstrap is linear in this $L^{p+1}$ bound.  Since
\[
E_p^*=\frac{p-1}{2(p+1)}\|\phi_p\|_{L^{p+1}}^{p+1},
\qquad \alpha_p\to n,
\]
the stated estimate for $C_{n,p}$ follows.
\end{proof}
To retain the boundary factor in the far field, we estimate the Green
potential using its source moment and the source norm on a half-annulus.
\begin{lemma}\label{lem7}
Let $G$ be the Dirichlet Green kernel of $\R^n_+$, $n\ge3$, and put
\[
 M(f):=\int_{\R^n_+}y_n|f(y)|\,dy,\qquad
 A_r^+:=\{y\in\R^n_+:r/2<|y|<2r\}.
\]
Suppose $M(f)<\infty$ and let $v=G[f]$. For $r=|x|$ one has
\begin{equation}\label{e86}
 \frac{|v(x)|}{x_n}\le CM(f)r^{-n}
 +Cr^{1-(n+1)/Q}\|f\|_{L^Q(A_r^+;y_ndy)},\qquad Q>n+1,
\end{equation}
whenever the norm on the right is finite. The unweighted version is
\begin{equation}\label{e87}
 \frac{|v(x)|}{x_n}\le CM(f)r^{-n}
 +Cr^{1-n/q}\|f\|_{L^q(A_r^+)},\qquad q>n.
\end{equation}
Assume, in addition, that for $r\ge2$ either
\[
 \|f\|_{L^Q(A_r^+;y_ndy)}\le Ar^{\nu+(n+1)/Q}
 \quad\hbox{or}\quad
 \|f\|_{L^q(A_r^+)}\le Ar^{\nu+n/q},\qquad \nu<-n-1,
\]
with the corresponding exponent as above. Set
$\eta=-\nu-n-1>0$ and $\kappa=\frac12\min\{1,\eta\}$.
Then, uniformly in the angular variable,
\begin{equation}\label{e88}
\begin{gathered}
 \sup_{|x|=r,\ x_n>0}
 \left|\frac{r^nv(x)}{x_n}-c_f\right|
 \le C(M(f)+A)r^{-\kappa},\qquad r\ge4,\\
 c_f=\frac{2}{n|B_1|}\int y_nf(y)\,dy.
\end{gathered}
\end{equation}
Under the unweighted hypothesis, the remainder also satisfies
\begin{equation}\label{eq:green-sobolev-tail}
 \|v(r\,\cdot)-c_fr^{1-n}|\cdot|^{-n}(\cdot)_n\|_{W^{2,q}(A_1^+)}
 \le C(M(f)+A)r^{1-n-\kappa},\qquad r\ge8.
\end{equation}
The constants depend on the displayed exponents and the dimension, and
are uniform on bounded families of sources with these bounds.
\end{lemma}
\begin{proof}
\emph{The boundary-factor bound.}
The kernel and its two useful bounds are
\begin{equation}\label{e89}
\begin{gathered}
G(x,y)=\frac{|x-y|^{2-n}-|x-y^*|^{2-n}}{n(n-2)|B_1|},\\
G(x,y)\le C|x-y|^{2-n},\qquad
G(x,y)\le Cx_ny_n|x-y|^{-n}.
\end{gathered}
\end{equation}
Indeed, for $A=|x-y|^2$, integration in the reflected normal variable
shows that
\begin{equation}\label{eq:green-boundary-factor-integral}
 G(x,y)=\frac{2x_ny_n}{n|B_1|}
             \int_0^1(A+4t x_ny_n)^{-n/2}\,dt.
\end{equation}
The regions $|y|<r/2$ and $|y|>2r$ contribute at most
$CM(f)r^{-n}$ to $|v(x)|/x_n$. In $A_r^+$ set $d=x_n$.
On $|x-y|<d/2$ one has $y_n\asymp d$; the Newtonian bound and
H\"older's inequality give
\[
 d^{-1}\int_{|x-y|<d/2}G(x,y)|f(y)|\,dy
 \le C d^{1-(n+1)/Q}\|f\|_{L^Q(A_r^+;y_ndy)}.
\]
On the complement use the second bound in \eqref{e89}. On a shell
$|x-y|\asymp\ell\ge d/2$, the measure $y_ndy$ is at most
$C\ell^{n+1}$. The H\"older bound is therefore
$C\ell^{1-(n+1)/Q}\|f\|_{L^Q(A_r^+;y_ndy)}$.
Its dyadic sum up to $Cr$ is bounded by the second term of \eqref{e86},
since $Q>n+1$. Without the measure $y_ndy$, the near and far shells
instead contribute $d^{1-n/q}$ and $\ell^{1-n/q}$; their sum proves
\eqref{e87} for $q>n$.

\smallskip
\noindent\emph{The far-field coefficient and remainder.}
If $|y|\le T$ and $r\ge2T$, then
\[
 A+4t x_ny_n\ge r^2/4,\qquad
 |A+4t x_ny_n-r^2|\le3rT\quad(0\le t\le1).
\]
Indeed, $A\ge(r-T)^2\ge r^2/4$, and
\[
 A+4tx_ny_n-r^2=-2x'\cdot y'+(4t-2)x_ny_n+|y|^2.
\]
Since $|4t-2|\le2$, the Cauchy--Schwarz inequality bounds the first two
terms by $2rT$, while $|y|^2\le T^2\le rT/2$.
Applying the mean-value estimate to
\eqref{eq:green-boundary-factor-integral}, we obtain
\begin{equation}\label{e90}
 \left|\frac{r^nG(x,y)}{x_n}-\frac{2y_n}{n|B_1|}\right|
 \le C\frac{T}{r}y_n.
\end{equation}
By H\"older's inequality, either tail assumption implies
$\int_{A_s^+}y_n|f|\le CA s^{-\eta}$ for $s\ge2$. Summing the
annuli $s=2^jT$ yields
\[
 \int_{|y|>T}y_n|f|\le CA\sum_{j\ge0}(2^jT)^{-\eta}
                    \le CA T^{-\eta},\qquad T\ge2.
\]
The part $|y|>T$ outside $A_r^+$ is estimated by this tail moment,
and the middle annulus by \eqref{e86} or \eqref{e87}. Hence
\[
 \left|\frac{r^nv(x)}{x_n}-c_f\right|
 \le C\left(M(f)\frac{T}{r}+AT^{-\eta}+Ar^{-\eta}\right).
\]
Taking $T=\sqrt r$ proves \eqref{e88}. In the unweighted case the
remainder solves $-\Delta(v-c_f|x|^{-n}x_n)=f$ away from zero.
Scaled interior and Dirichlet estimates on slightly enlarged annuli,
using $r^{\nu+2}=r^{1-n-\eta}$, prove \eqref{eq:green-sobolev-tail}.
\end{proof}

\begin{lemma}\label{lem8}
Let $p\in((n+1)/(n-1),(n+2)/(n-2))$ and let $u$ be a positive
fast-decaying singular solution of \eqref{e3}, with coefficient $c_u>0$.
Then
\begin{equation}\label{e92}
 u(x)\le C|x|^{-n}x_n,\qquad |x|\ge1.
\end{equation}
On each compact $p$-interval, $C$ is uniform when $c_u$ is bounded above.
Moreover, for such a family, there are constants $C,\kappa>0$ such that
\begin{equation}\label{eq:cell-uniform-far-remainder}
 \sup_{|x|=r,\ x_n>0}
 \left|\frac{r^nu(x)}{x_n}-c_u\right|\le Cr^{-\kappa},\qquad r\ge4.
\end{equation}
The corresponding scaled $W^{2,q}$ remainder is uniform as well.
\end{lemma}
\begin{proof}
\emph{The source moment and Green representation.}
The individual end asymptotics imply $M(u^p)<\infty$: the radial
integrands have orders $r^{n-a_p-2}$ at zero and $r^{n-(n-1)p}$ at
infinity, where $a_p=2/(p-1)<n-1$ and $(n-1)p>n+1$.
Choose a smooth radial cutoff $\chi\ge0$, equal to one near zero and
zero for large arguments. Testing~\eqref{eq:halfspace-very-weak} with
$\zeta_R(x)=x_n\chi(|x|/R)$ gives
\[
\begin{gathered}
 \int x_nu^p\chi(|x|/R)\,dx=-\int u\Delta\zeta_R\,dx,\\
 \Delta\zeta_R=\frac{x_n}{R^2}
 \left[\chi''(r/R)+\frac{n+1}{r/R}\chi'(r/R)\right].
\end{gathered}
\]
Only a far-field annulus contributes to the right-hand side. The
individual normalized expansion there and radial integration imply
\[
\begin{aligned}
 \lim_{R\to\infty}-\int u\Delta\zeta_R\,dx
 &=-c_u\int_{S^{n-1}_+}\theta_n^2\,d\sigma
       \int_0^\infty\bigl(s\chi''+(n+1)\chi'\bigr)\,ds\\
 &=\frac{n|B_1|}{2}c_u.
\end{aligned}
\]
Consequently,
\begin{equation}\label{eq:cell-total-moment}
 M(u^p)=\int u^px_n\,dx=\frac{n|B_1|}{2}c_u.
\end{equation}
The same very weak identity, tested against a regularized Dirichlet
Green kernel with an outer cutoff, yields $u=G[u^p]$. The outer cutoff
terms vanish by the individual fast decay. No inner boundary term is
introduced because the test identity includes the origin.

\smallskip
\noindent\emph{Decay at infinity.}
Fix the auxiliary exponent $Q=n+2$, independently of the unweighted
exponent $q>n$ used elsewhere. The initial bound \eqref{e79} is
$u\le A r^{-a_p}\theta_n$. If at some stage
$u\le A r^{-b}\theta_n$ for $r\ge1$, $a_p\le b\le n-1$, then
\[
 \|u^p\|_{L^Q(A_r^+;x_ndx)}
 \le C A^{p(1-1/Q)}M(u^p)^{1/Q}r^{-bp(1-1/Q)},\qquad r\ge2.
\]
Equation~\eqref{e86} improves the exponent to
\[
 b_{\rm new}=\min\{n-1,\ p(1-1/Q)b+(n+1)/Q-2\}.
\]
Before the minimum reaches $n-1$, its untruncated increment is
\[
 [p(1-1/Q)-1](b-a_p)+\frac{n-1-a_p}{Q}
 \ge\frac{n-1-a_p}{Q}>0.
\]
The slope is positive for $Q=n+2$ and the present range of $p$.
On a compact $p$-interval the increment has a common positive lower
bound, so finitely many iterations, uniformly in $p$, establish \eqref{e92}.
The fixed annuli $1\le r\le2$ are controlled by \eqref{e79} at each step.
All constants depend only on the compact parameter interval and the bound
in \eqref{eq:cell-total-moment}.

Now $u^p\le Cr^{-(n-1)p}$ for $r\ge1$. If $p\ge\underline p>(n+1)/(n-1)$,
apply the unweighted part of Lemma~\ref{lem7} with
$\nu=-(n-1)\underline p<-n-1$. This proves
\eqref{eq:cell-uniform-far-remainder} and its Sobolev version.
\end{proof}
\subsection{Spectrum of the linearized profile equation}\label{sec:spectrum}
We use the spectral problem~\eqref{e7}, denote its eigenpairs by
$(\lambda_p^j,\psi_p^j)$, and take $\psi_p^1$ positive and
$L^2(S^{n-1}_+)$-normalized. Since $\phi_p$ is invariant under rotations of the first $n-1$
coordinates, decompose the angular operator by spherical-harmonic degree
in $L^2(S^{n-2})$. Write
\[
 \theta=(\cos s\,\omega,\sin s),\qquad
 0<s<\frac{\pi}{2},\qquad \omega\in S^{n-2}.
\]
For $m=0,1,2,\ldots$, let
$\{Y_{m,a}\}_{a=1}^{d_m}$ be an orthonormal basis of spherical harmonics
of degree $m$ on $S^{n-2}$, so that
\begin{equation}\label{e94}
 -\Delta_{S^{n-2}}Y_{m,a}=\nu_mY_{m,a},
 \qquad \nu_m=m(m+n-3).
\end{equation}
Choose $Y_{0,1}$ to be constant and, for $m=1$, choose the
basis proportional to the coordinate functions on $S^{n-2}$.  The
linearized angular operator preserves each degree-$m$ sector.  Thus its
eigenfunctions may be chosen in the separated form
\[
 \psi_p^{m,\ell,a}(\omega,s)
 =Y_{m,a}(\omega)\Psi_p^{m,\ell}(s),
 \qquad \ell=1,2,\ldots,
\]
where the associated eigenvalue, denoted by
$\lambda_p^{m,\ell}$, is independent of the multiplicity index $a$.  In
these coordinates the one-dimensional equation is
\begin{equation}\label{e94a}
 -\frac1{(\cos s)^{n-2}}
 \frac{d}{ds}\left((\cos s)^{n-2}\frac{d\Psi}{ds}\right)
 +\frac{\nu_m}{\cos^2s}\Psi
 +\beta_p\Psi-p\Phi_p^{p-1}\Psi
 =\lambda\Psi,
\end{equation}
where
$\Phi_p(s):=\phi_p(\cos s\,\omega,\sin s)$ is independent of $\omega$.
At $s=0$ one has the Dirichlet condition $\Psi(0)=0$; at
$s=\pi/2$ one imposes the regularity condition inherited from a smooth
function on the sphere.  This formulation keeps the
multiplicity of the degree-one spherical harmonics explicit.

By the compactness and uniqueness argument in Section~\ref{sec:stability},
$\phi_p$ depends continuously on $p$ in $C^2(\overline{S^{n-1}_+})$
on compact $p$-intervals. It follows from the min--max principle that
the eigenvalues in each fixed angular sector are continuous. The following spectral
gap determines the boundary data for the inner problem.

\begin{proposition}\label{prop7}
    Let $\{(\psi_p^j,\lambda_p^j)\}_{j=1}^\infty$ denote the eigenpairs from Section~\ref{sec:spectrum}. Then
    \begin{equation}\label{e95}
        \lambda_p^1<0,\qquad \lambda_p^j\ge\alpha_p+1>0,
        \quad j\in\{2,3,\ldots\}.
    \end{equation}
    Here
    $\lambda_p^1=\lambda_p^{0,1}$,
    $\lambda_p^{1,1}=\alpha_p+1$ with multiplicity $n-1$, while
    $\lambda_p^{0,\ell}>\alpha_p+1$ for $\ell\ge2$ and
    $\lambda_p^{m,1}>\alpha_p+1$ for $m\ge2$.
\end{proposition}
The proof is given in Appendix~\ref{app:proof-prop7}.

\begin{lemma}\label{lem:profile-analytic}
For
\[
\frac{n+1}{n-1}<p<\frac{n+2}{n-2},
\]
the map
\[
p\longmapsto\phi_p
\]
is real analytic locally in \(W^{2,q}(S^{n-1}_+)\), for every fixed
\(q>n\). For the angular operator
\[
\mathcal A_p
:=
-\Delta_{S^{n-1}}+\beta_p-p\phi_p^{p-1},
\]
the map $p\mapsto(\lambda_p^1,\psi_p^1)$ is locally real analytic with
values in
\[
 \R\times\bigl(W^{2,q}(S^{n-1}_+)\cap W^{1,q}_0(S^{n-1}_+)\bigr),
\]
where $\lambda_p^1$ is the lowest eigenvalue and $\psi_p^1$ is its
positive $L^2$-normalized eigenfunction.
\end{lemma}

The proof is given in Appendix~\ref{app:proof-lem:profile-analytic}.
Only the first eigenpair is selected analytically: it defines the
coefficient functional $c_p^1$ and the projection $P_p$ in
\eqref{e103}. The higher modes enter through the analytic inverse of
the boundary-value operator in Lemma~\ref{lem9}; no analytic choice of
their individual eigenfunctions is used.

\begin{remark}\label{rem:root-relation}
Set $a_p=2/(p-1)$.  If $\mu_j^\pm$ are the roots of
\begin{equation}\label{e102a}
\mu^2-\alpha_p\mu-\lambda_p^j=0,
\end{equation}
then the indicial roots in the original variables are
\begin{equation}\label{e102b}
\gamma_j^\pm=-a_p+\mu_j^\pm.
\end{equation}
The roots also satisfy
\begin{equation}\label{e102c}
\frac14(\alpha_p^2+4\lambda_p^j)
=\frac{(n-2)^2}{4}-\beta_p+\lambda_p^j.
\end{equation}
The numbers $\mu_j^\pm$ are the characteristic roots of the constant
coefficient modal equation in the cylindrical variable $t=\log r$.
Real roots are ordered by $\mu_j^-\le\mu_j^+$; otherwise they are
conjugate with real part $\alpha_p/2$. The corresponding solutions are
$e^{\mu_j^\pm t}$, or $e^{\mu_jt}$ and $te^{\mu_jt}$ at a repeated root.
For comparison, the Dirichlet harmonic model has indicial roots
\begin{equation}\label{eq:harmonic-indicial-roots}
 \ell,\quad 2-n-\ell,\qquad \ell=1,2,\ldots,
\end{equation}
because the angular eigenvalues are $\ell(\ell+n-2)$ on the odd
spherical harmonics. In particular $(1-n,1)$ is root-free. The
fast-decaying cell has this harmonic limiting operator at infinity.
\end{remark}

\subsection{An analytic global continuation theorem}
We use the following boundary-point version of analytic global
continuation, due to Garc\'ia--Haziot
\cite[Appendix~\textup{A}, Theorem~\textup{A}.2]{gh-bif} and based on
Chen--Walsh--Wheeler \cite[Theorem~\textup{6.1}]{cww-bif}. We state it here to
make the hypotheses and the three conclusions used below explicit.
\begin{theorem}[Analytic continuation from a boundary point]
\label{thm:boundary-continuation}
Let $X,Y$ be real Banach spaces, let $\mathcal O\subset\R\times X$ be
open with $(0,0)\in\partial\mathcal O$, and let
$\mathcal F:\mathcal O\to Y$ be real analytic. Assume:
\begin{enumerate}[label=(\Roman*)]
\item For some $\lambda_*>0$ there is a continuous local solution graph
\[
 \mathscr C_{\rm loc}
 =\{(\lambda,\widetilde x(\lambda)):0<\lambda<\lambda_*\}
 \subset\mathcal O\cap\mathcal F^{-1}(0),
 \qquad \widetilde x(\lambda)\longrightarrow0\quad(\lambda\downarrow0).
\]
\item $D_x\mathcal F(\lambda,\widetilde x(\lambda)):X\to Y$ is an
isomorphism for every $0<\lambda<\lambda_*$.
\item At every zero in $\mathcal O$, $D_x\mathcal F:X\to Y$ is
Fredholm of index zero.
\item There are bounded closed sets $Q_j\subset\mathcal O$ with
$\mathcal O=\bigcup_jQ_j$ such that
$Q_j\cap\mathcal F^{-1}(0)$ is compact for every $j$.
\end{enumerate}
Then $\mathscr C_{\rm loc}$ is contained in a continuous curve
\[
 \mathscr C=\{(\lambda(s),x(s)):s\in(0,\infty)\}
 \subset\mathcal O\cap\mathcal F^{-1}(0)
\]
with the following properties:
\begin{enumerate}[label=(\alph*)]
\item For every $j$ there is $s_j$ such that
$(\lambda(s),x(s))\notin Q_j$ whenever $s>s_j$.
\item Near each point, the curve admits a real-analytic
parametrization after a local change of its parameter.
\item $(\lambda(s),x(s))\notin\mathscr C_{\rm loc}$ for all sufficiently
large $s$.
\end{enumerate}
\end{theorem}
For its proof we refer to
\cite[Appendix~\textup{A}]{gh-bif} and \cite[Section~\textup{6}]{cww-bif}.

\section{Analytic continuation of half-space cells}\label{sec:analytic-continuation}

We match an inner solution near the homogeneous profile to an exterior
solution with the fixed far-field coefficient supplied by Lemma~\ref{lem14}.
The resulting analytic branches prove existence and are used in
Section~\ref{sec:generic-nondegeneracy}. Throughout this section and
Appendix~\ref{app:continuation-estimates}, $n\ge3$ is the normal dimension.

For the continuation construction we fix the far-field coefficient $c>0$
selected by Lemma~\ref{lem14}. All analytic continuation is performed
at this fixed coefficient; normalization to coefficient $1$ is made
separately at each parameter.

Fix $q>n$ and an exponent cutoff $\epstr$, independent of the gluing
scale $\ve$. For each sufficiently small
$\epstr>0$, set
\[
p_-:=\frac{n+1}{n-1}+\epstr,
\qquad
p_+:=\frac{n+2}{n-2}-\epstr,
\]
and fix an exterior source weight $\nu_0=\nu_0(\epstr)$ such that
\[
-(n-1)p_-<\nu_0<-n-1.
\]
Thus the continuation is carried out for $p\in(p_-,p_+)$.  The auxiliary
parameters are chosen in the order
\[
\begin{gathered}
\epstr\ \longrightarrow\ \nu_0\ \longrightarrow\ \eta_{{\rm an},0}
\ \longrightarrow\ \rho_0=c_*\eta_{{\rm an},0}
\ \longrightarrow\ \eta_{\rm an}\ \longrightarrow\ c,\\
0<\eta_{\rm an}\ll\rho_0\le\eta_{{\rm an},0}\ll1.
\end{gathered}
\]
Here $\eta_{{\rm an},0}$ bounds the deviation of the inner profile from
$\phi_p$ in $W^{2,q}_{\rm ul}$, and $\rho_0$ bounds its boundary data.
The smaller $C^2$ bound $\eta_{\rm an}$ in Lemma~\ref{lem14}, together
with the trace and first-mode derivative estimates, places every
coefficient-$c$ cell strictly inside these neighborhoods; the associated
admissible set $\mathcal V_{\eta_{{\rm an},0}}$ is defined in
\eqref{eq:matching-admissible-set}. The smallness conditions are
imposed in the continuation proof, with $c$ chosen last and
$\epstr,\nu_0$ kept fixed.

The exterior graph space $W$ is defined in Section~\ref{sec:continuation-outer};
$Q_-=S^{n-1}_+\times(-\infty,0)$, and $W^{2,q}_{\ul}$ denotes the
space with bounded $W^{2,q}$ norms on unit cylinder strips.

\begin{theorem}\label{thm4}
For every sufficiently small $\epstr>0$, there is a continuous curve
$\{(p_s,u_s)\}_{s>0}$ with
\[
p_s\in(p_-,p_+),\qquad
u_s\in W^{2,q}_{-\frac{2}{p_s-1},-n+1;0}(\R^n_+).
\]
The branch is continuous in fixed coordinates: the exterior restriction belongs to the fixed space $W$, while the inner part is represented by the Emden--Fowler profile
\[
w_s(\theta,t):=e^{\frac{2}{p_s-1}t}u_s(e^t\theta),\qquad t\le0,
\]
in the fixed space $W^{2,q}_{\rm ul}(\overline{Q_-})$.  The curve satisfies:
    \begin{enumerate}[label=(\roman*)]
\item \label{it1} $u_s$ is a strong solution of~\eqref{e3} with $p=p_s$ for every $s>0$.
\item \label{it2}
\[
u_s(x)=|x|^{-\frac{2}{p_s-1}}\bigl(\phi_{p_s}(x/|x|)+o(1)\bigr)
\qquad\text{as }|x|\to0.
\]
\item \label{it3}
\[
u_s(x)=(c+o(1))|x|^{-n}x_n
\qquad\text{as }|x|\to\infty.
\]
\item \label{it4} For every $s_0>0$, there exists an injective continuous map
\[
s^*:(-1,1)\to(0,\infty),\qquad s^*(0)=s_0,
\]
such that, in the fixed outer/inner coordinates just described,
$\sigma\mapsto p_{s^*(\sigma)}$, the corresponding exterior restriction,
and inner Emden--Fowler profile are real analytic.
\item \label{it5}
\[
\lim_{s\downarrow0}p_s=p_-,
\qquad
\lim_{s\to\infty}p_s=p_+.
\]

    \end{enumerate}
\end{theorem}

The proof proceeds by constructing the inner boundary map and then
continuing the exterior equation. Uniqueness of cells throughout the
full exponent interval is not asserted. Generic H\textup{(i)} is proved
separately in Section~\ref{sec:generic-nondegeneracy};
Lemma~\ref{lem:parameter-projection} then establishes local finiteness of the
exceptional exponents on each branch.
\subsection{Matching the inner and outer solutions}
We match the inner and outer solutions at $r=1$. In the cylinder,
both characteristic roots of the first angular mode of the linearized
equation have positive
real part, whereas every higher mode has one positive and one negative
root. Boundedness as $t\to-\infty$ therefore leaves two first-mode
data and one datum in each higher mode. We prescribe the value and
the first-mode derivative; the inner problem then determines the
remaining derivatives.

Let $\{\psi_p^j\}_{j=1}^\infty$ and $\{\lambda_p^j\}_{j=1}^\infty$ be as in Section~\ref{sec:spectrum}. Define the coefficient functionals $c_p^j$, $j\in\{1,2,\ldots\}$, and a projection operator $P_p$ by
\begin{gather}
c_p^j[f]:=\langle f,\psi_p^j\rangle,\quad f\in H^{-1}\left(S^{n-1}_+\right),\quad j\in\{1,2,\ldots\},\\
P_pf:=f-c_p^1[f]\psi_p^1,\quad f\in H^{-1}\left(S^{n-1}_+\right).\label{e103}
\end{gather}
\begin{enumerate}[label=(\roman*)]
    \item \label{it10} For small $z$ in the trace space $Z$ of Section~\ref{sec:inner-problem} and $b\in\R$, the inner problem determines $w=w[p,z,b]$:
    \begin{equation}\label{e104}
   \left\{ \begin{aligned}
        \partial_t^2w-\alpha_p\partial_tw-\beta_pw+\Delta_\theta w+w^p&=0& &\text{in}& &Q_-,\\
        w(\theta,t)=\phi_p(\theta)&+o\left(1\right)& &\text{as}& &t\to -\infty,\\
       w&=0& &\text{on}& &\partial_L Q_-,\\
        w(\,\cdot\,,0)=\phi_p(\theta)&+z& &\text{on}& & S^{n-1}_+,\\
        c^1_p[\partial_t w(\,\cdot\,,0)]&=b,
    \end{aligned}\right.
\end{equation}
where $Q_-:=S^{n-1}_+\times(-\infty,0)$, $\partial_L Q_-:=\partial S^{n-1}_+\times(-\infty,0)$. Then set
\[
B[p,z,b]:=P_p\partial_t w[p,z,b](\,\cdot\,,0).
\]
\item \label{it11} For small $\eta_{\rm an}>0$, solve the outer problem
\begin{equation}\label{e105}
\left\{\begin{aligned}
    -\Delta u&=u^p\quad\text{in}\quad  \R^n_+\setminus\overline{B_1(0)},\\
    P_p\partial_\lambda u|_{S^{n-1}_+}&=B\left[p,u|_{S^{n-1}_+}-\phi_p,c_p^1[\partial_\lambda u|_{S^{n-1}_+}]\right],\\
    u&=0\quad\text{on}\quad\partial\R^n_+\setminus\overline{B_1(0)},\\
    u(x)&=(c+o(1))|x|^{-n}x_n\quad\text{as}\quad|x|\to\infty,
\end{aligned}\right.
\end{equation}
where $\partial_\lambda u:=x\,\cdot\,\nabla u+\frac{2}{p-1}u$. Extend $u$ to $\R^n_+$ by
    \begin{equation}\label{e106}
\begin{gathered}
 z=u|_{S^{n-1}_+}-\phi_p,\qquad b=c_p^1[\partial_\lambda u|_{S^{n-1}_+}],\\
 u(e^t\theta):=e^{-\frac2{p-1}t}w[p,z,b](\theta,t),
 \qquad t<0,\quad\theta\in S^{n-1}_+.
\end{gathered}
    \end{equation}
At $r=1$, the reconstruction has
$\partial_ru=w_t-\frac2{p-1}w$ and $\partial_\lambda u=w_t$.
Thus the common value, the prescribed first-mode derivative and the
projected matching equation determine the full physical normal derivative.
The resulting $u$ is a fast-decaying boundary singular solution of~\eqref{e3}.
\end{enumerate}
To apply this local inner construction to the entire continuation
branch, we choose the far-field coefficient so that the inner profiles
remain close to $\phi_p$. The next lemma provides this choice.
\begin{lemma}\label{lem14}
Let $\epstr>0$ and $\eta_{\rm an}>0$. Then there exists $c>0$ such that, if $p\in\left(\frac{n+1}{n-1},\frac{n+2}{n-2}-\epstr\right)$ and $w\in C^2(Q)\cap C\left(\overline{Q}\right)$ solves~\eqref{e73} with
\[
\lim_{t\to-\infty}w(\theta,t)=\phi_p(\theta),\quad \lim_{t\to+\infty}e^{\left(n-1-\frac{2}{p-1}\right)t}w(\theta,t)=c'\theta_n,\quad \theta\in S^{n-1}_+,
\]
for some $c'\ge c$, then
\begin{equation}\label{e108}
    \|w(\,\cdot\,,t)-\phi_p\|_{C^2\left(\overline{S^{n-1}_+}\right)}<\eta_{\rm an},\quad t<0.
\end{equation}
\end{lemma}
\begin{proof}
    By \eqref{e79}, there exists $\epstr'>0$ such that if $p\le\frac{n+1}{n-1}+\epstr'$, then
    \[
    \|w(\,\cdot\,,t)\|_{C^2\left(\overline{S^{n-1}_+}\right)}<\frac{\eta_{\rm an}}{2},\quad \|\phi_p\|_{C^2\left(\overline{S^{n-1}_+}\right)}<\frac{\eta_{\rm an}}{2},\quad t\in\R.
    \]
    Then \eqref{e108} holds for all $t\in\R$. Assume also that $p>\frac{n+1}{n-1}+\epstr'$.

\proofstep{1}{Construct two translations of the counterexample sequence.}
    If the assertion fails, there are $c_i>1$ tending to infinity,
    $p_i\in((n+1)/(n-1)+\epstr',(n+2)/(n-2)-\epstr)$, and
    $w_i\in C^2(Q)\cap C(\overline Q)$ solving~\eqref{e73} with $p=p_i$,
    such that, with $a_{p_i}=2/(p_i-1)$,
    \[
    \lim_{t\to-\infty}w_i(\theta,t)=\phi_{p_i}(\theta),\qquad
    \lim_{t\to+\infty}e^{(n-1-a_{p_i})t}w_i(\theta,t)=c_i\theta_n,
    \]
    but $\|w_i(\cdot,\tau_i)-\phi_{p_i}\|_{C^2(\overline{S^{n-1}_+})}
    \ge\eta_{\rm an}$ for some $\tau_i<0$. Translating by $\tau_i$
    preserves the equation and the singular limit, and changes the
    far-field coefficient to
    \[
    \widehat c_i=c_i e^{-(n-1-a_{p_i})\tau_i}\ge c_i.
    \]
    Relabeling the translated functions and coefficients as $w_i,c_i$,
    we retain both limits above and $c_i\to\infty$, with
    \[
    \|w_i(\cdot,0)-\phi_{p_i}\|_{C^2(\overline{S^{n-1}_+})}
    \ge\eta_{\rm an}.
    \]
    Set
    \[
    t_i:=\frac{\log c_i}{n-1-a_{p_i}},\qquad
    \widetilde w_i(\theta,t):=w_i(\theta,t+t_i).
    \]
    The uniform bound \eqref{e79} and standard boundary elliptic
    estimates give, after passing to a subsequence,
    \[
    p_i\to p_\infty\in
    \left[\frac{n+1}{n-1}+\epstr',\frac{n+2}{n-2}-\epstr\right],
    \qquad w_i\to w_*,\quad\widetilde w_i\to w^*
    \]
    in $C^2$ on every compact closed strip of $\overline Q$.
    This second translation normalizes the far-field coefficient to $1$.
    Since $p_i>\frac{n+1}{n-1}+\epstr'$,
    \[
    t_i\ge\frac{1}{n-1}\log c_i\to\infty
    \quad\text{as}\quad i\to\infty.
    \]
    \proofstep{2}{Establish a uniform energy gap.}
    There exists $\gamma>0$ such that
    \begin{equation}\label{e109}
        E_{p_i}(w_i;0)\le E^*_{p_i}-\gamma,\qquad i=1,2,\ldots .
    \end{equation}
    Otherwise, after passing to a subsequence (with the same notation),
    \[
    E_{p_i}(w_i;0)\to E^*_{p_\infty}\quad\text{as}\quad i\to\infty.
    \]
    The limit $w_*$ solves \eqref{e73} at $p=p_\infty$. Since
    \[
    E_{p_i}(w_i;0)\le E_{p_i}(w_i;t)\le E^*_{p_i},\qquad t\le0,
    \]
    its energy is identically $E^*_{p_\infty}>0$ on $t\le0$.
    Equation~\eqref{e76} makes $w_*$ stationary, and
    \cite[Theorem~\textup{1.1}]{bpv1} identifies it with $\phi_{p_\infty}$.
    This contradicts the $C^2$ separation at $t=0$ and proves
    \eqref{e109}.

    \par\medskip
    By \eqref{e76} and \eqref{e109},
    $E_{p_i}(\widetilde w_i;t)=E_{p_i}(w_i;t+t_i)
    \le E^*_{p_i}-\gamma$ for $t\ge-t_i$. Since $t_i\to\infty$,
    passage to the limit implies
    \[
    E_{p_\infty}(w^*;t)\le E^*_{p_\infty}-\gamma,\quad t\in\R.
    \]
    Let $U^*$ be the corresponding solution in $\R^n_+$.
    This energy bound forces $w^*\equiv0$. Indeed, if the origin were a nonremovable boundary
    singularity of $U^*$, then \cite[Theorem~\textup{1.2}]{bpv1} would imply
    $w^*(\,\cdot\,,t)\to\phi_{p_\infty}$ as $t\to-\infty$, and hence
    $E_{p_\infty}(w^*;t)\to E^*_{p_\infty}$, contradicting the strict
    bound above.  Thus the origin is removable.  The half-space Liouville
    theorem \cite{bdn,gs} then excludes a nonzero positive classical
    solution in $\R^n_+$ with zero Dirichlet boundary data.  Hence
    $w^*\equiv0$ in $Q$.

\proofstep{3}{Obtain a positive lower bound at a fixed point.}
The corresponding functions
    $\widetilde u_i(r\theta)=r^{-a_{p_i}}\widetilde w_i(\theta,\log r)$
    all have far-field coefficient $1$. Lemma~\ref{lem8} provides constants
    $C,\kappa>0$, independent of $i$, such that
    \[
      \left|R^{n-1}\widetilde u_i(Re_n)-1\right|\le CR^{-\kappa}.
    \]
    Fix $R$ so large that $CR^{-\kappa}\le1/2$. Then
    \[
    \widetilde u_i(Re_n)\ge R^{1-n}/2>0,
    \]
    contradicting $\widetilde w_i\to0$ on compact strips.
\end{proof}

\subsection{The inner problem}\label{sec:inner-problem}
We solve the inner problem to construct the analytic map $B[p,z,b]$.
To apply the implicit function theorem, first consider its linearization
    \begin{equation}\label{e112}
        \left\{
        \begin{aligned}
            \partial_{t}^2\omega-\alpha_p\partial_t\omega-\beta_p\omega+\Delta_\theta\omega+p\phi_p^{p-1}\omega&=f& &\text{in}& &Q_-,\\
            \omega&=0& &\text{on}& &\partial_L Q_-,\\
            \omega(\,\cdot\,,0)&=z& &\text{on}& & S^{n-1}_+,\\
            c_p^1[\partial_t\omega(\,\cdot\,,0)]&=b.
        \end{aligned}
        \right.
    \end{equation}
For $\ell\in\mathbb Z_{\ge0}$, set
\[
\|w\|_{\ell,q}
:=\sup_{t\le0}
\|w\|_{W^{\ell,q}(S^{n-1}_+\times(t-1,t))},
\qquad
W^{\ell,q}_{\ul}(\overline{Q_-})
:=\{w:\|w\|_{\ell,q}<\infty\},
\]
and $L^q_{\ul}(\overline{Q_-}):=W^{0,q}_{\ul}(\overline{Q_-})$.
The closed Banach subspace for the nonlinear unknown is
\[
 E_-:=\{v\in W^{2,q}_{\rm ul}(\overline{Q_-}):
                                  v|_{\partial_LQ_-}=0\},\qquad
 \|v\|_{E_-}:=\|v\|_{2,q}.
\]
For $\ell\in\mathbb Z_{\ge1}$, with $\alpha=1-\frac nq$,
$W^{\ell,q}_{\ul}(Q_-)\subset BC^{\ell-1,\alpha}(\overline{Q_-})$. Set
\begin{align*}
    Z&:=W^{2-q^{-1},q}\left(S^{n-1}_+\right)\cap W^{1,1}_0\left(S^{n-1}_+\right)\\
    &=\left\{w(\,\cdot\,,0):w\in W^{2,q}_\ul\left(\overline{Q_-}\right),w|_{\partial_LQ_-}=0\right\},\\
    Z_+&:=\left\{z\in Z:\inf_{\theta\in S^{n-1}_+}\theta_n^{-1}z(\theta)>0\right\},\\
    \|z\|_Z&:=\inf\left\{\|w\|_{2,q}:w\in W^{2,q}_\ul\left(\overline{Q_-}\right),\ w|_{\partial_LQ_-}=0,\ w(\,\cdot\,,0)=z\right\}.
\end{align*}
\begin{lemma}\label{lem9}
Let $p\in(p_-,p_+)$. For every $f\in L^q_{\rm ul}(\overline{Q_-})$, $z\in Z$, and $b\in\R$, problem~\eqref{e112} admits a unique solution
\[
\omega=\omega[p,z,f,b]\in W^{2,q}_{\rm ul}(\overline{Q_-}),
\]
and
\begin{equation}\label{e113}
\|\omega\|_{2,q}
\le C_{n,\epstr,q}\bigl(\|f\|_{0,q}+\|z\|_Z+|b|\bigr),
\end{equation}
where $C_{n,\epstr,q}$ is independent of $p\in(p_-,p_+)$.
\end{lemma}
The proof is given in Appendix~\ref{app:proof-lem9}.

We now establish nonlinear solvability and analytic dependence for the inner problem.
\begin{proposition}\label{prop8}
Fix the truncated $p$-interval and put
\[
 \sigma_0:=\frac12\inf_{p\in(p_-,p_+)}\operatorname{Re}\mu_1^-(p)>0.
\]
There is $c_*>0$, depending on the linear estimates of
Lemma~\ref{lem9}, their weighted version at $\sigma_0$, and the
embedding constant in Remark~\ref{rem5}, such that for each sufficiently small solution radius
$\eta_{\rm an}>0$, with $\rho_{\eta_{\rm an}}=c_*\eta_{\rm an}$, one may take
\[
 \mathcal U_{\eta_{\rm an}}=\{(p,z):p\in(p_-,p_+),\ \|z\|_Z<\rho_{\eta_{\rm an}}\}.
\]
There is a real-analytic map
\[
w:\mathcal U_{\eta_{\rm an}}\times(-\rho_{\eta_{\rm an}},\rho_{\eta_{\rm an}})
\longrightarrow W^{2,q}_{\rm ul}(\overline{Q_-})
\]
such that $\phi_p+z\in Z_+$ for $(p,z)\in\mathcal U_{\eta_{\rm an}}$, and $w[p,z,b]$ is the unique solution of~\eqref{e104} satisfying
\[
\|w[p,z,b]-\phi_p\|_{2,q}\le\eta_{\rm an}.
\]
If
\[
\pi[p,z,b]:=w[p,z,b]-\phi_p-\omega[p,z,0,b],
\]
then
\begin{equation}\label{e117}
\|\partial_z\pi[p,z,b]\|_{Z\to W^{2,q}_{\rm ul}}
+\|\partial_b\pi[p,z,b]\|_{2,q}
\le C\eta_{\rm an}
\end{equation}
for $(p,z,b)\in\mathcal U_{\eta_{\rm an}}\times(-\rho_{\eta_{\rm an}},\rho_{\eta_{\rm an}})$.
\end{proposition}
\begin{proof}
Work on the fixed zero-lateral-trace space $E_-$ and write
$w=\phi_p+v$. Put
$\mathcal N_p(v)=(\phi_p+v)^p-\phi_p^p-p\phi_p^{p-1}v$. Then
\[
 \mathcal F(p,z,b,v):=v-\omega[p,z,-\mathcal N_p(v),b]=0,
 \qquad D_v\mathcal F(p,0,0,0)=I_{E_-}.
\]
For $v$ with zero lateral trace, boundary division gives
$\|v/\phi_p\|_\infty\le C\|v\|_{2,q}$. Put $h=v/\phi_p$.
For $\|h\|_\infty\le1/2$, the exact factorization and Taylor estimate are
\[
 \mathcal N_p(v)=\phi_p^p\bigl((1+h)^p-1-ph\bigr),\qquad
 |\mathcal N_p(v)|\le C\phi_p^p|h|^2.
\]
Consequently,
\[
 \|\mathcal N_p(v)\|_{0,q}\le C\|v\|_{2,q}^2,\qquad
 \|D\mathcal N_p(v)\|_{E_-\to L^q_{\rm ul}}
       \le C\|v\|_{2,q}.
\]
The quotient $(\phi_p+v)/\theta_n$ is then uniformly positive. Lemmas~\ref{lem:profile-analytic} and~\ref{lem9}
make $\mathcal F$ analytic, including for noninteger exponents, so the analytic
implicit function theorem applies.

On the truncated interval set
$R_0f=\omega[p,0,f,0]$. With $C_{\rm lin}$ the constant in \eqref{e113}, choose $2C_{\rm lin}c_*\le1/2$, then
take $\eta_{\rm an}$ so small that
$\|R_0D\mathcal N_p(v)\|\le1/2$ on its ball. The contraction theorem
ensures uniqueness throughout the ball and
$\|v\|_{2,q}\le C(\|z\|_Z+|b|)$ for
$\|z\|_Z,|b|<c_*\eta_{\rm an}$. Finally,
\[
 \pi=-R_0\mathcal N_p(v),\qquad
 \|D\mathcal N_p(v)\|_{E_-\to L^q_{\rm ul}}
 \le C\eta_{\rm an}.
\]
Differentiating the contraction equation bounds
$\partial_zv,\partial_bv$ uniformly and proves \eqref{e117}.
The weighted construction in Remark~\ref{rem5} below proves the
singular-end limit required in~\eqref{e104} for this bounded solution.
\end{proof}
For the later matching argument we also need the inverse at the
inner solution $w=w[p,z,b]$. Write
\[
\begin{gathered}
 D_0=\partial_t^2-\alpha_p\partial_t-\beta_p+\Delta_\theta+p\phi_p^{p-1},\\
 V_w=p(w^{p-1}-\phi_p^{p-1}),\qquad D_w=D_0+V_w,
 \qquad R_0f=\omega[p,0,f,0].
\end{gathered}
\]
The zero lateral trace and Hopf bound imply
$|w-\phi_p|\le C\eta_{\rm an}\theta_n$ and $\phi_p\ge c\theta_n$.
Writing $w=\phi_p(1+h)$, we have $\|h\|_\infty\le C\eta_{\rm an}$ and
\[
 |V_w|=p\phi_p^{p-1}|(1+h)^{p-1}-1|
 \le C\eta_{\rm an}\phi_p^{p-1}\le C\eta_{\rm an}.
\]
This bound also holds for $1<p<2$, despite the failure of the scalar
map $s\mapsto s^{p-1}$ to be Lipschitz at zero.
Shrink the radius uniformly so that $\|R_0V_w\|<1/2$. The Neumann
series yields
\begin{equation}\label{eq:actual-inner-source}
 R_w=(I+R_0V_w)^{-1}R_0,
 \qquad \|R_wf\|_{2,q}\le C\|f\|_{0,q}.
\end{equation}
The identity $\eta=R_0(f-V_w\eta)$ shows that $R_w$ is the inverse with boundary conditions
\[
 D_w\eta=f,\quad \eta|_{\partial_LQ_-}=0,\quad
 \eta(\cdot,0)=0,\quad c_p^1[\partial_t\eta(\cdot,0)]=0.
\]
Lemma~\ref{lem9} establishes the converse and hence uniqueness in
$W^{2,q}_{\rm ul}$, without imposing decay at $t=-\infty$ for a general
$L^q_{\rm ul}$ source.

We can now define the boundary map and its linearization. For any function space $X$ continuously embedded in $H^{-1}(S^{n-1}_+)$, let $P_pX$ be the subspace of functions $\varphi\in X$ such that $c_p^1[\varphi]=0$, equipped with the same norm as $X$. Set
    \begin{align*}
    Z'&:=W^{1-q^{-1},q}_0\left(S^{n-1}_+\right)\\
    &=\left\{\partial_tw(\,\cdot\,,0):w\in W^{2,q}_\ul\left(\overline{Q_-}\right),w|_{\partial_LQ_-}=0\right\},\\
    \|\varphi\|_{Z'}&:=\inf\left\{\|w\|_{2,q}:w\in W^{2,q}_\ul\left(\overline{Q_-}\right),\ w|_{\partial_LQ_-}=0,\ \partial_tw(\,\cdot\,,0)=\varphi\right\}.
\end{align*}
We then define an analytic map $B:\mathcal{U}_{\eta_{\rm an}}\times(-\rho_{\eta_{\rm an}},\rho_{\eta_{\rm an}})\to Z'$ by
    \[
    B[p,z,b]:=P_p\partial_tw[p,z,b](\,\cdot\,,0)\in P_pZ',
    \]
    where $P_p$ is as in \eqref{e103}. Define the linearized operator $C_p:Z\to P_pZ'$ by
    \[
    C_pz:=P_p\partial_t\omega[p,z,0,0](\cdot,0).
    \]
Since $P_p\omega[p,0,0,b]=0$,
\begin{equation}\label{e118}
\begin{gathered}
P_p\partial_t\omega[p,z,0,b](\cdot,0)=C_pz,\\
\partial_zB[p,z,b]h=C_ph+P_p(\partial_z\partial_t\pi[p,z,b]h)(\cdot,0).
\end{gathered}
\end{equation}
\begin{remark}\label{rem4}
For $h\in P_pZ'$, choose, by the definition of the trace norm, a function
$w_h\in W^{2,q}_{\ul}(\overline{Q_-})$ with zero lateral trace,
$\partial_tw_h(0)=h$, and $\|w_h\|_{2,q}\le2\|h\|_{Z'}$; take $w_0=0$.
Let $\chi=1$ near $0$, with compact support in $(-\infty,0]$, and set
\[
 H_h(\theta,t)=\chi(t)\bigl(w_h(\theta,t)-w_h(\theta,2t)\bigr),
 \qquad
 f_h=(\partial_t^2-\alpha_p\partial_t-\beta_p+
                    \Delta_\theta+p\phi_p^{p-1})H_h.
\]
Then $H_h(0)=0$, $\partial_tH_h(0)=-h$, and
$c_p^1[\partial_tH_h(0)]=0$, while
$\|f_h\|_{0,q}\le C\|h\|_{Z'}$.
Uniqueness in Lemma~\ref{lem9} identifies
$H_h=\omega[p,0,f_h,0]$. Thus
\[
 \partial_t\omega[p,0,f_h,0](0)=-h,
 \qquad \|f_h\|_{0,q}\le C\|h\|_{Z'}.
\]
The first-mode condition here concerns only the terminal derivative.
\end{remark}

\begin{remark}\label{rem5}
The exponential weight defines the spaces
\[
\begin{aligned}
 L^q_{{\rm ul},\sigma}(Q_-)
   &:=\{f:e^{-\sigma t}f\in L^q_{\rm ul}(Q_-)\},\\
 E_{-,\sigma}&:=\{v:e^{-\sigma t}v\in E_-\},
 \qquad \|v\|_{E_{-,\sigma}}:=\|e^{-\sigma t}v\|_{2,q}.
\end{aligned}
\]
The source norm is $\|f\|_{L^q_{{\rm ul},\sigma}}
:=\|e^{-\sigma t}f\|_{0,q}$. Suppose
$f\in L^q_{{\rm ul},\mu}(Q_-)$ for some $\mu>0$.
Lemma~\ref{lem9} gives, for
\[
0<\widetilde\mu<\min\{\mu,\operatorname{Re}\mu_1^-\},
\]
\[
\|e^{-\widetilde\mu t}\omega[p,z,f,b]\|_{2,q}
\le C_{\widetilde\mu}
\bigl(\|e^{-\mu t}f\|_{0,q}+\|z\|_Z+|b|\bigr).
\]
Set $\sigma=\widetilde\mu$ and solve the equation conjugated by
$e^{\sigma t}$. Its roots are $\mu_j^\pm-\sigma$, and its terminal
first-mode derivative is $b-\sigma z_1$, where $z_1=c_p^1[z]$.
At a repeated root $r=\alpha_p/2$, the source kernel satisfies
\[
\int_0^\infty h e^{-(r-\sigma)h}\,dh=(r-\sigma)^{-2}<\infty.
\]
The high-mode spectral gap, spectral summation, and strip
regularity estimates in Appendix~\ref{app:proof-lem9} provide the same
$W^{2,q}_{\rm ul}$ control for the conjugated equation. Multiplying
this solution by $e^{\sigma t}$ and using the unweighted uniqueness
in Lemma~\ref{lem9} identifies it with $\omega[p,z,f,b]$. Thus this
solution belongs to $E_{-,\sigma}$ and satisfies the estimate. On compact
parameter intervals the constants are uniform provided
$\inf_p(\operatorname{Re}\mu_1^--\sigma)>0$.
If $\mu_1^-<\mu_1^+$ are real and $\mu>\mu_1^-$, $\mu\ne\mu_1^+$, then for some $a_0\in\R$,
\[
\left\|e^{-\mu''t}
\bigl(\omega[p,z,f,b]-a_0e^{\mu_1^-t}\psi_p^1\bigr)\right\|_{2,q}
\le C\bigl(\|e^{-\mu t}f\|_{0,q}+\|z\|_Z+|b|\bigr),
\]
where $\mu''=\min\{\mu,\mu_1^+\}$.

For the nonlinear inner solution $v=w[p,z,b]-\phi_p$ the equation is
\[
 D_0v=-\mathcal N_p(v),\qquad
 \mathcal N_p(v)=(\phi_p+v)^p-\phi_p^p-p\phi_p^{p-1}v.
\]
Boundary division and the quadratic estimate in the proof of
Proposition~\ref{prop8} imply, for functions in its small unweighted ball,
\[
 \|\mathcal N_p(v)-\mathcal N_p(\widetilde v)\|_{L^q_{{\rm ul},\sigma}}
 \le C\eta_{\rm an}\|v-\widetilde v\|_{E_{-,\sigma}}.
\]
For $\sigma>0$, write $y=e^{-\sigma t}v$. The product rule gives
\[
 v_t=e^{\sigma t}(y_t+\sigma y),\qquad
 v_{tt}=e^{\sigma t}(y_{tt}+2\sigma y_t+\sigma^2y).
\]
Together with the angular and mixed derivative formulas and
$e^{\sigma t}\le1$ on $Q_-$, this yields
\[
 \|v\|_{2,q}\le C_{\rm emb}(\sigma)\|v\|_{E_{-,\sigma}},
 \qquad C_{\rm emb}(\sigma)\le C(1+\sigma+\sigma^2).
\]
Take $\sigma=\sigma_0$ from Proposition~\ref{prop8}. Both
$C_{\rm emb}(\sigma_0)$ and the weighted linear constant $C_{\sigma_0}$
are fixed uniformly on the truncated interval. Choose $c_*$ also to satisfy
\[
 4C_{\rm emb}(\sigma_0)C_{\sigma_0}c_*\le1.
\]
Then the ball of radius $2C_{\sigma_0}(\|z\|_Z+|b|)$ in
$E_{-,\sigma_0}$ lies in the unweighted radius-$\eta_{\rm an}$ ball
whenever $\|z\|_Z,|b|<c_*\eta_{\rm an}$.
Since $c_*$ is fixed before the solution radius, we may decrease
$\eta_{\rm an}$ further; the weighted linear inverse then makes
$v\mapsto\omega[p,z,-\mathcal N_p(v),b]$ a contraction in
$E_{-,\sigma_0}$. Its fixed point lies in the unweighted ball, where
uniqueness in Proposition~\ref{prop8} identifies it with
$w[p,z,b]-\phi_p$.
Thus exponential convergence is proved without assuming weighted
membership of the original bounded solution.
Iteration with the quadratic source then achieves any rate below
$\operatorname{Re}\mu_1^-$. If the first-mode roots $\mu_1^\pm$ are real and distinct,
this yields, for some $A_p\in\R$ and every
$\mu_1^-<\gamma<\min\{2\mu_1^-,\mu_1^+\}$,
\[
 w[p,z,b]=\phi_p+A_pe^{\mu_1^-t}\psi_p^1
             +O_{W^{2,q}_{\rm ul}}(e^{\gamma t}).
\]

For such a fixed background $w$, a Jacobi function $\eta$ with traces
$(z_1,b_1)$ solves $D_w\eta=0$, equivalently
\begin{equation}\label{e120}
 \eta=\omega[p,z_1,-V_w\eta,b_1],
 \qquad V_w=p(w^{p-1}-\phi_p^{p-1}).
\end{equation}
The small weighted perturbation first establishes a positive decay rate for
$\eta$. Bootstrapping the product $V_w\eta$ improves this to any rate
below $\mu_1^-$. Extracting the first mode in the linear estimate above
therefore yields
\[
 \eta=a_0e^{\mu_1^-t}\psi_p^1
             +O_{W^{2,q}_{\rm ul}}(e^{\gamma t})
\]
for the same range of $\gamma$. Only the strict real-part decay estimate
is used when the first roots are repeated or nonreal; the displayed
single slow-mode expansions are restricted to the real, distinct case.
\end{remark}

\subsection{The outer problem}\label{sec:continuation-outer}
We next solve the exterior problem using the boundary map $B$.
For its linearization, we first add a first-mode boundary condition
to obtain a coercive energy problem. We then remove that condition by
projection and impose zero far-field coefficient, obtaining the
index-zero reference operator. Consider the auxiliary problem
    \begin{equation}\label{e123}
\left\{\begin{aligned}
    -\Delta v&=g\quad\text{in}\quad \R^n_+\setminus\overline{B_1(0)},\\
    v&=0\quad\text{on}\quad\partial\R^n_+\setminus\overline{B_1(0)},\\
    \partial_\lambda v|_{S^{n-1}_+}&=C_pv|_{S^{n-1}_+}+\frac{2}{p-1}c_p^1[v|_{S^{n-1}_+}]\psi_p^1+h.
\end{aligned}\right.
    \end{equation}
Its projected variant has the same equation and flat boundary condition,
but the spherical boundary condition and far-field normalization are
    \begin{equation}\label{e122}
    \left\{\begin{aligned}
    P_p\partial_\lambda v|_{S^{n-1}_+}&=C_pv|_{S^{n-1}_+}+h,\\
    v(x)&=o(1)|x|^{-n}x_n\quad\text{as }|x|\to\infty.
    \end{aligned}\right.
    \end{equation}
     The solution formula for problem~\eqref{e116} gives
    \[
    C_pz=\sum_{j=2}^\infty \mu_p^{j,+}c_p^j[z]\psi_p^j
    \]
    for $z\in W^{2-q^{-1},q}\cap W^{1,1}_0\left(S^{n-1}_+\right)$. Define energy weak solutions to~\eqref{e123} as follows. Let $\mathcal{D}^{1,2}_0(\R^n_+\setminus B_1(0))$ be a completion of $C^\infty_c(\R^n_+)|_{\R^n_+\setminus\overline{B_1(0)}}$ with respect to the norm $\|\nabla\,\cdot\,\|_{L^2\left(\R^n_+\setminus\overline{B_1(0)}\right)}$.
For $\Sigma_+=S^{n-1}_+$, the critical Dirichlet trace space
$H^{1/2}_{00}(\Sigma_+)$ consists of the $H^{1/2}$ functions with finite norm
\begin{equation}\label{eq:critical-dirichlet-trace}
 \|z\|_{H^{1/2}_{00}}^2
 :=\|z\|_{H^{1/2}}^2+
 \int_{\Sigma_+}\frac{|z(\theta)|^2}
 {\operatorname{dist}(\theta,\partial\Sigma_+)}\,d\sigma.
\end{equation}
By the critical trace characterization
\cite[Chapter~\textup{1}, Theorem~\textup{11.7}]{lm},
$H^{1/2}_{00}=[L^2,H^1_0]_{1/2}$ is the space whose zero extension
belongs to $H^{1/2}(S^{n-1})$. For
$v\in\mathcal D^{1,2}_0(\R^n_+\setminus B_1)$, zero extension across
the physical plane, the trace theorem, and the angular Dirichlet
Poincar\'e inequality give
\begin{equation}\label{eq:energy-critical-trace}
 \|v|_{\Sigma_+}\|_{H^{1/2}_{00}}^2
 \le C\|v\|_{H^1(\{1<|x|<2,\ x_n>0\})}^2
 \le C\|\nabla v\|_{L^2(\R^n_+\setminus B_1)}^2.
\end{equation}

Define the bilinear form $\L_p$ on $(\mathcal{D}^{1,2}_0(\R^n_+\setminus B_1(0)))^2$ by
    \[
    \L_p(v,\varphi):=\int_{\R^n_+\setminus\overline{B_1(0)}}\nabla v\,\cdot\,\nabla\varphi\,dx+\sum_{j=2}^\infty \left(\mu_p^{j,+}-\frac{2}{p-1}\right)c_p^j[v|_{S^{n-1}_+}]c_p^j[\varphi|_{S^{n-1}_+}],
    \]
    for $v,\varphi\in \mathcal{D}^{1,2}_0(\R^n_+\setminus B_1(0))$.
Since $\partial_\lambda=\partial_r+2/(p-1)$ on $r=1$, the
boundary condition in \eqref{e123} is
\[
 \partial_rv=C_pv-\frac2{p-1}P_p(v|_{S^{n-1}_+})+h.
\]
The exterior normal there is $-\partial_r$, so Green's formula gives
the boundary datum with a minus sign in the weak formulation below.
Projection removes the auxiliary first-mode condition; the zero
far-field coefficient in \eqref{e122} is imposed separately.
    \begin{definition}\label{def1}
    Let $g$ and $h$ be measurable functions on
    $\R^n_+\setminus\overline{B_1(0)}$ and $S^{n-1}_+$, respectively.
    Suppose the functional
    \[
    \langle\g,\varphi\rangle:=
    \int_{\R^n_+\setminus\overline{B_1(0)}}g\varphi\,dx
    -\int_{S^{n-1}_+}h\varphi|_{S^{n-1}_+}\,d\sigma
    \]
    is bounded on $\mathcal D^{1,2}_0(\R^n_+\setminus B_1(0))$, so that
    $\g\in(\mathcal D^{1,2}_0(\R^n_+\setminus B_1(0)))^*$.
 A function $v\in \mathcal{D}^{1,2}_0(\R^n_+\setminus B_1(0))$ is a weak solution of~\eqref{e123} if
        \begin{equation}\label{e124}
        \L_p(v,\varphi)=\langle\g,\varphi\rangle\ \text{for all}\ \varphi\in \mathcal{D}^{1,2}_0(\R^n_+\setminus B_1(0)).
        \end{equation}
    \end{definition}
    \begin{lemma}\label{lem10}
The form $\mathcal L_p$ is bounded and coercive on
$\mathcal D^{1,2}_0(\R^n_+\setminus B_1)$, uniformly for $p\in[p_-,p_+]$.
For every bounded energy functional $\mathbf g$, equation~\eqref{e124}
has a unique energy weak solution.
\end{lemma}

\begin{proof}
After a fixed sufficiently large positive shift, the angular operator has form domain
$H^1_0(\Sigma_+)$ with a uniformly equivalent norm on $[p_-,p_+]$.
Interpolation identifies its quarter-power domain with
$H^{1/2}_{00}(\Sigma_+)$. Since
$\mu_p^{j,+}\asymp(1+\lambda_p^j)^{1/2}$ for large $j$, the spectral
theorem and \eqref{eq:energy-critical-trace} imply, for $z=v|_{\Sigma_+}$,
\begin{equation}\label{eq:critical-spectral-bound}
 \sum_{j=2}^\infty \mu_p^{j,+}|c_p^j[z]|^2
 \le C\|z\|_{H^{1/2}_{00}(\Sigma_+)}^2
 \le C\|\nabla v\|_{L^2(\R^n_+\setminus B_1)}^2.
\end{equation}
Thus $C_p:H^{1/2}_{00}\to(H^{1/2}_{00})^*$ and the energy form
$\mathcal L_p$ are bounded.

Write $a_p=2/(p-1)$ and $z=v|_{\Sigma_+}$. For a smooth function in
the dense energy class, radial integration by parts gives
\[
\begin{aligned}
 &\int_{|x|>1,\ x_n>0}|\nabla v|^2\,dx-(n-1)\|z\|_2^2\\
 &\quad=\int_1^\infty\!\!\int_{\Sigma_+}
 r^{n-1}\left|v_r+\frac{n-1}{r}v\right|^2\,d\sigma\,dr\\
 &\qquad+\int_1^\infty\!\!\int_{\Sigma_+}
 r^{n-3}\bigl(|\nabla_\theta v|^2-(n-1)|v|^2\bigr)\,d\sigma\,dr
 \ge0.
\end{aligned}
\]
The last step follows from the hemisphere Poincar\'e inequality;
density extends the identity to the energy space. Since $\mu_p^{j,+}>0$ for $j\ge2$,
\[
 \mathcal L_p(v,v)\ge
 \|\nabla v\|_2^2-a_p\|z\|_2^2
 \ge\left(1-\frac{a_p}{n-1}\right)\|\nabla v\|_2^2.
\]
The coefficient on the right is uniformly positive on $[p_-,p_+]$.
Lax--Milgram now proves the assertion.
\end{proof}

We next recover the decay needed for continuation from the energy
solution. Put $\nu:=\nu_0$, the fixed source weight chosen before
Theorem~\ref{thm4}. For $\ell\in\mathbb Z_{\ge0}$, set
    \begin{gather*}
    W^{\ell,q}_\nu\left(\R^n_+\setminus\overline{B_1(0)}\right):=\left\{\varphi:|\varphi|_{\ell,q;\nu}<\infty\right\},\\
    |\varphi|_{\ell,q;\nu}:=\sup_{r\ge1}r^{-\nu}\|\varphi(r\,\,\cdot\,)\|_{W^{\ell,q}\left(B_2(0)\setminus\overline{B_1(0)}\right)},\\
    W^{\ell,q}_{\nu;0}\left(\R^n_+\setminus B_1(0)\right):=\left\{\varphi\in W^{\ell,q}_\nu\left(\R^n_+\setminus\overline{B_1(0)}\right):\varphi|_{\partial\R^n_+\setminus\overline{B_1(0)}}=0\right\},\\
    L^q_\nu\left(\R^n_+\setminus\overline{B_1(0)}\right):=W^{0,q}_\nu\left(\R^n_+\setminus\overline{B_1(0)}\right).
    \end{gather*}
    Set also
    \begin{gather*}
    W:=\left\{v\in W^{2,q}_{-n+1;0}\left(\R^n_+\setminus B_1(0)\right):-\Delta v\in L^q_\nu\left(\R^n_+\setminus\overline{B_1(0)}\right)\right\},\\
    \|v\|_{W}:=|v|_{2,q;-n+1}+|-\Delta v|_{0,q;\nu}.
    \end{gather*}
For $\ell\ge2$ and $\varphi\in W^{\ell,q}_{\nu;0}$, the embedding $W^{\ell,q}\hookrightarrow C^1$ implies
\[
\sup |x|^{-\nu+1}x_n^{-1}|\varphi(x)|\le C|\varphi|_{\ell,q;\nu},
\]
so the set
\[
\left\{\varphi\in W^{\ell,q}_{\nu;0}:
\inf_{x\in\R^n_+\setminus\overline{B_1(0)}}
|x|^{-\nu+1}x_n^{-1}\varphi(x)>0\right\}
\]
is open in $W^{\ell,q}_{\nu;0}$.
    \begin{lemma}\label{lem11}
Let $\nu<-n-1$, $g\in L^q_\nu(\R^n_+\setminus\overline{B_1})$, and $h\in Z'$. If
\[
v\in\mathcal D^{1,2}_0(\R^n_+\setminus B_1)
\]
is a weak solution of~\eqref{e123}, then $v\in W$, $v$ is a strong solution, and
\begin{equation}\label{e126}
|v|_{2,q;-n+1}
\le C_{n,q,\nu,\epstr}
\left(
|g|_{0,q;\nu}+\|h\|_{Z'}
+\|\nabla v\|_{L^2(\R^n_+\setminus B_1)}
\right).
\end{equation}
The constant is uniform for $p\in(p_-,p_+)$. The limit
\[
\langle I,v\rangle
:=\lim_{|x|\to\infty}|x|^nx_n^{-1}v(x)
\]
exists and defines a continuous functional $I\in W^*$.
\end{lemma}
    The proof is given in Appendix~\ref{app:proof-lem11}.

Set
\[
\begin{gathered}
 E:=L^q_\nu(\R^n_+\setminus\overline{B_1})\times Z',\qquad
 \mathcal E_p:=L^q_\nu(\R^n_+\setminus\overline{B_1})\times P_pZ',\qquad
 F:=\ker I\subset W.
\end{gathered}
\]
Define $M:W\to E$, $M_1:W\to\mathcal E_p$ and
$M_2:F\to\mathcal E_p$ by
\[
\begin{aligned}
 Mv&:=\left(-\Delta v,\,
 \partial_\lambda v|_{S^{n-1}_+}-C_p[v|_{S^{n-1}_+}]
 -\frac{2}{p-1}c_p^1[v|_{S^{n-1}_+}]\psi_p^1\right),\\
 M_1&:=\pi M,\qquad M_2:=M_1|_F,\qquad
 \pi(g,h):=(g,P_ph).
\end{aligned}
\]
For $(g,h)\in E$, the bounds $g\in L^{2n/(n+2)}$ (since
$\nu_0<-n-1$) and $h\in Z'\subset L^2$, together with the Sobolev
and trace inequalities, make \eqref{e124} an energy problem.
Lemma~\ref{lem10} provides the unique energy weak solution; Lemma~\ref{lem11}
then places it in $W$. Conversely, $W$ embeds in the energy space. Thus
\[
M:W\longrightarrow E
\quad\hbox{is a bounded isomorphism}.
\]
The projected problem~\eqref{e122}, including its zero far-field
coefficient, is precisely $M_2v=(g,h)$ with $v\in F$.

\begin{proposition}\label{prop9}
The operators $M:W\to E$, $M_1:W\to\mathcal E_p$ and
$M_2:F\to\mathcal E_p$ are Fredholm, with
\[
\operatorname{ind}M=0,\qquad
\operatorname{ind}M_1=1,\qquad
\operatorname{ind}M_2=0.
\]
\end{proposition}
\begin{proof}
Since $M$ is an isomorphism and $\pi$ is surjective with kernel
$\operatorname{span}\{(0,\psi_p^1)\}$, index additivity yields indices
$0$ and $1$. The functional $I$ takes value $1$ on a cutoff of
$|x|^{-n}x_n$ equal to it near infinity, so $F=\ker I$ has codimension
one. Restriction to $F$ lowers the index to zero.
\end{proof}
The exterior equation and the inner matching condition can now be
combined in one analytic map. Let $\mathcal U_{\eta_{\rm an}}$ be as in
Proposition~\ref{prop8}, and define
    \begin{equation}\label{e132}
    \begin{gathered}
        A:\mathcal{V}_{\eta_{\rm an}}\to \mathcal E_p,\\
        \begin{aligned}
        A[p,u]:=\bigl(&-\Delta u-u^p,\,
        P_p\partial_\lambda u|_{S^{n-1}_+}\\
        &-B[p,u|_{S^{n-1}_+}-\phi_p,
        c_p^1[\partial_\lambda u|_{S^{n-1}_+}]]\bigr),
        \qquad u\in\mathcal{V}_{\eta_{\rm an}},
        \end{aligned}
    \end{gathered}
    \end{equation}
    where
    \begin{equation}\label{eq:matching-admissible-set}
    \begin{aligned}
    \mathcal{V}_{\eta_{\rm an}}:=\Bigl\{(p,u)&\in
    (p_-,p_+)\times W:\\
        &(p,u|_{S^{n-1}_+}-\phi_p)\in\mathcal{U}_{\eta_{\rm an}},\qquad
        \left|c_p^1[\partial_\lambda u|_{S^{n-1}_+}]\right|<\rho_{\eta_{\rm an}},\\
        &\langle I,u\rangle=c,\qquad
        \inf_{x\in\R^n_+\setminus\overline{B_1(0)}}
        |x|^{n}x_n^{-1}u(x)>0\Bigr\}.
    \end{aligned}
    \end{equation}
The trace conditions place the data in the analytic domain of $B$.
For the volume term, put $d(x)=|x|^{-n}x_n$. The weighted embedding
implies $|v|\le C\|v\|_W d$, while $u/d$ is bounded above and below
away from zero. Apply the positive-quotient argument of
Lemma~\ref{lem:positive-boundary-power}; the strict margin
$\nu_0+(n-1)p_->0$ also absorbs the logarithms at infinity from
$p$-derivatives. Thus $u^p$ is jointly analytic into $L^q_{\nu_0}$.
Proposition~\ref{prop8} and the analytic trace projections establish the
same conclusion for the boundary term.
     $u$ is a solution to problem~\eqref{e105} if and only if $A[p,u]=0$. The domain of $A$ lies in an affine subspace parallel to $\R\times F$, and $A[p,\,\cdot\,]$ takes values in $\mathcal{E}_p$.

    The linearized operator $\partial_u A[p,u]:F\to\mathcal{E}_p$ for the map $A$ is
    \begin{equation}\label{e133}
    \begin{aligned}
    &\partial_uA[p,u]v\\
        &=\left(-\Delta v-pu^{p-1}v,P_p\partial_\lambda v|_{S^{n-1}_+}
        -\partial_zB\left[p,u|_{S^{n-1}_+}-\phi_p,c_p^1\left[\partial_\lambda u|_{S^{n-1}_+}\right]\right]v|_{S^{n-1}_+}\right.\\
        &\qquad\left.-\partial_bB\left[p,u|_{S^{n-1}_+}-\phi_p,c_p^1\left[\partial_\lambda u|_{S^{n-1}_+}\right]\right]c_p^1[\partial_\lambda v|_{S^{n-1}_+}]\right).
    \end{aligned}
    \end{equation}
    \begin{proposition}\label{prop10}
Fix a sufficiently small $\epstr>0$ and the associated truncated $p$-interval and exterior weight $\nu_0$. Then there exists a small $\eta_{{\rm an},0}>0$ such that for every $(p,u)\in\mathcal{V}_{\eta_{{\rm an},0}}$, the operator $\partial_uA[p,u]:F\to\mathcal{E}_p$ is a Fredholm operator of index $0$.
\end{proposition}
\begin{proof}
By \eqref{e133}, \eqref{e117}, and \eqref{e118},
\[
(\partial_uA[p,u]-M_2)v
=\left(-pu^{p-1}v,\Pi_p[u]v\right),
\]
where
\[
z_u:=u|_{S^{n-1}_+}-\phi_p,
\qquad
b_u:=c_p^1[\partial_\lambda u|_{S^{n-1}_+}],
\]
and
\[
\begin{aligned}
\Pi_p[u]v:={}&-P_p\partial_z\partial_t\pi[p,z_u,b_u]
\bigl(v|_{S^{n-1}_+}\bigr)\big|_{t=0}\\
&-\partial_bB[p,z_u,b_u]\,
c_p^1[\partial_\lambda v|_{S^{n-1}_+}].
\end{aligned}
\]
By the trace estimate,
$\|\Pi_p[u]v\|_{Z'}\le C\eta_{{\rm an},0}\|v\|_W$.
Small bounded perturbations preserve the Fredholm index, so
$M_2+(0,\Pi_p[u])$ is Fredholm of index zero provided $\eta_{{\rm an},0}$ is sufficiently small.

For the potential term, the exterior weight satisfies
$\nu_0>-(n-1)p_-$.  For $r\ge2$, the weighted $W^{2,q}_{-n+1}$ bounds imply
$|u|+|v|\le Cr^{1-n}$ on rescaled annuli, hence
\[
r^{-\nu_0}
 \|u^{p-1}v(r\,\cdot)\|_{L^q(B_2^+\setminus\overline{B_1^+})}
\le Cr^{-\eta}\|v\|_W,
\qquad
\eta:=\nu_0+(n-1)p_->0.
\]
Rellich compactness on bounded annuli and this vanishing tail show
that $v\mapsto pu^{p-1}v$ is compact from $F$ to $L^q_{\nu_0}$.
Compact perturbation therefore preserves the index zero.
\end{proof}

To obtain an initial invertible point for continuation, we relate
the kernel of this exterior linearization to the global Jacobi fields.
\begin{lemma}[Linear matching]\label{lem:linear-matching}
Let $(p,u)\in\mathcal V_{\eta_{{\rm an},0}}$ satisfy $A[p,u]=0$,
with the radius small enough for \eqref{eq:actual-inner-source}.
Extend $u$ using $w[p,z,b]$. Here $a=2/(p-1)$,
$z=u|_{\Sigma_+}-\phi_p$ and $b=c_p^1[\partial_\lambda u]$.
Restriction to the exterior is a bijection between global Jacobi fields
\[
 (-\Delta-pu^{p-1})\varphi=0,\qquad
 \varphi\in W^{2,q}_{-a,1-n;0}(\R^n_+),
 \qquad \varphi=o(r^{-n}x_n)\quad(r\to\infty),
\]
and $\ker(\partial_uA[p,u]:F\to\mathcal E_p)$.
\end{lemma}
\begin{proof}
For an exterior kernel element $v$, put $\xi=v|_{\Sigma_+}$,
$d=c_p^1[\partial_\lambda v]$, and
\[
 \eta=\partial_z w[p,z,b]\xi+\partial_b w[p,z,b]d.
\]
Differentiation of the inner problem shows $D_w\eta=0$,
$\eta(0)=\xi$, $c_p^1[\eta_t(0)]=d$, and
\[
 P_p\eta_t(0)=\partial_zB[p,z,b]\xi+\partial_bB[p,z,b]d.
\]
The boundary component of \eqref{e133} therefore implies
$\eta_t(0)=\partial_\lambda v$. Define
$\varphi(r\theta)=r^{-a}\eta(\theta,\log r)$ for $r<1$ and
$\varphi=v$ for $r>1$. The values agree, and
\[
 \partial_r\varphi|_{r=1}=\eta_t(0)-a\eta(0)
 =\partial_\lambda v-av=\partial_rv|_{r=1}.
\]
There is no interface source, and the inner estimates ensure
$|\varphi|_{2,q;-a,1-n}\le C\|v\|_W$.
Since $v\in F=\ker I$, the far-field conclusion of
Lemma~\ref{lem11} ensures the required little-o condition.

Conversely, the equation and $\nu_0+(n-1)p>0$ put the exterior
restriction of a field in the stated class in $W$, and its coefficient
is zero. In the inner cylinder, subtracting the two displayed
homogeneous variations leaves a bounded solution with zero terminal
value and first-mode derivative. The zero-data uniqueness in
\eqref{eq:actual-inner-source} makes it zero. The remaining derivative
matching is exactly the boundary component of \eqref{e133}.
\end{proof}

\subsection{Completion of the proof of Theorem~\ref{thm1}}\label{sec:proof-cell-existence}
\subsubsection{Verification of the continuation hypotheses}
We now verify hypotheses \textup{(I)}--\textup{(IV)} of
Theorem~\ref{thm:boundary-continuation} for the fixed-coefficient
matching map. Conclusions
\textup{(a)} and \textup{(c)} of Theorem~\ref{thm:boundary-continuation} provide escape
from each exhaustion set and eventual departure from the initial graph.
Conclusion \textup{(b)} provides the local analytic parametrizations in
Theorem~\ref{thm4}\,\ref{it4}. At a joining point, the homeomorphic
change of the curve parameter is supplied by the pairing of analytic
arcs in \cite[Theorem~\textup{6.1}, pp.~563--564]{cww-bif}; restricting and
rescaling its interval defines the injective map $s^*$ in that assertion.
The stronger non-degeneracy H\textup{(i)} requires the separate argument in
Section~\ref{sec:generic-nondegeneracy}.

\subsubsection{Compactness in the fixed spaces}
The fixed-coefficient domain and target are identified with Banach spaces as follows.

\begin{remark}\label{rem6}
The coefficient functional $I$ is independent of $p$. Set $F=\ker I$
and choose a fixed $e\in W$ with $\langle I,e\rangle=1$. The level set
$\{u\in W:\langle I,u\rangle=c\}$ is an affine Banach space, so the
translation $u\mapsto u-ce$ identifies it with $F$. Thus writing
\[
u=ce+v,\qquad v\in F,
\]
identifies the fixed-coefficient domain with the open set
\[
\mathcal U_*:=\{(p,v)\in\R\times F:(p,ce+v)\in\mathcal V_{\eta_{{\rm an},0}}\}.
\]
Only the boundary target space depends on $p$. We trivialize it by
the fixed projection associated with $\theta_n$: set
\[
P_0Z'=\left\{h\in Z':\int_{S^{n-1}_+}h\theta_n\,d\sigma=0\right\}.
\]
Let $\langle\cdot,\cdot\rangle$ denote the $L^2(\Sigma_+)$ pairing and set
\[
 P_0h=h-\frac{\langle h,\theta_n\rangle}{\|\theta_n\|_2^2}\theta_n,
 \qquad
 S_pg=g-\frac{\langle g,\psi_p^1\rangle}
                 {\langle\theta_n,\psi_p^1\rangle}\theta_n.
\]
Then $P_0:P_pZ'\to P_0Z'$ is an isomorphism with inverse $S_p$.
The denominator is positive and uniformly bounded away from zero on
compact parameter intervals. Lemma~\ref{lem:profile-analytic} therefore
makes $S_p$ analytic. With
\[
 \mathcal E_0=L^q_\nu(\R^n_+\setminus\overline{B_1})\times P_0Z',
 \qquad
 \widehat A(p,v)=\bigl(A_1[p,ce+v],P_0A_2[p,ce+v]\bigr),
\]
the zeros agree, and
\[
 \partial_v\widehat A=
 \begin{pmatrix}I&0\\0&P_0\end{pmatrix}\partial_uA\big|_F.
\]

Thus the kernel and Fredholm index are unchanged. We will apply
Theorem~\ref{thm:boundary-continuation} with $X=F$, $Y=\mathcal E_0$ and open domain
$\mathcal U_*$, translated at the initial endpoint as specified below.
\end{remark}

When no confusion is possible, we retain the notation $(p,u)$ and $A$
for these fixed-coordinate objects; the continuation theorem acts on $\mathcal U_*\subset\R\times F$
with target $\mathcal E_0$.

We prove compactness of the exterior restrictions for solutions with
fixed coefficient $c$ and $p\in[p_-,p_+]$, including sequences
approaching the boundary of the admissible set. Let $(p_j,u_j)$ be
such a sequence. Put $K(x)=|x|^{-n}x_n$
and $\mathcal A=\{1<|x|<2,\ x_n>0\}$.
The local estimates and Lemma~\ref{lem8} give, after passing to a subsequence,
$p_j\to p_\infty$ and strong $W^{2,q}$ convergence on each compact
half-annulus. The same lemma provides
\begin{equation}\label{eq:fixed-c-compactness}
 r^{n-1}\|u_j(r\,\cdot)-cK(r\,\cdot)\|_{W^{2,q}(\mathcal A)}
 \le Cr^{-\kappa},\qquad r\ge2,
\end{equation}
with $C,\kappa>0$ independent of $j$. The common leading terms cancel
in differences, so the weighted $W^{2,q}$ tails are uniformly small.
Moreover, since $u_j\le Cr^{1-n}$ and $p_j\ge p_-$, the source satisfies
\begin{equation}\label{eq:source-tail-bound}
 r^{-\nu_0}\|u_j^{p_j}(r\,\cdot)\|_{L^q(\mathcal A)}
 \le Cr^{-\nu_0-(n-1)p_j}
 \le Cr^{-\eta},\qquad \eta=\nu_0+(n-1)p_->0.
\end{equation}
Local convergence of the sources and this tail bound ensure convergence
in $L^q_{\nu_0}$. First taking the exterior radius large and then
letting $j,\ell\to\infty$ shows that $(u_j)$ is Cauchy in $W$.
The resulting limit satisfies
\[
u_j\longrightarrow \widehat u\quad\hbox{in }W,
\qquad \langle I,\widehat u\rangle=c,
\]
where the coefficient identity follows from continuity of $I$.

Proposition~\ref{prop10} verifies the Fredholm hypothesis \textup{(III)} of
Theorem~\ref{thm:boundary-continuation}. To verify \textup{(IV)}, set $\mathbb X=\R\times F$ and
\begin{equation}\label{eq:continuation-exhaustion}
 Q_j=\left\{z\in\mathcal U_*:\|z\|_{\mathbb X}\le j,\quad
       \operatorname{dist}_{\mathbb X}(z,\mathbb X\setminus\mathcal U_*)\ge j^{-1}\right\}.
\end{equation}
These bounded closed sets increase to exhaust $\mathcal U_*$. The
preceding strong compactness, closedness of $Q_j$, and continuity of
$\widehat A$ ensure compactness of $Q_j\cap\widehat A^{-1}(0)$, as required
by \textup{(IV)}.

\subsubsection{The initial stable branch}
For the initial arc required by Theorem~\ref{thm:boundary-continuation}, take the solutions of \eqref{e3} constructed in \cite{dmp} for $p$ close to $\frac{n+1}{n-1}$.
    \begin{lemma}\label{lem13}
There exists $\epstr>0$ such that, for every
\[
\frac{n+1}{n-1}<p<\frac{n+1}{n-1}+2\epstr,
\]
the following hold:
\begin{enumerate}[label=(\roman*)]
\item \label{it29} Problem~\eqref{e3} subject to \eqref{e6} has a unique solution
\[
u_p\in C^2(\overline{\R^n_+}\setminus\{0\}).
\]
\item \label{it30} If
\[
v\in W^{2,q}_{-\frac{2}{p-1},-n+1;0}(\R^n_+)
\]
satisfies
\begin{equation}\label{e134}
\begin{gathered}
\begin{cases}
-\Delta v-pu_p^{p-1}v=0&\text{in }\R^n_+,\\
v=0&\text{on }\partial\R^n_+\setminus\{0\},
\end{cases}\\
v(x)=o(|x|^{-n}x_n)\qquad(|x|\to\infty),
\end{gathered}
\end{equation}
then $v\equiv0$.
\end{enumerate}
\end{lemma}

\begin{proof}
By \eqref{e79} and $\phi_p\to0$ as
        $p\downarrow\frac{n+1}{n-1}$, after shrinking the interval of $p$
        there exists a constant
        \[
        0<\kappa<\left(\frac{n-2}{2}\right)^2
        \]
        such that, for every fast-decaying solution $u$ under consideration,
        \begin{equation}\label{e135}
          pu^{p-1}(x)\le\kappa |x|^{-2},\qquad
          p\phi_p^{p-1}(\theta)\le\kappa.
        \end{equation}
        Shrink the interval further so that the characteristic roots
        $\mu_1^\pm$ of the first mode are real and distinct.  With
        \[
        a_p:=\frac2{p-1},\qquad
        \alpha_p=2a_p-(n-2),
        \]
        Proposition~\ref{prop7} and continuity from the lower endpoint imply
        \[
        0<\mu_1^-<\frac{\alpha_p}{2}<\mu_1^+,
        \qquad 2\mu_1^-<\mu_1^+.
        \]
        Set
        \[
        \gamma_0:=\frac32\mu_1^-.
        \]
        Then $\mu_1^-<\gamma_0<\min\{2\mu_1^-,\mu_1^+\}$.  In particular,
        \begin{equation}\label{e135a}
        -a_p+\mu_1^+> -\frac{n-2}{2}.
        \end{equation}

The existence of $u_p$ follows from \cite{dmp}. Set
\[
w_p(\theta,t):=e^{\frac{2}{p-1}t}u_p(e^t\theta).
\]
Remark~\ref{rem5} provides
        \begin{equation}\label{e136}
        u_p-\left(u_p^*+A_p|x|^{-\frac{2}{p-1}+\mu_1^-}\psi_p^1\left(\frac{x}{|x|}\right)\right)\zeta\in W^{2,q}_{-\frac{2}{p-1}+\gamma_0,-n+1;0}(\R^n_+).
        \end{equation}
Here $\zeta\in C^\infty_c\left(\overline{\R^n_+}\right)$ is such that $\zeta=1$ on $B_1^+(0)$, and $u_p^*$ is the homogeneous solution introduced in Section~\ref{sec:prelim}.

\proofstep{1}{Show that the first-mode coefficient is nonzero.}
To determine the sign of $A_p$, set
        \[
        V_p(\theta,t):=w_p(\theta,t)-\phi_p(\theta),\qquad
        b_1(t):=\int_{S^{n-1}_+}V_p(\theta,t)\psi_p^1(\theta)\,d\sigma.
        \]
        Projecting the equation for $V_p$ onto the positive first
        eigenfunction gives
        \begin{equation}\label{e137}
        b_1''-\alpha_pb_1'-\lambda_p^1b_1=-F_p(t),
        \end{equation}
        where
        \[
        F_p(t):=\int_{S^{n-1}_+}
        \bigl[(\phi_p+V_p)^p-\phi_p^p-p\phi_p^{p-1}V_p\bigr]
        \psi_p^1\,d\sigma\ge0.
        \]
        The inequality follows from the strict convexity of $s\mapsto s^p$.
        Moreover, $F_p\not\equiv0$, since $w_p\to0$ as $t\to+\infty$
        whereas $\phi_p\not\equiv0$.  Remark~\ref{rem5} implies
        $F_p(t)=O(e^{2\mu_1^-t})$ as $t\to-\infty$, while $F_p$ is bounded
        as $t\to+\infty$.  Since both roots of
        $\mu^2-\alpha_p\mu-\lambda_p^1=0$ are positive, variation of
        constants in \eqref{e137}, using boundedness of $b_1$ at $+\infty$,
        yields
        \begin{equation}\label{e137a}
        A_p=-\frac1{\mu_1^+-\mu_1^-}
        \int_{-\infty}^{+\infty}e^{-\mu_1^-s}F_p(s)\,ds<0.
        \end{equation}
        Thus $A_p\neq0$.  Every other fast-decaying solution satisfying \eqref{e6} has the same negative sign for its first coefficient.

\proofstep{2}{Prove uniqueness of the cell.}
Let $\overline u$ be another solution of \eqref{e3} satisfying \eqref{e6}. As in \eqref{e136},
        \[
        \overline{u}-\left(u_p^*+\overline{A}_p|x|^{-\frac{2}{p-1}+\mu_1^-}\psi_p^1\left(\frac{x}{|x|}\right)\right)\zeta\in W^{2,q}_{-\frac{2}{p-1}+\gamma_0,-n+1;0}(\R^n_+),
        \]
        for some $\overline{A}_p<0$. Using the scaling invariance of \eqref{e3}, choose $\rho>0$ so that
        \[
        \overline u_\rho(x):=
        \rho^{\frac{2}{p-1}}\overline u(\rho x)
        \]
        has the same coefficient $A_p$; explicitly,
        $\rho^{\mu_1^-}\overline A_p=A_p$. Then
        \[
        \overline{\varphi}:=\overline{u}_\rho-u_p
        \in W^{2,q}_{-\frac{2}{p-1}+\gamma_0,-n+1;0}(\R^n_+).
        \]
        Cancellation of the common $\mu_1^-$ coefficient
        improves the decay through the linear equation.  In cylinder variables
        \[
        H(\theta,t):=e^{a_pt}\overline\varphi(e^t\theta)
        \]
        satisfies
        \[
        D_0H=-q_p(\theta,t)H,
        \]
        where, with $\overline w_\rho(\theta,t):=
        e^{a_pt}\overline u_\rho(e^t\theta)$,
        \[
        q_p(\theta,t):=
        \frac{\overline w_\rho(\theta,t)^p-w_p(\theta,t)^p}
        {\overline w_\rho(\theta,t)-w_p(\theta,t)}
        -p\phi_p(\theta)^{p-1},
        \]
        with the quotient interpreted by continuity when the denominator
        vanishes.  Remark~\ref{rem5}, the boundary $C^1$ estimate, and the Hopf lemma imply
        $w_p/\phi_p=1+O(e^{\mu_1^-t})$ and
        $\overline w_\rho/\phi_p=1+O(e^{\mu_1^-t})$ up to the lateral
        boundary. Factoring out $\phi_p$ yields
        $q_p=O(e^{\mu_1^-t})$ uniformly there, also for $1<p<2$.
        Since the $e^{\mu_1^-t}\psi_p^1$ coefficient of $H$ is zero,
        the second weighted estimate in Remark~\ref{rem5} shows that
        $H=O(e^{\gamma t})$ improves to every exponent below
        $\min\{\gamma+\mu_1^-,\mu_1^+\}$.  Starting from
        $\gamma_0=3\mu_1^-/2$, finitely many iterations establish
        \begin{equation}\label{e137b}
        \overline\varphi\in
        W^{2,q}_{-\frac{2}{p-1}+\gamma_*,-n+1;0}(\R^n_+)
        \quad\text{for some}\quad
        \frac{\alpha_p}{2}<\gamma_*<\mu_1^+.
        \end{equation}
        The improved exponent and the fast decay at infinity imply
        \[
        -a_p+\gamma_*>-\frac{n-2}{2},
        \qquad \overline\varphi\in\mathcal D^{1,2}_0(\R^n_+).
        \]
        We can therefore use the energy identity. The difference satisfies
        \[
        -\Delta\overline\varphi=p\widetilde u^{p-1}\overline\varphi,
        \]
        where $\widetilde u$ lies pointwise between $\overline u_\rho$ and
        $u_p$.  Estimate \eqref{e135} is invariant under scaling, hence
        $p\widetilde u^{p-1}\le\kappa |x|^{-2}$.  Testing by
        $\overline\varphi$ and using Hardy's inequality gives
        \[
        \int_{\R^n_+}|\nabla\overline\varphi|^2
        \le \kappa\int_{\R^n_+}|x|^{-2}\overline\varphi^2
        \le \frac{4\kappa}{(n-2)^2}
        \int_{\R^n_+}|\nabla\overline\varphi|^2.
        \]
        Since $4\kappa/(n-2)^2<1$, Hardy's inequality forces
        $\overline\varphi\equiv0$. Hence $\overline u_\rho\equiv u_p$. Comparing their far-field coefficients implies
        \[
        \rho^{\frac2{p-1}+1-n}=1.
        \]
        Since $\frac2{p-1}+1-n\ne0$ in the present range, $\rho=1$, and therefore $\overline u\equiv u_p$. The proof of assertion~\ref{it29} is complete.

\proofstep{3}{Show that the normalized linearization has trivial kernel.}
For assertion~\ref{it30}, Remark~\ref{rem5} provides a slow-mode coefficient
$a_0$ for $v$. By \eqref{e136}, the scaling Jacobi field
$Z_p=\partial_\lambda u_p$ has coefficient $\mu_1^-A_p\ne0$.
Consequently
\[
 \varphi:=v-\frac{a_0}{\mu_1^-A_p}Z_p
 \in W^{2,q}_{-a_p+\gamma_0,1-n;0}(\R^n_+)
\]
solves the homogeneous linearized equation and has no slow mode.
Its cylindrical perturbation is
$p(w_p^{p-1}-\phi_p^{p-1})=O(e^{\mu_1^-t})$, uniformly up to the
lateral boundary by the boundary-factor estimate above. The bootstrap
leading to \eqref{e137b} and the ensuing energy-space and Hardy argument
therefore apply unchanged and prove $\varphi=0$.
Finally, differentiating the normalized far-field expansion using
scaled elliptic estimates yields
\[
 \lim_{|x|\to\infty}|x|^nx_n^{-1}Z_p(x)=a_p+1-n\ne0.
\]
The zero far-field coefficient in \eqref{e134} now yields
\[
 0=\frac{a_0}{\mu_1^-A_p}(a_p+1-n),
\]
so $a_0=0$ and $v=0$.
    \end{proof}

\subsubsection{Continuation of the branch and its endpoints}
\begin{proof}[Proof of Theorem~\ref{thm4}: continuation with coefficient $c$]
\proofstep{1}{Choose the neighborhood and the coefficient.}
Fix the truncated interval of Lemma~\ref{lem13}. Choose $\eta_{{\rm an},0}$
for Propositions~\ref{prop8}, \ref{prop10} and
\eqref{eq:actual-inner-source}, then set $\rho_0=c_*\eta_{{\rm an},0}$.
Before choosing $c$, we record the first-mode derivative estimate.
Let $w$ solve \eqref{e73} on $Q$, and put $v=w-\phi_p$ and
$b_1(t)=c_p^1[v(\cdot,t)]$. Projection of the
equation onto $\psi_p^1$ gives
\[
 b_1''-\alpha_pb_1'-\lambda_p^1b_1
 =-c_p^1[\mathcal N_p(v)].
\]
If $\sup_{t\le0}\|v(\cdot,t)\|_{C^2}\le\eta$, boundary division and
the quadratic remainder estimate imply
$\|b_1\|_{L^\infty(-1,0)}\le C\eta$ and
$\|c_p^1[\mathcal N_p(v)]\|_{L^\infty(-1,0)}\le C\eta^2$.
The one-sided local estimate for this scalar ODE therefore implies
\begin{equation}\label{eq:first-mode-derivative-bound}
 |c_p^1[\partial_tw(\cdot,0)]|=|b_1'(0)|
 \le C\bigl(\|b_1\|_{L^\infty(-1,0)}
       +\|c_p^1[\mathcal N_p(v)]\|_{L^\infty(-1,0)}\bigr)
 \le C_*\eta.
\end{equation}
All constants are uniform on the fixed truncated exponent interval.

Choose $\eta_{\rm an}>0$ so small that
    \begin{equation}\label{e138}
\begin{gathered}
 \|z\|_{C^2(\overline{S^{n-1}_+})}<2\eta_{\rm an}
       \ \Longrightarrow\ \|z\|_Z<\rho_0/2,\\
 2C_*\eta_{\rm an}<\rho_0/2,
\end{gathered}
    \end{equation}
    uniformly for $p$ in the truncated interval, and let $c>0$ be as in Lemma~\ref{lem14} with this value of $\eta_{\rm an}$. After rescaling, the far-field normalization in~\eqref{e6} becomes
    \begin{equation}\label{e139}
      u(x)=(c+o(1))|x|^{-n}x_n\quad\text{as}\quad |x|\to+\infty.
    \end{equation}
    In the remainder of the proof, $u_p$ denotes the rescaled
    representative from Lemma~\ref{lem13} satisfying \eqref{e139}.
    If $u$ is a solution of~\eqref{e3} with conditions~\eqref{e6}, then, after the normalization \eqref{e139}, Lemma~\ref{lem14}, the first-mode derivative estimate~\eqref{eq:first-mode-derivative-bound},
and \eqref{e138} imply
    \[
    (p,u|_{\R^n_+\setminus\overline{B_1(0)}})\in
    \mathcal{V}_{\eta_{{\rm an},0}}.
    \]
    The exterior positivity in the definition of $\mathcal V_{\eta_{{\rm an},0}}$ follows from positivity of $u$, the boundary Hopf lemma on compact annuli, and the positive coefficient $c$ at infinity.

\proofstep{2}{Construct the initial branch and identify its endpoint.}
Set
\[
\mathcal A^+
:=
\left\{\left(p,u_p|_{\R^n_+\setminus\overline{B_1(0)}}\right)\right\}_{p\in\left(\frac{n+1}{n-1}+\epstr,\frac{n+1}{n-1}+2\epstr\right)}.
\]
Lemma~\ref{lem:linear-matching} extends a kernel element of
$\partial_uA[p,u_p]$ to a global Jacobi field with zero far-field
coefficient. Lemma~\ref{lem13}\,\ref{it30}, invariant under fixed dilation,
makes it zero. Index zero (Proposition~\ref{prop10}) therefore implies
invertibility. By the analytic implicit function theorem and
Lemma~\ref{lem13}\,\ref{it29}, we obtain the analytic initial arc.

For the boundary endpoint required by Theorem~\ref{thm:boundary-continuation},
observe that $p_-$ lies strictly inside the
near-endpoint existence and uniqueness interval of Lemma~\ref{lem13}.
For $p_j\downarrow p_-$, inner compactness and
\eqref{eq:fixed-c-compactness} produce a coefficient-$c$ limit in $W$.
They apply at $p_-$ because it is an interior exponent of the uniqueness
interval and $\nu_0+(n-1)p_->0$. The inner data retain strict admissibility
margins, so continuity of the inner map and Remark~\ref{rem5} determine the
singular profile. Lemma~\ref{lem13}\,\ref{it29} identifies every limit with
$u_{p_-}$; hence the full family converges:

\begin{equation}\label{eq:initial-boundary-endpoint}
 u_p|_{\R^n_+\setminus\overline{B_1}}\longrightarrow
 u_{p_-}|_{\R^n_+\setminus\overline{B_1}}\quad\hbox{in }W,
 \qquad p\downarrow p_-.
\end{equation}
In the coordinates of Remark~\ref{rem6}, set
$v_0=u_{p_-}|_{\R^n_+\setminus\overline{B_1}}-ce\in F$.
Then $(p,u_p-ce)\to(p_-,v_0)\in\partial\mathcal U_*$, since
$p_-$ is on the parameter face. Translate by
$\lambda=p-p_-$, $x=v-v_0$ and
$\mathfrak F(\lambda,x)=\widehat A(p_-+\lambda,v_0+x)$.
The initial graph now tends to $(0,0)$ on the domain boundary.
Lemma~\ref{lem13}, matching, and the analytic implicit function theorem
verify \textup{(I)}--\textup{(II)}; translation preserves the Fredholm property and
compact exhaustion, so all hypotheses of Theorem~\ref{thm:boundary-continuation} hold.

\proofstep{3}{Continue the branch globally.}
The continuation theorem provides a continuous global branch
\[
    \mathcal A^+\subset\mathcal A=\{(p_s,u_s)\}_{s>0}
    \subset\mathcal V_{\eta_{{\rm an},0}}.
\]
Extend $u_s$ to $\R^n_+$ by \eqref{e106}; then $u_s$ solves \eqref{e3}.
Its singular asymptotic follows from the inner construction; its
far-field asymptotic follows from membership in
$\mathcal V_{\eta_{{\rm an},0}}$ and \eqref{e139}. The local analytic reparametrization, together with the
analytic inner solution map, proves assertions~\ref{it1}--\ref{it4}.
The initial endpoint limit has already been established. It remains
to prove
\begin{equation}\label{e141}
\lim_{s\to\infty}p_s=\frac{n+2}{n-2}-\epstr.
\end{equation}
By Theorem~\ref{thm:boundary-continuation}\textup{(a)},
\begin{equation}\label{eq:global-leaves-bounded-closed}
 \text{for every }j,\quad (p_s,v_s)\notin Q_j
 \quad\text{for all sufficiently large }s.
\end{equation}

\proofstep{4}{Show that the exponent approaches the upper endpoint.}
The local estimates, \eqref{eq:fixed-c-compactness}, and
\eqref{eq:source-tail-bound}, with
$\nu_0>-(n-1)((n+1)/(n-1)+\epstr)$, give
\begin{equation}\label{eq:global-W-bound}
\sup_{s>0}\|u_s\|_W<\infty.
\end{equation}
Lemma~\ref{lem14}, \eqref{e138}, and
\eqref{eq:first-mode-derivative-bound} place the normalized data
\eqref{e139} strictly inside the domain of $B$:
\[
 \|u_s|_{S^{n-1}_+}-\phi_{p_s}\|_Z<\rho_0/2,
 \qquad
 |c_{p_s}^1[\partial_\lambda u_s]|
 \le C_*\eta_{\rm an}<\rho_0/2.
\]
Set $m(u)=\inf_{\R^n_+\setminus\overline{B_1}}|x|^nx_n^{-1}u(x)$.
If $m(u_j)\to0$ along the branch, \eqref{eq:fixed-c-compactness}
gives, after passage to a subsequence, $u_j\to\widehat u$ in $W$
and locally in $C^1$. Since $I\in W^*$, $\langle I,\widehat u\rangle=c>0$.
The uniform remainder \eqref{eq:cell-uniform-far-remainder} gives
$|x|^nx_n^{-1}u_j\ge c/2$ for $|x|\ge R$, with $R$ fixed.
On $1\le|x|\le R$, the strong maximum principle and Hopf lemma give
a positive minimum for $\widehat u/x_n$, including its extension to the
flat boundary. Local $C^1$ convergence transfers this lower bound to
$u_j$, a contradiction. Thus
\begin{equation}\label{eq:global-positive-margin}
\inf_{s>0}m(u_s)>0.
\end{equation}

The parameter also stays uniformly above $p_-$ as $s\to\infty$.  If not, there
would be $s_j\to\infty$ such that
\[
p_{s_j}\longrightarrow\frac{n+1}{n-1}+\epstr.
\]
For all large $j$ these parameters lie on the interval of uniqueness in
Lemma~\ref{lem13}\,\ref{it29}. After the fixed-coefficient normalization,
the corresponding solutions therefore belong to $\mathcal A^+$.
This contradicts Theorem~\ref{thm:boundary-continuation}\textup{(c)}, which excludes return
to that initial graph for large $s$. Thus there is $\eta_->0$ such that
\begin{equation}\label{eq:global-lower-p-margin}
p_s\ge\frac{n+1}{n-1}+\epstr+\eta_-
\qquad\text{for all sufficiently large }s.
\end{equation}

Suppose now that \eqref{e141} fails.  Then for some $\eta_+>0$ there is a
sequence $s_j\to\infty$ such that
\[
p_{s_j}\le\frac{n+2}{n-2}-\epstr-\eta_+.
\]
The uniform $W$-bound and \eqref{eq:fixed-c-compactness} make these
points precompact in $\mathbb X$. The strict margins for the inner datum and
the first-mode parameter, \eqref{eq:global-positive-margin}, and
\eqref{eq:global-lower-p-margin} put every limit in $\mathcal U_*$.
Their closure is therefore a compact subset of $\mathcal U_*$; it has
bounded norm and positive distance from $\mathbb X\setminus\mathcal U_*$.
It is contained in some $Q_j$, contradicting
\eqref{eq:global-leaves-bounded-closed}. Hence \eqref{e141} holds.

We have obtained the coefficient-$c$ branch, locally analytic in the
fixed exterior space $W$ and, by Proposition~\ref{prop8}, in the inner
space $W^{2,q}_{\ul}(\overline{Q_-})$. Together with the two endpoint
limits, this proves the theorem.
\end{proof}

\begin{remark}[Normalization]\label{rem:normalization-continuity}
For a continuous coefficient-$c$ family in these fixed coordinates, put
\[
 a_s=\frac2{p_s-1},\qquad
 \lambda_s=c^{1/(n-1-a_s)},\qquad
 \widetilde u_s(x)=\lambda_s^{a_s}u_s(\lambda_sx).
\]
Then $\widetilde u_s$ has coefficient $1$.
Writing $h_s=\log\lambda_s$, $w_s=e^{a_st}u_s(e^t\theta)$ and
$v_s=e^{(n-1)t}u_s(e^t\theta)$ on the full cylinder gives
\[
 \widetilde w_s(\theta,t)=w_s(\theta,t+h_s),\qquad
 \widetilde v_s(\theta,t)=c^{-1}v_s(\theta,t+h_s).
\]
Continuity in the fixed uniformly local $W^{2,q}$ spaces follows
from Sobolev translation on finitely many matching strips and uniform
smallness of the tails of $w_s-\phi_{p_s}$ and $v_s-c\theta_n$
(Remark~\ref{rem5} and Lemma~\ref{lem8}). The equation
$-\Delta\widetilde u_s=\widetilde u_s^{p_s}$ and the positive margin
$\nu_0+(n-1)p_s>0$ control the source tails, giving continuity in $W$.
At fixed $p$, dilation identifies the two-ended weighted spaces and
conjugates the Jacobi operators, preserving their kernels.
Analyticity is used only before normalization.

\end{remark}

\begin{proof}[Proof of Theorem~\ref{thm1}]
Fix $p$ in the strict subcritical interval and choose a truncation
with $p\in(p_-,p_+)$. The endpoint limits and continuity in
Theorem~\ref{thm4} give some $s$ with $p_s=p$.
The inner construction and Remark~\ref{rem5} ensure a positive decay
rate for $w_s-\phi_p$ in the uniformly local $W^{2,q}$ norm.
Since $q>n$, boundary division by $\theta_n$ and
$\phi_p\ge c_p\theta_n$ convert this into the relative singular
asymptotic in \eqref{e6}. Lemma~\ref{lem8} establishes the relative
far-field asymptotic, uniformly up to the flat boundary.
Remark~\ref{rem:normalization-continuity} normalizes its coefficient to
$1$. Proposition~\ref{prop6} proves the global boundary-factor upper
bound, preserved by this dilation. This proves all the assertions.
\end{proof}

\section{Non-degeneracy of the half-space cells}\label{sec:model-linear}
We prove non-degeneracy for the stable cells, their local threshold
continuation, and the global analytic branches, without assuming
global uniqueness.

\subsection{Non-degeneracy assumptions and weights}\label{sec:assumption-H}
Fix a normal dimension $n\ge3$, an exponent in the strict subcritical
range, and a coefficient-$1$ cell $u_1$. Set
$L_p=L_1:=\Delta_X+pu_1^{p-1}$.
H\textup{(i)} is the two-ended injectivity condition for $L_p$; H\textup{(ii)} is its
product-domain analogue. We derive H\textup{(ii)} from shifted injectivity
of $L_p-|\xi|^2$.

We use the general two-ended spaces defined in
Section~\ref{sec:general-model-spaces}. Put $a_p=2/(p-1)$ and fix $q>n$. Choose nonindicial weights satisfying
\begin{equation}\label{e59-twoend}
\operatorname{Re}\gamma_1^+(p)<\mu<2-n+a_p,
\qquad 1-n<\nu<-\frac{n-2}{2}.
\end{equation}
The upper bound on $\mu$ permits a compatible low weight in the gluing
problem. The interval is nonempty because
$-a_p<\operatorname{Re}\gamma_1^-$ and
$\operatorname{Re}\gamma_1^-+\operatorname{Re}\gamma_1^+=2-n$;
the first inequality follows from \eqref{e16} for real roots and from
$a_p>(n-2)/2$ for nonreal roots.

\medskip
\noindent{\bf Assumption H.}
Put $\mathcal O_k:=\R^k\times
(\overline{\R^n_+}\setminus\{0\})$ and
$\langle Y\rangle=(1+|Y|^2)^{1/2}$. For arbitrary real weights
$\mu,\nu$, define the following function classes:
\begin{equation}\label{eq:model-pointwise-class}
\begin{aligned}
 \mathscr X^m_{\mu,\nu}(\mathcal O_k)
 :=\biggl\{\Phi\in W^{2,q}_{\rm loc}(\mathcal O_k):\ &
 \Phi|_{X_n=0}=0,\\[-2pt]
 &\sup_{Y,X}\frac{|\Phi(Y,X)|}
 {\langle Y\rangle^m\omega_{\mu,\nu}(|X|)}<\infty\biggr\},
 \qquad m\ge0.
\end{aligned}
\end{equation}
For $k=0$ this is $W^{2,q}_{\rm loc}\cap C^0_{\mu,\nu}$ with zero
trace. The trace is imposed on the physical boundary, that is, on
$X_n=0$, $X\ne0$; the singular edge $\{X=0\}$ is controlled by the weights. Local elliptic
regularity permits increasing $q$ for solutions when needed.
For Assumption~H, take the weights in~\eqref{e59-twoend} and require:
\begin{enumerate}
\item[\textup{(i)}] The scalar kernel in this class is trivial; that is,
\[
 \Phi\in\mathscr X^0_{\mu,\nu}(\mathcal O_0),\quad
 (\Delta_X+pu_1^{p-1})\Phi=0
 \quad\Longrightarrow\quad\Phi=0.
\]
\item[\textup{(ii)}] For every finite $m\ge0$, the product kernel is trivial:
\[
 \Phi\in\mathscr X^m_{\mu,\nu}(\mathcal O_k),\quad
 (\Delta_Y+\Delta_X+pu_1^{p-1})\Phi=0
 \quad\Longrightarrow\quad\Phi=0.
\]
\end{enumerate}

\begin{remark}[Compatibility with the gluing weights]\label{rem:weight-compatibility}
For every pair in \eqref{e59-twoend}, the interval
\[
-a_p<\delta<
\min\{0,1-a_p,\operatorname{Re}\gamma_1^-,2-n-\mu\}
\]
is nonempty: each upper endpoint is strictly larger than $-a_p$.
One may choose $\delta$ and $\delta_*=2-n-\delta$ away from the
indicial spectrum. Then $\mu<\delta_*$ and \eqref{e59} holds.
Conversely, the requirements \eqref{e59} and $\mu<\delta_*$
imply the upper bound in \eqref{e59-twoend}. Thus the model and gluing
weights can be fixed simultaneously, including the restriction
$\delta<1-a_p$ used in the nonlinear error estimates.
For the gluing problem we take $q>N=n+k$. The local regularity in H
can be raised to any finite $q$ by elliptic estimates on half-annuli,
so this choice does not change its kernel class.
\end{remark}

\begin{remark}[The scaling Jacobi field]\label{rem:scaling-field}
The scaling field $Z_p=x\cdot\nabla u_1+a_pu_1$ satisfies $L_pZ_p=0$.
Scaled elliptic estimates differentiate the far-field expansion to obtain
$Z_p=(a_p+1-n+o(1))|x|^{-n}x_n$, with nonzero coefficient since $a_p<n-1$.
Thus fixed far-field normalization excludes $Z_p$. H\textup{(i)}, instead,
excludes it by the stronger inner weight: its decay $Z_p=O(r^{1-n})$ already satisfies the far-field bound
$|Z_p|\le Cr^\nu$ because $\nu>1-n$. Thus, whenever H\textup{(i)} holds,
the nonzero field $Z_p$ must fail the inner bound $|Z_p|\le Cr^\mu$.
At the threshold this is reflected in Proposition~\ref{prop:threshold-normalized-isomorphism}.

\end{remark}

\subsection{Massive normal operators and Fourier testing}\label{sec:massive-tools}
In this subsection we collect the tools which reduce H\textup{(ii)} to the
injectivity of the shifted normal operators $L_1-|\xi|^2$. They are used
again for the threshold branch in Section~\ref{sec:near-threshold-H},
for the homogeneous model in Section~\ref{s8}, and in
Appendix~\ref{app:fixed-regularity}.

We first record the elementary mode analysis of the conic models.
Let $V_0$ be either $p\phi_p^{p-1}$ or $0$ on $S^{n-1}_+$, let
$\widetilde\lambda_1<\widetilde\lambda_2\le\cdots$ be the Dirichlet
eigenvalues of $-\Delta_\theta-V_0$, and let $\gamma_j^\pm$ be the
corresponding indicial roots; for $V_0=0$ these are the harmonic roots
\eqref{eq:harmonic-indicial-roots}.

\begin{lemma}[Conic modes]\label{lem:conic-modes}
Let $c\ge0$, and let
$\Psi\in W^{2,q}_{\rm loc}(\overline{\R^n_+}\setminus\{0\})$ have zero
trace on $\partial\R^n_+\setminus\{0\}$ and solve
\[
 \bigl(\Delta+|X|^{-2}V_0(X/|X|)-c\bigr)\Psi=0
 \quad\text{in }\R^n_+.
\]
\begin{enumerate}[label=(\roman*)]
\item If $c=0$, $\sigma$ is nonindicial and $|\Psi|\le C|X|^\sigma$ in
$\R^n_+$, then $\Psi=0$.
\item If $c>0$, $\sigma>\operatorname{Re}\gamma_1^-$, and
$|\Psi|\le C|X|^\sigma$ for $|X|\le1$ and $|\Psi|\le C|X|^M$ for
$|X|\ge1$ and some $M$, then $\Psi=0$.
\end{enumerate}
\end{lemma}
\begin{proof}
Put $d=(n-2)/2$ and let $\psi_j$ be the eigenfunctions. The
coefficients $\Psi_j(r)=\int_{S^{n-1}_+}\Psi(r\theta)\psi_j\,d\sigma$
inherit the bounds of $\Psi$, and $\Psi_j=r^{-d}g_j$, where
\[
 g_j''+\frac1rg_j'-\Bigl(c+\frac{\nu_j^2}{r^2}\Bigr)g_j=0,\qquad
 \nu_j^2=d^2+\widetilde\lambda_j,\qquad
 \gamma_j^\pm=-d\pm\nu_j.
\]
If $c=0$, then $\Psi_j$ is a combination of $r^{\gamma_j^\pm}$, or of
$r^{-d}$ and $r^{-d}\log r$ when $\nu_j=0$. At zero the bound excludes
the terms whose exponent has real part below $\sigma$, and at infinity
those with real part above $\sigma$. Since $\sigma$ is nonindicial, no
term remains.

If $c>0$, then $g_j=AI_{\nu_j}(\sqrt c\,r)+BK_{\nu_j}(\sqrt c\,r)$.
The polynomial bound at infinity gives $A=0$. Near zero,
\[
 r^{-d}K_{\nu_j}(\sqrt c\,r)=
 \begin{cases}
 C r^{\gamma_j^-}(1+o(1)),&\nu_j>0,\\
 r^{-d}(-\log r+O(1)),&\nu_j=0,\\
 r^{-d}(C_+r^{\nu_j}+C_-r^{-\nu_j}+O(r^2)),&\nu_j\in i(0,\infty),
 \end{cases}
\]
with nonzero leading coefficients. Each behavior is incompatible with
$O(r^\sigma)$, since
$\operatorname{Re}\gamma_j^-\le\operatorname{Re}\gamma_1^-<\sigma$.
Thus $B=0$ in every mode, and $\Psi=0$.
\end{proof}

To deduce H\textup{(ii)} from shifted injectivity, we construct model inverses
whose smooth frequency dependence yields Fourier test functions. The
construction uses the scalar equations on the half-space. Throughout
this construction,
\begin{equation}\label{eq:scalar-model-potentials}
 B_\xi=\Delta_X-|\xi|^2+V(X),\qquad
 V\in\{pu_1^{p-1},\ p|X|^{-2}\phi_p^{p-1},\ 0\}.
\end{equation}
The first two potentials have the same indicial spectrum
$\mathfrak I_p$ at zero, whereas the zero potential has the harmonic
spectrum~\eqref{eq:harmonic-indicial-roots}. The positive mass requires the same
exterior weight for the solution and source. With a smooth cutoff
$\chi=1$ near zero, define
\begin{equation}\label{eq:massive-two-end-spaces}
\begin{aligned}
 \mathcal Y_{s,t}&=\{f:\chi f\in\mathcal H_{s-2},\ (1-\chi)f\in\mathcal H_t\},\\
 \mathcal D_{s,t}(\xi)&=\{v:\chi v\in\mathcal H_s,\ (1-\chi)v\in\mathcal H_t,
                    \ B_\xi v\in\mathcal Y_{s,t},\ v|_{\partial\R^n_+\setminus0}=0\}.
\end{aligned}
\end{equation}
Here, explicitly,
\begin{equation}\label{eq:normal-Hilbert-spaces}
 \mathcal H_b:=\left\{f:\int_{\R^n_+}r^{-2b-n}|f|^2\,dX<\infty\right\},
 \qquad
 \|f\|_{\mathcal H_b}^2=\int_{\R^n_+}r^{-2b-n}|f|^2\,dX.
\end{equation}
The source norm is the sum of its inner and outer norms, and
\[
 \|v\|_{\mathcal D_{s,t}(\xi)}
 :=\|\chi v\|_{\mathcal H_s}
   +\|(1-\chi)v\|_{\mathcal H_t}
   +\|B_\xi v\|_{\mathcal Y_{s,t}}.
\]
We use the closed distributional realization with zero Dirichlet trace
on the physical boundary away from zero. Local Dirichlet regularity
identifies it with $H^2$ functions on compact closed half-annuli.
Thus $\mathcal D_{s,t}(\xi)$ is a Banach space with this graph norm. In particular the inclusion $\mathcal D_{s,t}\hookrightarrow\mathcal Y_{s,t}$ is
bounded: at zero this follows from the factor $r^2$, and at infinity
both spaces use $t$. Thus the domains are independent of $\xi$ locally,
with equivalent norms, since
\[
 B_\xi-B_{\xi_0}=-(|\xi|^2-|\xi_0|^2)I.
\]

\begin{lemma}[Scalar realizations with positive mass]\label{lem:scalar-massive}
Let $B_\xi$ have one of the potentials in
\eqref{eq:scalar-model-potentials}. For the cell potential assume
$r^2V-p\phi_p^{p-1}=O(r^\sigma)$ at zero for some $\sigma>0$;
at infinity $V=O(r^{-2})$. For
$\xi\in K\Subset\R^k\setminus\{0\}$, the graph realization
\eqref{eq:massive-two-end-spaces} is Fredholm whenever $s$ avoids the
real parts of the inner indicial roots. Under the unweighted sesquilinear pairing
\[
 \langle f,h\rangle_0:=\int_{\R^n_+}f\overline h\,dX,
\]
the range annihilator, that is,
\[
 (\operatorname{Ran}B_\xi)^{\rm ann}
 :=\{\ell\in\mathcal Y_{s,t}^*:\ell(B_\xi v)=0
                     \text{ for all }v\in\mathcal D_{s,t}\},
\]
is identified with the kernel at complementary weights
\[
 (s,t)^\vee=(2-n-s,-n-t).
\]
This annihilator is the dual of the cokernel because the range is
closed. For every finite exterior weight $t$, the graph domain is
locally independent of $\xi$. The kernel is unchanged when $t$ varies
over finite values, or when $s$ varies within a connected component of
$\R\setminus\mathfrak R$, where $\mathfrak R$ is the set of real parts
of the relevant inner roots.
\end{lemma}
\begin{proof}
For the cells used below, a fixed dilation and Lemma~\ref{lem14}
place the inner profile in the neighborhood of the construction.
Remark~\ref{rem5} gives
$w-\phi_p=O_{W^{2,q}_{\rm ul}}(e^{\sigma't})$ for some $\sigma'>0$.
Since $q>n$, the H\"older or Lipschitz bound for the power $p-1$ yields
\[
 r^2V-p\phi_p^{p-1}=p(w^{p-1}-\phi_p^{p-1})=O(r^\sigma),
 \qquad \sigma=\sigma'\min\{1,p-1\}>0,
\]
and a fixed dilation preserves this rate.

\proofstep{1}{Construct the end inverses.}
Near zero, multiplication by $r^2$ and $\tau=\log r$ give the limiting
cylinder operator $P=\partial_\tau^2+(n-2)\partial_\tau-\mathcal A$,
where $\mathcal A=-\Delta_{S^{n-1}_+}-p\phi_p^{p-1}$ (zero potential
for the harmonic model). In angular mode $j$, $e^{-s\tau}Pe^{s\tau}$
has Fourier symbol
\[
 (i\lambda+s)^2+(n-2)(i\lambda+s)-\widetilde\lambda_j,
 \qquad \lambda\in\R.
\]
Here $\widetilde\lambda_j$ are the Dirichlet eigenvalues of $\mathcal A$.
The nonindicial condition makes these symbols invertible and angular
spectral decomposition gives the cylinder inverse with its graph
estimate, also for repeated or nonreal roots. The terms
$-|\xi|^2r^2$ and $r^2V-p\phi_p^{p-1}$ are small for sufficiently
small $r$.

At infinity, $U_tv=r^{-t-n/2}v$ identifies $\mathcal H_t$ with
$L^2(dX)$ and, for $a=t+n/2$,
\begin{equation}\label{eq:massive-polynomial-conjugation}
 U_t(\Delta-|\xi|^2)U_t^{-1}
 =\Delta-|\xi|^2+2a\frac{X}{r^2}\cdot\nabla
                +\frac{a(a+n-2)}{r^2}.
\end{equation}
The $O(r^{-1})$ first-order term and $O(r^{-2})$ potential are
absorbed by the positive-mass Dirichlet estimate outside a large ball,
uniformly for $\xi\in K$. Patching these end inverses with local
Dirichlet parametrices gives a two-sided parametrix: cutoff
commutators are compact by local Rellich compactness, and the
potential errors by the same argument and their small tails.
This proves Fredholmness with estimates uniform on compact nonzero
frequency sets; only boundedness of the angular potential is used.
Moreover, $B_\xi-B_{\xi_0}=|\xi_0|^2-|\xi|^2$ maps the domain into
$\mathcal Y_{s,t}$, so its graph domain is locally independent of $\xi$,
with locally equivalent graph norms.

\proofstep{2}{Identify the range annihilator.}
Under the unweighted integral, $\mathcal H_b^*=\mathcal H_{-n-b}$.
Applying this at the two ends of $\mathcal Y_{s,t}$ gives
\begin{equation}\label{eq:massive-green-dual}
 (s,t)^\vee=(2-n-s,-n-t).
\end{equation}
If $h\in\mathcal Y_{s,t}^*$ annihilates the range, interior tests give
$B_\xi h=0$. Zero-value tests with arbitrary normal derivative on the
physical boundary give its zero Dirichlet trace by transposition;
local boundary regularity then gives $H^2$ regularity there.
Together with the dual weights and $B_\xi h=0$, this puts
$h$ in $\ker(B_\xi:\mathcal D_{s^\vee,t^\vee}\to
\mathcal Y_{s^\vee,t^\vee})$.

Conversely, Green's formula extends to the graph domain by endpoint
cutoffs. Let $\eta_{\rho,R}=1$ on $2\rho<r<R$, vanish for
$r<\rho$ and $r>2R$, and set $A_d=\{d/2<r<4d\}\cap\R^n_+$.
Scaled Caccioppoli estimates at zero and unit-scale estimates at
infinity give, uniformly for $\xi\in K$,
\[
\begin{split}
 \|[B_\xi,\eta_{\rho,R}]v\|_{\mathcal Y_{s,t}}
 \le C\bigl(&\|v\|_{\mathcal H_s(A_\rho)}
       +\|B_\xi v\|_{\mathcal H_{s-2}(A_\rho)}\\
       &+R^{-1}\{\|v\|_{\mathcal H_t(A_R)}
                   +\|B_\xi v\|_{\mathcal H_t(A_R)}\}\bigr)
 \longrightarrow0.
\end{split}
\]
Here $[B_\xi,\eta]v=2\nabla\eta\cdot\nabla v+(\Delta\eta)v$;
the inner source weight $s-2$ compensates for both cutoff derivatives,
whereas the exterior source has weight $t$.
Thus $\eta_{\rho,R}v\to v$ in the graph norm. Dirichlet approximation
on compact half-annuli supplies a smooth graph core, on which Green's
formula pairs $B_\xi v$ to zero against the complementary kernel.
Continuity of the dual pairing extends this identity to the whole
domain, proving the claimed annihilator identification. Closed range
then identifies it with the dual of the cokernel.

Every kernel element has exponential decay at infinity. To justify
this from a polynomial exterior graph bound, put
$c_K=\inf_{\xi\in K}|\xi|>0$ and choose
$0<\gamma<\gamma_1<c_K$. Take $M$ large enough that
$r^{-M}v\in L^2(r>1)$, and let $\chi_R$ vanish for $r<R$ and equal one
for $r>2R$. We test the homogeneous equation against $q_T^2\overline v$, where
$q_T=\chi_R r^{-M}e^{\gamma_1\min\{r,T\}}$ is a truncated weight. Integration by parts gives
\[
 \int\bigl(|\nabla(q_Tv)|^2+(|\xi|^2-V)|q_Tv|^2\bigr)\,dX
       =\int|\nabla q_T|^2|v|^2\,dX.
\]
After absorbing the weight derivatives, for $R$ sufficiently large,
\[
 \|q_Tv\|_{H^1(r>2R)}^2
 \le C_{K,\gamma_1,M,R}\|v\|_{L^2(R<r<2R)}^2,
 \qquad T>2R.
\]
Here $|\nabla\log(r^{-M}e^{\gamma_1\min\{r,T\}})|
\le\gamma_1+M/r$ and $V=O(r^{-2})$, so the positive gap
$|\xi|^2-V-(\gamma_1+M/r)^2$ absorbs the weight terms. The estimate
is first obtained with an additional outer cutoff, which can be removed
using the polynomial graph bound. Letting $T\to\infty$ and using
$r^M e^{-(\gamma_1-\gamma)r}\le C$ proves $e^{\gamma r}v\in L^2(r>2R)$.
Local Dirichlet estimates yield pointwise exponential decay and the
corresponding first-derivative estimates. Related positive-solution
comparison arguments appear in \cite[Lemmas~\textup{4.2} and~4.5]{ios}.
Thus changing the finite polynomial exterior weight does not change
the kernel. At the inner end kernels are unchanged by
moving the weight inside a root-free interval, using the cylinder
inverse above and the small scaled perturbation.

\end{proof}

We record the common Gram correction and Fourier testing arguments.
\begin{lemma}[Gram correction]\label{lem:gram-correction}
Let $D_{\rm lo},D_{\rm hi},Y_{\rm lo},Y_{\rm hi}$ be fixed complex
Banach spaces. Suppose that the bounded operators
\[
 B_\xi^{\rm lo}:D_{\rm lo}\to Y_{\rm lo},\qquad
 B_\xi^{\rm hi}:D_{\rm hi}\to Y_{\rm hi}
\]
are Fredholm, the first is surjective, and the second is injective.
Set $K_\xi=\ker B_\xi^{\rm lo}$. Assume that a continuous
sesquilinear pairing $\langle\cdot,\cdot\rangle_
\xi:Y_{\rm hi}\times K_\xi\to\C$, linear in its first argument,
identifies $K_\xi$ with the range annihilator of $B_\xi^{\rm hi}$;
explicitly, the conjugate-linear map
\[
 J_\xi\psi:f\longmapsto\langle f,\psi\rangle_\xi
\]
is a bijection from $K_\xi$ to
\[
 \{\ell\in Y_{\rm hi}^*:\ell(B_\xi^{\rm hi}y)=0
                                  \text{ for every }y\in D_{\rm hi}\}.
\]
Let $C_\xi:D_{\rm lo}\to Y_{\rm hi}$ be bounded and assume that
$G_\xi(\kappa,\psi)=\langle C_\xi\kappa,\psi\rangle_\xi$ is a positive
definite Hermitian form on $K_\xi$. Then
\begin{equation}\label{eq:gram-corrected-operator}
 \mathcal T_\xi:D_{\rm lo}\times D_{\rm hi}\longrightarrow
                  Y_{\rm lo}\times Y_{\rm hi},\qquad
 \mathcal T_\xi(x,y)=(B_\xi^{\rm lo}x,\ B_\xi^{\rm hi}y-C_\xi x)
\end{equation}
is a bounded isomorphism. If the three operators depend smoothly on
$\xi$ in operator norm on an open frequency set, then so does
$\mathcal T_\xi^{-1}$, and every finite-order derivative of the inverse
is bounded on compact subsets. Over real Banach spaces, replace the
sesquilinear pairing and Hermitian form by a bilinear pairing and a
symmetric positive definite form.
\end{lemma}
\begin{proof}
Given $(f,g)$, first fix $x_0\in D_{\rm lo}$ satisfying
$B_\xi^{\rm lo}x_0=f$; such a solution exists by surjectivity. Then
\[
 B_\xi^{\rm lo}x=f
 \iff B_\xi^{\rm lo}(x-x_0)=0
 \iff x\in x_0+K_\xi,\qquad K_\xi=\ker B_\xi^{\rm lo}.
\]
Thus all solutions of the first equation, and only those solutions,
are $x=x_0+\kappa$ with $\kappa\in K_\xi$.
Since the high range is closed, the
annihilator hypothesis makes the second equation solvable precisely when
\[
 \langle C_\xi\kappa,\psi\rangle_\xi
 =-\langle g+C_\xi x_0,\psi\rangle_\xi
 \qquad(\psi\in K_\xi).
\]
The finite-dimensional positive definite Gram form determines
$\kappa$ uniquely, and high injectivity determines $y$ uniquely.
Thus $\mathcal T_\xi$ is bijective, and the bounded inverse theorem
applies. In particular the final pair is independent of the choice of
$x_0$. Smooth dependence follows from inversion on the fixed product
spaces, or from
\[
 \mathcal T_\xi^{-1}-\mathcal T_{\xi_0}^{-1}
 =-\mathcal T_\xi^{-1}(\mathcal T_\xi-\mathcal T_{\xi_0})
                      \mathcal T_{\xi_0}^{-1}.
\]
No smoothly varying basis of $K_\xi$ is needed.
\end{proof}

We next formulate the Fourier argument for the scalar and triangular
systems used below. Use $\mathcal O_k$ from Assumption~H, with
$k\ge1$, and fix $q>n+k$.
For vectors of real exponents $\mathbf b=(b_i)$, $\mathbf t=(t_i)$,
define $\mathscr X^m_{\mathbf b,\mathbf t}(\mathcal O_k)$ to consist of
$U\in W^{2,q}_{\rm loc}(\mathcal O_k)$ with zero physical Dirichlet
trace and finite norm
\begin{equation}\label{eq:fourier-solution-class}
\begin{split}
 \|U\|_{\mathscr X^m_{\mathbf b,\mathbf t}}
 :=\sum_i\bigg[&\sup_{Y,\,0<r\le1}
  \langle Y\rangle^{-m}r^{-b_i}
       (|U_i|+r|\nabla_XU_i|)\\
 &+\sup_{Y,\,r\ge1}\langle Y\rangle^{-m}r^{-t_i}
       (|U_i|+|\nabla_XU_i|)\bigg].
\end{split}
\end{equation}
The first-derivative bounds in this class can also be verified from
local elliptic estimates in each application. In particular, the
exterior exponents $t_i$ may be enlarged to accommodate these bounds.
Only a finite polynomial exterior bound is used in the next lemma.

\begin{lemma}[Fourier testing]\label{lem:fourier-testing}
Let $\mathcal L=\Delta_Y+\Delta_X+V(X)$, where $V$ is real, locally
bounded up to the flat boundary away from zero, $V=O(r^{-2})$ at zero,
and $V$ has at most polynomial growth at infinity. Consider either
the scalar equation $\mathcal LU=0$, or the system
\begin{equation}\label{eq:fourier-triangular-system}
 \mathcal LU_1=0,\qquad \mathcal LU_2-W(X)U_1=0,
\end{equation}
where $W$ is real, locally bounded away from zero,
$|W|\le Cr^{-d_c}$ for $r\le1$, and has at most polynomial growth at
infinity. Suppose that $U\in\mathscr X^m_{\mathbf b,\mathbf t}$ for
some finite $m\ge0$, and that the equations admit the usual local
Dirichlet elliptic estimates on closed normal half-annuli.

Fix the unitary Fourier transform in $Y$, with kernel
$(2\pi)^{-k/2}e^{-iY\cdot\xi}$, and put
$B_\xi=\Delta_X-|\xi|^2+V(X)$. For every smooth data vector
$F\in C_c^\infty((\R^k\setminus\{0\})\times\R^n_+)$, assume the existence of a zero-trace
test solution, smooth and compactly supported in $\xi$ away from zero,
of
\[
\begin{array}{ll}
 B_\xi\Psi=F,&\text{in the scalar case},\\
 B_\xi\Psi_1-W\Psi_2=F_1,\quad B_\xi\Psi_2=F_2,
 &\text{in the paired case}.
\end{array}
\]
Assume there are inner exponents $c_i$ such that, for every
multi-index $\alpha$, uniformly on the compact frequency support,
\begin{equation}\label{eq:fourier-test-estimates}
\begin{aligned}
 |\partial_\xi^\alpha\Psi_i|+
       r|\nabla_X\partial_\xi^\alpha\Psi_i|
   &\le C_\alpha r^{c_i},&&0<r\le1,\\
 |\partial_\xi^\alpha\Psi_i|+
       |\nabla_X\partial_\xi^\alpha\Psi_i|
   &\le C_\alpha e^{-\gamma r},&&r\ge1
\end{aligned}
\end{equation}
for some $\gamma>0$. Require
\begin{equation}\label{eq:fourier-inner-exponent-conditions}
\begin{aligned}
 b_i+c_i+n-2&>0&&\text{for each component},\\
 b_1+c_2+n-d_c&>0&&\text{in the paired case}.
\end{aligned}
\end{equation}
Then the partial Fourier transform satisfies
$\operatorname{supp}_\xi\widehat U\subset\{0\}$. Consequently,
\[
 U(Y,X)=\sum_{|\alpha|\le\lfloor m\rfloor}Y^\alpha U_\alpha(X)
\]
with distributional normal coefficients. In particular, a solution
uniformly bounded in $Y$ is independent of $Y$.
\end{lemma}
\begin{proof}
Let $F$ be arbitrary as in the statement and let $\Psi$ solve its
transposed normal problem. The inverse transforms
$f=F^\vee$ and $v=\Psi^\vee$ satisfy the transposed physical-space
equation. Integration by parts in $\xi$ provides, for every $N$, bounds
of the form
\[
 |v_i|+r|\nabla_Xv_i|
 \le C_N\langle Y\rangle^{-N}r^{c_i}\quad(r\le1),
\]
and corresponding bounds with $e^{-\gamma r}$ for $r\ge1$.
First insert a tangential cutoff and integrate on
$\{\rho<|X|<R,\ X_n>0\}$. With $\rho,R$ fixed, the local elliptic
estimates and $N>m+k+2$ allow removal of the tangential cutoff.
Indeed, with $\chi_T(Y)=\chi(Y/T)$, every $Y$-derivative of the cutoff
falls on the test function, since
$\mathcal L(\chi_Tv)=\chi_T\mathcal Lv+2\nabla_Y\chi_T\cdot\nabla_Yv
+(\Delta_Y\chi_T)v$; the inverse transforms $\nabla_Yv$ have the same
rapid decay in $Y$, so no bound on $\nabla_YU$ is needed.
Both zero traces eliminate the flat-boundary terms. The remaining
inner Green terms are bounded by
\[
 C\sum_i\rho^{b_i+c_i+n-2}
       \int_{\R^k}\langle Y\rangle^{m-N}\,dY=o(1).
\]
The outer terms vanish by exponential test decay and the polynomial
solution bounds. In the paired case the coupling integral at zero is
absolutely convergent because
\[
 \int_{r<1}|WU_1\overline{v_2}|\,dX\,dY
 \le C\int_0^1r^{b_1+c_2-d_c+n-1}\,dr<\infty.
\]
Hence Green's formula, followed by $\rho\downarrow0$ and
$R\to\infty$, gives
\[
 \langle\widehat U,F\rangle
 =\int_{\mathcal O_k}U\cdot\overline f\,dX\,dY=0.
\]
Here the paired volume terms cancel using
$\mathcal LU_1=0$ and $\mathcal LU_2=WU_1$. Thus $\widehat U$ vanishes
against every test away from zero frequency. Pairing with a compactly
supported normal test and applying the point-support theorem
\cite[Theorem~\textup{2.3.4}]{hormander} produces a polynomial in $Y$.
The common bound $\langle Y\rangle^m$ limits its degree to
$\lfloor m\rfloor$, uniformly in the normal test, and hence defines
the distributional coefficients $U_\alpha$. If $m=0$, only the constant
coefficient remains.
\end{proof}

\begin{lemma}[From shifted injectivity to H\textup{(ii)}]\label{lem:shift-to-product}
Let $u_1$ be a fast-decaying cell with the singular-end estimate in
Lemma~\ref{lem:scalar-massive}. Choose
\[
 \operatorname{Re}\gamma_1^+<\mu_0<\mu<2-n+a_p,\qquad
 1-n<\nu<-\frac{n-2}{2},
\]
with $\mu_0,\mu$ in the first root-free gap.
Suppose that, for every $\tau\ge0$,
\[
 v\in\mathscr X^0_{\mu_0,\nu}(\mathcal O_0),\quad
 (L_1-\tau)v=0\quad\Longrightarrow\quad v=0.
\]
Then H\textup{(i)} and H\textup{(ii)} hold with weights $(\mu,\nu)$, for every $k\ge1$.
\end{lemma}
\begin{proof}
\proofstep{1}{Construct test functions at nonzero frequencies.}
H\textup{(i)} follows by inclusion of the pointwise classes and $\tau=0$.
Fix $K\Subset\R^k\setminus\{0\}$ and set
\[
 B_\xi=L_1-|\xi|^2,\qquad \check\mu=2-n-\mu_0,\qquad t=-n/2.
\]
Write $\mathcal D_\sigma=\mathcal D_{\sigma,t}$ and
$\mathcal Y_\sigma=\mathcal Y_{\sigma,t}$.
Lemma~\ref{lem:scalar-massive} provides complementary high and low
realizations with inner weights $\mu_0$ and $\check\mu$.
A high kernel element has the pointwise bound $O(r^{\mu_0})$ by
dyadic elliptic estimates and decays exponentially at infinity, so
shifted injectivity makes it zero. The low realization is therefore
surjective.

Choose $0\le\chi\in C_c^\infty(\R^n_+\setminus\{0\})$ positive on an
open set. The triangular operator
\[
 \mathcal C_\xi=
 \begin{pmatrix}B_\xi&0\\-\chi&B_\xi\end{pmatrix}:
 \mathcal D_{\check\mu}\oplus\mathcal D_{\mu_0}
 \longrightarrow\mathcal Y_{\check\mu}\oplus\mathcal Y_{\mu_0}
\]
is an isomorphism by Lemma~\ref{lem:gram-correction}: Green duality
identifies the required kernel and annihilator, and
$\int\chi\kappa\overline\psi$ is positive definite by unique continuation.
The domains are locally fixed by Lemma~\ref{lem:scalar-massive}, so
the inverse is smooth on $K$.

For $\psi\in C_c^\infty(K\times(\R^n_+\setminus\{0\}))$, put
$(\eta,\zeta)=\mathcal C_\xi^{-1}(\psi,0)$.
Dyadic elliptic estimates, followed by finite local bootstrapping, give
\[
 |\partial_\xi^\alpha\eta|+
 r|\nabla_X\partial_\xi^\alpha\eta|
 \le C_\alpha r^{\check\mu},\qquad r\le1.
\]
At infinity these functions and their first derivatives decay
exponentially. Differentiating the equations adds only lower-order
sources in the same spaces, so induction establishes these bounds at every
finite order. The inner weights remain in fixed root-free gaps;
neither the number of derivatives nor the dimension $k$ requires
shrinking a previously fixed exponent interval.

\proofstep{2}{Show that the Fourier transform is supported at zero.}
Let $\Phi$ satisfy H\textup{(ii)}. The tests above have rapidly decreasing
inverse Fourier transforms, and their first component satisfies
$(\Delta_Y+L_1)\eta^\vee=\psi^\vee$. The inner Green exponent is
\begin{equation}\label{eq:common-frequency-boundary}
 \mu+\check\mu+n-2=\mu-\mu_0>0.
\end{equation}
By local elliptic estimates, the first scaled derivatives have the
same weights. Exponential test decay removes the outer terms.
Lemma~\ref{lem:fourier-testing} therefore implies
$\operatorname{supp}_\xi\widehat\Phi\subset\{0\}$.

\proofstep{3}{Show that all polynomial coefficients vanish.}
Lemma~\ref{lem:fourier-testing} and the polynomial growth imply

\begin{equation}\label{eq:zero-frequency-polynomial}
\Phi(Y,X)=\sum_{|\alpha|\le M_0}Y^\alpha a_\alpha(X),
\end{equation}
for some finite \(M_0\), initially with distributional coefficients in
\(X\).  Substitution into \((\Delta_Y+L_1)\Phi=0\) and local elliptic
regularity show that the coefficients are \(W^{2,q}_{\rm loc}\), and the
flat-boundary condition holds coefficientwise.  Evaluating the polynomial
on a fixed finite set of \(Y\)-points and using finite-dimensional norm
equivalence yields
\[
|a_\alpha(X)|
\le
C\omega_{\mu,\nu}(|X|).
\]
Since \(\Delta_Y\) lowers polynomial degree by two, every coefficient of
maximal \(Y\)-degree satisfies
\[
L_1a_\alpha=0
\]
and hence vanishes by H\textup{(i)}.  Descending in polynomial degree proves \(\Phi\equiv0\).
\end{proof}

\subsection{Non-degeneracy in the stable range}\label{sec:stable-nondegeneracy}
Let $u_1$ be the normalized cell constructed from the minimal branch
in Section~\ref{sec:stable-construction}. Semistability controls the quadratic form in
the shifted equation; comparison with the homogeneous profile will
provide the strict coercivity needed at zero frequency. The cell satisfies
\begin{equation}\label{eq:fast-cell-semistability}
Q_{u_1}(\zeta):=
\int_{\R^n_+}\left(
|\nabla\zeta|^2-pu_1^{p-1}\zeta^2
\right)\,dX\ge0
\qquad
\text{for every }\zeta\in C_c^\infty(\R^n_+).
\end{equation}
Section~\ref{sec:stable-construction} proves semistability of $V_\lambda$ on
$H_0^1(B_1^{n+2})$ by a positive-function argument. Lifting and blow-up
rescaling preserve the quadratic form up to a positive constant.
For $v=x_n\bar v$, one uses
\[
\int_{\R^n_+}|\nabla(x_n\xi)|^2
=
\int_{\R^n_+}x_n^2|\nabla\xi|^2
\]
for compactly supported $\xi$, while the potential terms transform in the
same way.  Finally, the Kelvin transform preserves the Dirichlet form and,
because $\kappa_p=p(n-2)-(n+2)$, also the linearized potential term.

The same construction preserves the order bound. With the notation of
Section~\ref{sec:stable-construction}, $0<\bar v_\lambda\le\lambda\bar v_\infty$ and
$R_j=m_j^{1/\sigma_p}$, where $\sigma_p=n-1-2/(p-1)$. Homogeneity gives
\[
\bar v_j(x)=m_j^{-1}\bar v_{\lambda_j}(x/R_j)
\le\lambda_j m_j^{-1}R_j^{\sigma_p}\bar v_\infty(x)
=\lambda_j\bar v_\infty(x).
\]
Passing to the limit, multiplying by $x_n$, and applying the Kelvin
transform yield $0<u_1\le u_\infty$. The final normalization by dilation
preserves this bound. Since $\Lambda_1(p)>0$, \eqref{e10}--\eqref{e12} imply
\[
Q_{u_1}(\zeta)\ge Q_{u_\infty}(\zeta)
\ge\Lambda_1(p)\int_{\R^n_+}|X|^{-2}\zeta^2\,dX,
\qquad \zeta\in C_c^\infty(\R^n_+).
\]
Since
$u_1^{p-1}\le u_\infty^{p-1}\le C|X|^{-2}$, the same estimate implies
\[
\int_{\R^n_+}|\nabla\zeta|^2\,dX
=Q_{u_1}(\zeta)+p\int_{\R^n_+}u_1^{p-1}\zeta^2\,dX
\le C_p Q_{u_1}(\zeta).
\]
Thus $Q_{u_1}$ controls both the homogeneous $H^1_0$ norm and the
$L^2(|X|^{-2}dX)$ norm.

\begin{lemma}\label{lem:stable-ground-state}
Assume
$\frac{n+1}{n-1}<p<p^*(n)$,
and let $u_1$ be the normalized fast-decaying cell obtained from the minimal branch in Section~\ref{sec:stable-construction}. Then Assumption~H\textup{(i)} holds for every pair $(\mu,\nu)$ satisfying \eqref{e59-twoend}.
\end{lemma}

\begin{proof}
Let $\Phi$ satisfy the equation, boundary condition, and bounds in H\textup{(i)}.
Choose a radial cutoff $\eta_{\rho,R}$ with
\[
\eta_{\rho,R}=1\quad\text{on }\{2\rho<|X|<R\},
\qquad
|\nabla\eta_{\rho,R}|\le C|X|^{-1}.
\]
Testing $L_1\Phi=0$ against $\eta_{\rho,R}^2\Phi$ gives
\begin{equation}\label{eq:stable-cutoff-identity}
Q_{u_1}(\eta_{\rho,R}\Phi)
=
\int_{\R^n_+}\Phi^2|\nabla\eta_{\rho,R}|^2\,dX.
\end{equation}
The two cutoff errors satisfy
\[
\int_{\R^n_+}\Phi^2|\nabla\eta_{\rho,R}|^2\,dX
\le C\bigl(\rho^{n+2\mu-2}+R^{n+2\nu-2}\bigr)\longrightarrow0,
\]
since $\mu>\gamma_1^+>-(n-2)/2$ and $\nu<-(n-2)/2$.
For any fixed compact half-annulus $A$, strict coercivity and
\eqref{eq:stable-cutoff-identity} therefore imply
\[
\Lambda_1(p)\int_A|X|^{-2}\Phi^2\,dX
\le Q_{u_1}(\eta_{\rho,R}\Phi)\longrightarrow0.
\]
Thus $\Phi=0$ on every such $A$.
\end{proof}

For the stable cell, the shifted hypothesis of Lemma~\ref{lem:shift-to-product} follows from
the same cutoff identity used at zero frequency.

\begin{lemma}\label{lem:stable-edge-nondeg}
Assume
$\frac{n+1}{n-1}<p<p^*(n)$,
and let $u_1$ be the normalized fast-decaying cell obtained from the minimal branch in Section~\ref{sec:stable-construction}. Then, for every $k\ge1$, Assumption~H\textup{(ii)} holds for every pair $(\mu,\nu)$ satisfying \eqref{e59-twoend}.
\end{lemma}

\begin{proof}
Let $\tau>0$.  If $(L_1-\tau)\Psi=0$ with the two-ended bounds in Assumption~H\textup{(i)}, then
\begin{equation}\label{eq:stable-shifted-cutoff}
Q_{u_1}(\eta_{\rho,R}\Psi)
+\tau\int_{\R^n_+}\eta_{\rho,R}^2\Psi^2\,dX
=
\int_{\R^n_+}\Psi^2|\nabla\eta_{\rho,R}|^2\,dX.
\end{equation}
By Lemma~\ref{lem:stable-ground-state}, the right-hand side tends to
zero. Since $Q_{u_1}\ge0$, this implies
\begin{equation}\label{eq:shifted-normal-injectivity-stable}
\ker(L_1-\tau)=\{0\},\qquad \tau>0.
\end{equation}

Choose $\mu_0<\mu$ in the same first root-free gap. The energy
argument above applies with $\mu_0$ in place of $\mu$, since
$\mu_0>\gamma_1^+>-(n-2)/2$. The case $\tau=0$ is
Lemma~\ref{lem:stable-ground-state}. The singular-end remainder was
verified in Lemma~\ref{lem:scalar-massive}, so
Lemma~\ref{lem:shift-to-product} implies H\textup{(ii)}.
\end{proof}

\subsection{Non-degeneracy at and near the stability threshold}
\label{sec:near-threshold-H}

At the stability threshold the two roots of the first angular mode coalesce and the homogeneous
coercivity constant $\Lambda_1(p)$ vanishes. The cylinder energy identity still proves H\textup{(i)},
but constructing a nearby branch also requires control of the repeated
mode. Once this branch is obtained, a scaling argument proves shifted
injectivity and Lemma~\ref{lem:shift-to-product} implies H\textup{(ii)}. Assume
\[
p_*:=p^*(n)<p_S(n):=\frac{n+2}{n-2},
\qquad
d:=\frac{n-2}{2}.
\]
Let \(u_*\) denote the coefficient-\(1\) threshold cell obtained in
Section~\ref{s2} as the stable-side limit \(p_j\uparrow p_*\).  Thus
\[
u_*(x)=(1+o(1))|x|^{-n}x_n
\qquad (|x|\to\infty).
\]
At \(p=p_*\), the two indicial roots of the first angular mode coalesce at
\[
\gamma_1^+=\gamma_1^-=-d.
\]
Fix once and for all an admissible pair \((\mu,\nu)\) for
\eqref{e59-twoend} at \(p=p_*\), chosen so that
\begin{equation}\label{eq:near-threshold-weights}
-d<\mu,
\qquad
1-n<\nu<-d,
\end{equation}
and \(\mu\) lies in the first root-free gap above \(-d\).  Since the
inequalities in \eqref{e59-twoend} are strict, after shrinking the
\(p\)-interval the same pair remains admissible for every nearby
parameter.

\subsubsection{Non-degeneracy at the threshold}
\begin{lemma}\label{lem:threshold-semistable-H}
The threshold cell \(u_*\) is semistable:
\[
Q_*(\zeta):=
\int_{\R^n_+}
\left(
|\nabla\zeta|^2-p_*u_*^{p_*-1}\zeta^2
\right)\,dx\ge0
\qquad
(\zeta\in C_c^\infty(\R^n_+)).
\]
For every admissible pair \((\mu,\nu)\) satisfying
\eqref{eq:near-threshold-weights}, with \(\mu\) in the first root-free gap
above \(-d\), the cell \(u_*\) satisfies Assumption~H\textup{(i)}.
\end{lemma}

\begin{proof}
The order bound and semistability in Section~\ref{sec:stable-nondegeneracy}
are dilation invariant and pass to the threshold limit, giving
\[
0<u_*\le u_{\infty,*},
\qquad
Q_*\ge0.
\]

For H\textup{(i)}, write
\[
a_*=\frac2{p_*-1},\qquad
w_*(\theta,t)=e^{a_*t}u_*(e^t\theta),\qquad
y(\theta,t)=e^{dt}\Phi(e^t\theta).
\]
Since \(p_*<p_S(n)\), Lemma~\ref{lem:profile-analytic} and the definition
of \(p_*\) imply
\[
\mathcal A_*:=-\Delta_\theta+d^2-p_*\phi_*^{p_*-1}\ge0.
\]
Since \(w_*\le\phi_*\),
\[
\mathcal V_*:=p_*(\phi_*^{p_*-1}-w_*^{p_*-1})\ge0,
\]
and \(L_*\Phi=0\) becomes
\[
-y_{tt}+\mathcal A_*y+\mathcal V_*y=0.
\]
The admissible weights ensure exponential decay of \(y\) and one
derivative at both cylinder ends.  Testing on
\((-T,T)\times S^{n-1}_+\) and letting \(T\to\infty\) yields
\[
\int_\R\left(
\|y_t(\cdot,t)\|_{L^2(S^{n-1}_+)}^2
+\langle\mathcal A_*y(\cdot,t),y(\cdot,t)\rangle_{L^2(S^{n-1}_+)}
+\int_{S^{n-1}_+}\mathcal V_*y^2\,d\sigma
\right)\,dt=0.
\]
Thus \(y_t=0\), and the end decay implies \(y\equiv0\).
\end{proof}

\subsubsection{Invertibility of the normalized linearization}
The preceding lemma excludes Jacobi fields with inner weight $\mu>-d$.
The matching linearization allows the repeated mode as well, so its
kernel must be examined before applying the implicit function theorem.
To place the threshold cell in the normalized analytic formulation,
let \(c>0\) be given by Lemma~\ref{lem14} and choose the unique
\(\rho_*>0\) such that
\[
\rho_*^{a_*+1-n}=c,\qquad
\bar u_*(x):=\rho_*^{a_*}u_*(\rho_*x).
\]
Then \(\bar u_*(x)=(c+o(1))|x|^{-n}x_n\) at infinity.
Lemma~\ref{lem14} supplies the required inner \(C^2\)-smallness, and
\eqref{eq:first-mode-derivative-bound} controls the first-mode derivative.  Positivity,
the Hopf lemma, and the far-field normalization therefore ensure
\[
(p_*,\bar u_*)\in\mathcal V_{\eta_{{\rm an},0}}.
\]
Dilation preserves Lemma~\ref{lem:threshold-semistable-H}.

\begin{proposition}\label{prop:threshold-normalized-isomorphism}
The normalized linearization at \((p_*,\bar u_*)\) is an isomorphism:
\[
\partial_uA[p_*,\bar u_*]:
F\longrightarrow\mathcal E_{p_*}.
\]
Hence there are \(\eta_{\rm th,1}>0\) and a real-analytic
coefficient-\(c\) branch
\[
p\longmapsto \bar u_p,
\qquad
|p-p_*|<\eta_{\rm th,1},
\qquad
\bar u_{p_*}=\bar u_*.
\]
After the natural dilation which changes the far coefficient from
\(c\) to \(1\), the coefficient-\(1\) branch is
\[
p\longmapsto u_p^{\rm th},
\qquad
u_{p_*}^{\rm th}=u_*,
\]
continuous in the fixed inner Emden--Fowler coordinates and the exterior
graph space $W$, and locally in \(C^2\).
\end{proposition}

\begin{proof}
\proofstep{1}{Derive the expansion at the repeated root.}
By Proposition~\ref{prop10}, the linearization is Fredholm of index
zero, so it remains to exclude its kernel. Lemma~\ref{lem:linear-matching}
extends a kernel element to a global Jacobi field $v$
with zero flat-boundary data.  Since \(v\in F=\ker I\),
\begin{equation}\label{eq:threshold-kernel-far}
v=o(|x|^{-n}x_n)
\qquad(|x|\to\infty).
\end{equation}

Set
\[
\alpha_*=2a_*-(n-2),
\qquad
\mu_*=\frac{\alpha_*}{2}=a_*-d>0.
\]
At \(p=p_*\) the two characteristic roots of the first mode coincide at \(\mu_*\).
The Green formulas in Lemma~\ref{lem9} remain valid for every
\(0<\sigma<\mu_*\): the first mode contains
\(e^{\mu_*t}\) and \(te^{\mu_*t}\), while the higher modes are separated
from \(\mu_*\) by Proposition~\ref{prop7}.  Therefore
\[
\bar w_*-\phi_*=O_{W^{2,q}_{\rm ul}}(e^{\sigma t}),
\]
and the same estimate holds for every kernel variation.

Choose $\sigma<\mu_*$ close enough to $\mu_*$ and
$0<\sigma_1<\min\{1,p_*-1\}\sigma$ so that
$\kappa:=\sigma+\sigma_1>\mu_*$. The boundary-factor estimate used
in \eqref{eq:actual-inner-source} implies
$p_*(\bar w_*^{p_*-1}-\phi_*^{p_*-1})=O(e^{\sigma_1t})$.
For $V(\theta,t)=e^{a_*t}v(e^t\theta)$, the cylindrical source is
\[
 f=-p_*(\bar w_*^{p_*-1}-\phi_*^{p_*-1})V
   =O_{L^q_{\rm ul}}(e^{\kappa t}).
\]
In particular, the source decays faster than the repeated-root exponent.
For $v_1=c_{p_*}^1[V]$ and $f_1=c_{p_*}^1[f]$,
\(\lambda_{p_*}^1=-\alpha_*^2/4\) and \(\mu_*=\alpha_*/2\), hence
\[
\Bigl(\frac{d}{dt}-\mu_*\Bigr)^2v_1=f_1.
\]
Since \(e^{-\mu_*t}f_1\) is integrable at \(-\infty\), variation of
constants gives
\[
v_1(t)=e^{\mu_*t}\left(
A_v+B_vt+
\int_{-\infty}^t(t-s)e^{-\mu_*s}f_1(s)\,ds
\right),
\]
and the contribution of the integral to \(v_1\) and to \(v_1'\) is
\(O(e^{\kappa t})\). For $j\ge2$, Proposition~\ref{prop7} implies
$\mu_j^+\ge\alpha_*+1>\mu_*$. Applying the Green formulas of
Lemma~\ref{lem9} with the spectral summation in \eqref{e116a} yields
\[
\begin{gathered}
 P_{p_*}V=V-c_{p_*}^1[V]\psi_{p_*}^1
 =O_{W^{2,q}_{\rm ul}}(e^{(\mu_*+\eta)t}),\\
 0<\eta<\min\{\kappa-\mu_*,\mu_2^+-\mu_*\}.
\end{gathered}
\]
Since \(q>n\), the local Sobolev embedding gives
\begin{equation}\label{eq:threshold-repeated-expansion}
V(\theta,t)
=
e^{\mu_*t}(A_v+B_vt)\psi_{p_*}^1(\theta)
+
O_{C^1}\!\left(e^{(\mu_*+\eta)t}\right).
\end{equation}

\proofstep{2}{Compare with the scaling field by a Wronskian identity.}
The same expansion holds for the scaling Jacobi field
\[
Z_*:=x\cdot\nabla\bar u_*+a_*\bar u_*,
\]
with coefficients \((A_Z,B_Z)\).  Since
\[
Z_*(x)
=
\bigl(c(a_*+1-n)+o(1)\bigr)|x|^{-n}x_n
\quad(|x|\to\infty),
\]
one must have \((A_Z,B_Z)\ne(0,0)\); otherwise the improved inner decay
from \eqref{eq:threshold-repeated-expansion} and
Lemma~\ref{lem:threshold-semistable-H} would force \(Z_*=0\).

Consider the Wronskian
\[
\mathcal W(r)
=
\int_{S_r^+}
\left(
v\,\partial_rZ_*-Z_*\,\partial_rv
\right)\,dS.
\]
Here $S_r^+=\{|x|=r,\ x_n>0\}$. It is constant by Green's identity.  From
\eqref{eq:threshold-kernel-far} and the exterior elliptic estimate,
\[
v=o(r^{1-n}),\qquad \nabla v=o(r^{-n}),
\]
so \(\mathcal W(r)\to0\) as \(r\to\infty\).  As \(r\downarrow0\),
\eqref{eq:threshold-repeated-expansion} implies
\[
\mathcal W(r)
=
(A_vB_Z-B_vA_Z)
\int_{S^{n-1}_+}(\psi_{p_*}^1)^2\,d\sigma
+o(1).
\]
Hence
\[
(A_v,B_v)=\lambda(A_Z,B_Z)
\]
for some \(\lambda\in\R\).  Thus
\[
\Phi:=v-\lambda Z_*
\]
has improved singular decay
\[
|\Phi|+|x||\nabla\Phi|
\le C|x|^{-d+\eta_1}
\]
for some \(\eta_1>0\), and \(\Phi=O(|x|^{1-n})\) at infinity.

\proofstep{3}{Exclude the scaling field by the far-field normalization.}
Choosing an admissible
\[
-d<\widetilde\mu<-d+\eta_1
\]
in the same root-free gap, Lemma~\ref{lem:threshold-semistable-H}
implies $\Phi=0$, hence $v=\lambda Z_*$. The zero far-field coefficient
of $v$ now gives
\[
0=\lambda c(a_*+1-n),\qquad c(a_*+1-n)\ne0.
\]
Thus $\lambda=0$ and $v=0$. Since the kernel is trivial and the index is zero, the operator is
surjective.

Remark~\ref{rem6} identifies domain and target with fixed Banach
spaces, so the analytic implicit function theorem provides the
coefficient-\(c\) branch \(p\mapsto\bar u_p\).  Finally set
\[
a_p=\frac2{p-1},\qquad
\rho(p)=c^{-1/(a_p+1-n)},\qquad
u_p^{\rm th}(x)=\rho(p)^{a_p}\bar u_p(\rho(p)x).
\]
Remark~\ref{rem:normalization-continuity} ensures continuity in the stated
fixed coordinates. Local elliptic estimates improve it to local $C^2$
convergence. After shrinking the interval, $\rho(p)$ and its inverse
are bounded, giving uniform two-ended estimates. Dilation conjugates
the shifted equations with $\tau$ multiplied by a positive scalar.
\end{proof}

\subsubsection{Injectivity of the shifted operator}
Lemma~\ref{lem:shift-to-product} requires injectivity of the shifted
operators on a common parameter interval for all $\tau\ge0$. We compare
the mass scale with the distance to the singular point. Write
\[
L_p:=\Delta+p(u_p^{\rm th})^{p-1}.
\]

\begin{proposition}\label{prop:near-threshold-shift}
There exists \(\eta_{\rm th,2}>0\) such that
\[
\ker(L_p-\tau)=\{0\}
\]
for every
\[
p\in[p_*,p_*+\eta_{\rm th,2}),
\qquad
\tau\ge0,
\]
where the kernel is taken in
$\mathscr X^0_{\mu,\nu}(\mathcal O_0)$ with weights
\eqref{eq:near-threshold-weights}. The same conclusion holds
simultaneously for any fixed finite collection of singular weights in
the first root-free gap above \(-d\), with the same far weight \(\nu\).
\end{proposition}

\begin{proof}
At \(p=p_*\), the case \(\tau=0\) is
Lemma~\ref{lem:threshold-semistable-H}.  For \(\tau>0\), testing
\((L_*-\tau)\Psi=0\) against \(\eta_{\rho,R}^2\Psi\) gives
\[
Q_*(\eta_{\rho,R}\Psi)
+\tau\int\eta_{\rho,R}^2\Psi^2
=
\int\Psi^2|\nabla\eta_{\rho,R}|^2
=
O(\rho^{n+2\mu-2})+O(R^{n+2\nu-2}),
\]
which tends to zero because \(\mu>-d\) and \(\nu<-d\).  Hence
\(\Psi=0\).

\proofstep{1}{Construct a normalized counterexample sequence.}
Assume uniform injectivity fails.  Then
\[
p_j\downarrow p_*,
\qquad
\tau_j\ge0,
\qquad
(L_{p_j}-\tau_j)\Psi_j=0,
\qquad
\sup_X\frac{|\Psi_j(X)|}{\omega_{\mu,\nu}(|X|)}=1.
\]
Choose $X_j$ with normalized value at least $1/2$ and put
$r_j=|X_j|$, $c_j=\tau_jr_j^2$, and
$\varrho_j=\operatorname{dist}(X_j,\partial\R^n_+)$. At scale $s_j=r_j$
when $c_j$ is bounded, or $s_j=\ell_j=\tau_j^{-1/2}$ when
$c_j\to\infty$, the scaled coefficients are bounded. On a corresponding
boundary patch $U_j$ away from $0$, weights are comparable and
\[
 \frac12\le\frac{|\Psi_j(X_j)|}{\omega_{\mu,\nu}(r_j)}
 \le C\frac{\varrho_j}{s_j}
       \sup_{U_j}\frac{|\Psi_j|}{\omega_{\mu,\nu}(r_j)}.
\]
Thus the boundary estimate prevents loss of the normalization at either scale.

\proofstep{2}{Exclude the case $\tau_jr_j^2\to\infty$.}
First suppose $c_j\to\infty$ and set
\[
\ell_j=\tau_j^{-1/2}\ll r_j,\qquad
W_j(Z)=\frac{\Psi_j(X_j+\ell_jZ)}{\omega_{\mu,\nu}(r_j)}.
\]
The uniform two-ended bounds on the threshold branch imply
\[
 |X|^2p_j(u_{p_j}^{\rm th})^{p_j-1}\le C.
\]
On every fixed rescaled ball,
$|X_j+\ell_jZ|\asymp r_j$, so the scaled potential is $O(c_j^{-1})$
and the weight ratio tends to $1$. The boundary estimate above preserves
a nonzero limit. Thus $W_j$ converges locally to a bounded nonzero
solution of $\Delta W-W=0$ in $\R^n$ or a Dirichlet half-space.
Odd reflection and the Fourier transform exclude both possibilities.

\proofstep{3}{Exclude the three limits when $\tau_jr_j^2$ is bounded.}
We may therefore assume $c_j\to c<\infty$. There remain three spatial
scales.

If $r_j\asymp1$, then $\tau_j$ is bounded, and a subsequence
produces a nonzero threshold kernel of $L_*-\tau_\infty$, already excluded.

If \(r_j\to0\), the uniform inner
asymptotics imply
\[
r_j^2p_j(u_{p_j}^{\rm th}(r_jZ))^{p_j-1}
\longrightarrow
p_*|Z|^{-2}\phi_*(Z/|Z|)^{p_*-1}.
\]
The functions \(r_j^{-\mu}\Psi_j(r_jZ)\) converge, after a diagonal extraction, to a
nonzero solution satisfying \(|\Psi(Z)|\le C|Z|^\mu\) of
\[
\left(
\Delta+p_*|Z|^{-2}\phi_*^{p_*-1}-c
\right)\Psi=0.
\]
Lemma~\ref{lem:conic-modes}, applied with $\sigma=\mu$, excludes this
limit: $\mu$ is nonindicial and $\mu>-d=\operatorname{Re}\gamma_1^-(p_*)$.

Finally let \(r_j\to\infty\).  Since
\[
u_{p_j}^{\rm th}(X)\le C|X|^{1-n},
\qquad
(n-1)(p_*-1)>2,
\]
the rescaled potential tends to zero and
\(r_j^{-\nu}\Psi_j(r_jZ)\) converges to a nonzero solution of
\[
(\Delta-c)\Psi=0
\quad\text{in }\R^n_+\setminus\{0\},
\qquad
|\Psi(Z)|\le C|Z|^\nu.
\]
Lemma~\ref{lem:conic-modes} with $V_0=0$ and $\sigma=\nu$ excludes this
limit, since $\nu\in(1-n,1)$ is nonindicial for the harmonic model and
$\nu>1-n=\operatorname{Re}\gamma_1^-$.

Each alternative is excluded whenever the singular weight lies in the
first fixed root-free gap above \(-d\) and \(1-n<\nu<-d\). Hence the
same argument applies to any fixed finite collection of such singular weights.
\end{proof}

\subsubsection{Non-degeneracy in the product domain}
For the threshold branch, the frequency argument establishes H\textup{(ii)} on a
parameter interval independent of the tangential dimension and the
allowed polynomial growth.

\begin{proposition}\label{prop:near-threshold-H}
There exists
\[
0<\eta_*<
\frac{n+2}{n-2}-p_*
\]
such that, for every
\[
p\in[p_*,p_*+\eta_*),
\]
the normalized threshold branch cell \(u_p^{\rm th}\) satisfies
Assumptions~H\textup{(i)} and H\textup{(ii)}, for every \(k\ge1\) in H\textup{(ii)}.
\end{proposition}

\begin{proof}
Choose $-d<\mu_0<\mu$ in the same first root-free gap and apply
Proposition~\ref{prop:near-threshold-shift} simultaneously at these
two weights. Fix its resulting interval $\eta_*(n)$ before choosing
$k$, the polynomial growth order $m$, a compact frequency set $K$,
or a derivative order $M>m+k+2$.
The inner construction, Remark~\ref{rem5}, and the bounded dilation
factors imply, uniformly on this interval,
\[
 p(u_p^{\rm th})^{p-1}
 =p r^{-2}\phi_p^{p-1}+O(r^{-2+\sigma'})
\]
for some $\sigma'>0$; for $p<2$ the H\"older continuity of the
power may decrease $\sigma'$. At infinity the potential is
$O(r^{-(n-1)(p-1)})$, with $(n-1)(p-1)>2$.
Lemma~\ref{lem:shift-to-product} therefore establishes H\textup{(i)} and H\textup{(ii)}:
its inner Green exponent is $\mu-\mu_0>0$, and all frequency derivative
orders are allowed after the interval has been fixed.
Local elliptic bootstrapping raises $q$ to any finite $q>N$
without changing the interval.

\end{proof}

\subsection{Generic non-degeneracy}\label{sec:generic-nondegeneracy}
Outside the intervals just treated, we use the analytic branches of
Section~\ref{sec:analytic-continuation}. The operator for H\textup{(i)} has a
stronger inner weight than the matching linearization, so it requires
its own Fredholm argument. We keep the far-field coefficient fixed at
$c$ while applying analyticity, and normalize the selected cells only
afterward.

\begin{proposition}
\label{prop:generic-H}
For the analytic continuation branch in Theorem~\ref{thm4}, there is a
locally finite set $\Sigma_H\subset(0,\infty)$ such that, for every
$s\notin\Sigma_H$, there exist exponents $\mu_s,\nu_s$ with
\[
\operatorname{Re}\gamma_1^+(p_s)<\mu_s,
\qquad
1-n<\nu_s<-\frac{n-2}{2},
\]
for which the branch representative $u_s$ satisfies Assumption H\textup{(i)} of
Section~\ref{sec:assumption-H}.  Moreover, $\mu_s$ may be chosen so close to
$\operatorname{Re}\gamma_1^+(p_s)$ that there is a nonindicial low weight
$\delta_s$ satisfying
\[
-\frac2{p_s-1}<\delta_s<
\min\left\{0,1-\frac2{p_s-1},\operatorname{Re}\gamma_1^-(p_s),
2-n-\mu_s\right\}.
\]
Thus $\delta_{*,s}:=2-n-\delta_s>\mu_s$, so the weights can be chosen to
satisfy \eqref{e59} and~\eqref{e59-twoend} with $n=N$ when Theorem~\ref{thm2}
is applied.
\end{proposition}
\begin{proof}
\proofstep{1}{Establish the Fredholm property in fixed weighted spaces.}
Fix a local analytic branch chart. Choose $\mu$ in the first open
indicial gap, sufficiently close to
$\operatorname{Re}\gamma_1^+(p_s)$ that $\mu<2-n+a_{p_s}$, and choose
$1-n<\nu<-(n-2)/2$. Remark~\ref{rem:weight-compatibility} then supplies
the nonindicial low weight $\delta_s$ in the statement. Shrinking the
chart keeps $(\mu,\nu)$ admissible without crossing any indicial root.

On the fixed spaces set
\[
\mathcal L_s:W^{2,q}_{\mu,\nu;0}(\R^n_+)
\longrightarrow L^q_{\mu-2,\nu-2}(\R^n_+),
\qquad
\mathcal L_s\Phi=-\Delta\Phi-p_su_s^{p_s-1}\Phi.
\]
At the singular end,
multiplication by $r^2$ and the logarithmic change of variables produce
the constant-coefficient cylinder operator with angular part
$-\Delta_\theta-p_s\phi_{p_s}^{p_s-1}$. Conjugation by $r^\mu$
has invertible mode symbols because $\mu$ avoids the real parts of
the indicial roots. At infinity the limiting operator is the
Dirichlet Laplacian, and $\nu\in(1-n,1)$ avoids its roots.
The resulting end inverses, joined with local elliptic parametrices
on a compact annulus, produce a Fredholm operator. The scaled potential
errors tend to zero at both ends; their small tails and local Rellich
compactness make them compact from the weighted domain to the source
space. The same end estimates, applied successively to a kernel
element, improve its weight across any root-free subinterval.
Consequently the kernels agree for the different choices of weights
within the stated gaps.

\proofstep{2}{Prove analytic dependence along the branch.}
Analyticity follows directly on the coefficient-$c$ branch.
Boundary division and the inner construction imply, locally in a
branch chart,
\[
\begin{gathered}
 0<m\le w_s/\theta_n\le M,\qquad 0<m\le b_s/\theta_n\le M,\\
 b_s(\theta,t)=e^{(n-1)t}u_s(e^t\theta),\qquad t\ge0,
\end{gathered}
\]
with analytic dependence of these quotients in $L^\infty$.
For the exterior quotient, positivity follows from the admissible set
and the far-field asymptotic; the $W$ norm controls boundary division.
At the two ends,
\[
\begin{gathered}
 r^2u_s^{p_s-1}=
 \begin{cases}
  w_s^{p_s-1},&r\le1,\\
  r^{-\beta(s)}b_s^{p_s-1},&r\ge1,
 \end{cases}\\
 \beta(s)=(n-1)(p_s-1)-2>0.
\end{gathered}
\]
Lemma~\ref{lem:positive-boundary-power} treats both boundary powers.
After fixing a positive lower bound $\beta_0$ for the decay margin,
\[
 \sup_{r\ge1}r^{-\beta_0}(\log r)^j\le j!\beta_0^{-j}
\]
controls every radial parameter derivative. Thus multiplication by
$p_su_s^{p_s-1}$ is analytic between the displayed fixed spaces, and
$\mathcal L_{s^*(\sigma)}$ is an analytic Fredholm family.

\proofstep{3}{Prove local finiteness of the exceptional parameters.}
By the uniqueness in Lemma~\ref{lem13}\textup{(i)}, the initial arc consists,
up to dilation, of the stable cells from Section~\ref{sec:stable-construction}.
Lemma~\ref{lem:stable-ground-state} therefore provides an injective point.
The fixed dilation to coefficient $c$ preserves the two-ended kernel
class, so this conclusion applies to $\mathcal L_s$.

On each connected analytic chart, finite-dimensional Fredholm
reduction shows that the noninjective set is either the whole chart
or locally discrete. Indeed, the reduced analytic matrix
$S(\sigma):\R^m\to\R^\ell$ has the same kernel; wherever it is injective, some $m\times m$
minor is nonzero, so noninjectivity is confined to isolated zeros. The
set where noninjectivity is locally discrete is therefore open and
closed, hence it is the whole connected chart once it is nonempty. This argument does
not require index zero.

On an open overlap, choose a point outside the first chart's discrete
exceptional set. The root-free weight improvement above identifies
the two Jacobi kernels there, giving an injective point in the next
chart. A finite chain of overlaps along each compact parameter
interval therefore propagates initial injectivity and proves that
$\Sigma_H$ is locally finite.

Finally, the scale-invariant interior and boundary estimates on dyadic
half-annuli place every Jacobi field satisfying the pointwise two-ended
bounds in H\textup{(i)} in the displayed weighted domain.  Thus triviality of the
weighted kernel is H\textup{(i)}.
\end{proof}

To obtain the exceptional set of exponents from the exceptional branch
parameters, we use the following elementary observation.

\begin{lemma}\label{lem:parameter-projection}
Let $p:(0,\infty)\to(p_-,p_+)$ be continuous, with
$p(s)\to p_-$ as $s\downarrow0$ and $p(s)\to p_+$ as $s\to\infty$.
If $\Sigma\subset(0,\infty)$ is locally finite, then for every compact
$K\Subset(p_-,p_+)$ the set
\[
\{p(s):s\in\Sigma,\ p(s)\in K\}
\]
is finite. Moreover $p((0,\infty))=(p_-,p_+)$.
\end{lemma}
\begin{proof}
The endpoint limits place $p^{-1}(K)$ in a compact subinterval of
$(0,\infty)$, which meets the locally finite set $\Sigma$ in only
finitely many points. This proves the first assertion. By the same limits and the intermediate value theorem, $p$ attains
every value in $(p_-,p_+)$; no monotonicity is needed.
\end{proof}

\begin{proposition}\label{prop:generic-exponents}
For $n\ge3$, let
\[
\mathcal E_n\subset\left(\frac{n+1}{n-1},\frac{n+2}{n-2}\right)
\]
be the set of exponents for which there is no coefficient-$1$ fast-decaying
cell satisfying H\textup{(i)} for any pair of weights satisfying \eqref{e59-twoend}.
Then $\mathcal E_n$ is locally finite. For every $p\notin\mathcal E_n$
in this interval there is a cell satisfying H\textup{(i)} with weights compatible
with the linear and nonlinear gluing estimates.
\end{proposition}
\begin{proof}
Fix a compact interval $K$ in the strict subcritical range. Choose one
sufficiently small truncation $\epstr$ so that $K\Subset(p_-,p_+)$,
and take its branch in Theorem~\ref{thm4}. By
Proposition~\ref{prop:generic-H}, every parameter outside a
locally finite set $\Sigma_H$ supplies H\textup{(i)}. Normalize each selected
cell to coefficient $1$; the fixed-parameter kernel isomorphism in
Remark~\ref{rem:normalization-continuity} preserves H\textup{(i)}, with
\[
-a_p<\delta<\min\{0,1-a_p,\operatorname{Re}\gamma_1^-,2-n-\mu\}.
\]
Thus the weights satisfy \eqref{e59-twoend} and
Remark~\ref{rem:weight-compatibility}, including the nonlinear restriction.
Every exponent in $K$ has a preimage by Lemma~\ref{lem:parameter-projection}.
If an exponent belongs to $\mathcal E_n\cap K$, each of its preimages
on this branch must belong to $\Sigma_H$. Consequently
\[
\mathcal E_n\cap K\subset
\{p_s:s\in\Sigma_H,\ p_s\in K\},
\]
and the right-hand side is finite by the same lemma. This proves the
claim on every compact interval using a single truncation containing it.
\end{proof}

The argument establishes generic H\textup{(i)}; it supplies no generic H\textup{(ii)}.

\section{Geometry and approximate solutions}\label{s6}

We use the half-space cell to construct an approximate solution near
the prescribed singular set, then solve for a correction on $\Omega$.
Let $\Omega\subset\R^N$ be a smooth bounded domain, let $0\le k<N-2$,
and set $n=N-k$. For $k=0$, let $S\subset\partial\Omega$ be finite;
for $k>0$, let $S$ be a smooth compact embedded $k$-dimensional
submanifold of $\partial\Omega$. Fix a sufficiently small tubular radius
$R>0$ for the Fermi coordinates below, independently of $\ve$, and then
choose $\ve_0>0$ sufficiently small; throughout, $0<\ve<\ve_0$.
Fix a fast-decaying half-space cell $u_1$ satisfying H\textup{(i)} if $k=0$,
or H\textup{(ii)} if $k>0$, as defined in Section~\ref{sec:model-linear}.

\subsection{Fermi coordinates near the singular set}\label{sec:fermi}
Choose Fermi coordinates near $S$ as in \cite[Section~\textup{3.1}]{dmp}.
For $y\in S$ in a coordinate neighborhood, let $E_{k+1},\ldots,E_{N-1}$
be a local orthonormal frame of the normal bundle of $S$ in
$\partial\Omega$, and let $E_N$ be the inward unit normal to
$\partial\Omega$. For small $\bar x=(x_{k+1},\ldots,x_{N-1})$ and
$x_N\ge0$, set
\[
\begin{aligned}
F(y,\bar x)&=\operatorname{Exp}_{y}
  \Big(\sum_{j=k+1}^{N-1}x_jE_j(y)\Big),\\
\bar F(y,X)&=F(y,\bar x)+x_NE_N(F(y,\bar x)),
\qquad X=(\bar x,x_N).
\end{aligned}
\]
The Laplacian then takes the form
\begin{equation}\label{e52}
\Delta_{\R^N}=\Delta_{NS}+O(|X|)\nabla^2+O(1)\nabla,
\end{equation}
where
$X:=(x_{k+1},\ldots,x_N)\in\R^n_+$ and
$\Delta_{NS}=\Delta_X+\Delta_S$ is the model Laplace--Beltrami operator on the normal bundle of $S$.

\subsection{The approximate solution}\label{sec:approximate-solution}
For the fixed radius $R$, set
\begin{equation*}
B_+(R):=\{(y,X)\in NS: 0<|X|<R,\ x_N>0\}.
\end{equation*}

For sufficiently small $\ve>0$, set
\begin{equation}\label{e54}
 u_\ve(y,X):=\ve^{-\frac{2}{p-1}}u_1\left(\frac{X}{\ve}\right),
 \qquad
 \bar u_\ve:=\eta_R u_\ve.
\end{equation}
Two orthonormal frames of the normal bundle in $\partial\Omega$ differ
by an element of $O(n-1)$, while the inward normal is fixed.
Equation~\eqref{eq:cell-rotation-invariance} therefore makes \eqref{e54}
independent of the frame.

Choose the cutoff to depend only on $|X|$,
with $0\le\eta_R\le1$, $\eta_R=1$ in $\overline B_+(R/2)$ and
$\eta_R=0$ in $\overline\Omega\setminus B_+(R)$.  Then $\bar u_\ve=0$ on $\partial\Omega\setminus S$.

We seek a solution $u=\bar u_\ve+v$ by solving the correction equation
\begin{equation}\label{e56}
\left\{\begin{array}{l}
L_\ve v+\mathcal E_\ve+\mathcal N(v)=0 \text{ in }\Omega\\
v=0 \text{ on }\partial\Omega\setminus S
\end{array}
\right.,
\end{equation}
where
\begin{equation*}
L_\ve:=\Delta+p\bar u_\ve^{p-1},\qquad \mathcal E_\ve:=\Delta\bar u_\ve+\bar u_\ve^p
\end{equation*}
and
\begin{equation*}
\mathcal N(v):=(\bar u_\ve+v)_+^p-\bar u_\ve^p-p\bar u_\ve^{p-1}v
\end{equation*}
is the nonlinear remainder. The positive part allows the equation to be
solved before positivity of $u$ is known.

To solve \eqref{e56}, we first construct a right inverse for $L_\ve$
with a bound independent of $\ve$. The relevant model depends on the
rescaling radius $r$:
\[
\begin{array}{c|c}
r/\ve\to0 & \text{homogeneous model}\\
r/\ve\to\lambda\in(0,\infty) & \text{cell model, after a fixed dilation}\\
r/\ve\to\infty & \text{harmonic model}.
\end{array}
\]
The cell linearization is $L_1$ of Section~\ref{sec:model-linear}.
All three limits enter the linear estimate below.

\section{Uniform linear estimates}
\label{s8}

We now solve the linear equation associated with \eqref{e56}. The
correction is sought at a low weight, where the kernel may be nontrivial.
We therefore select the solution orthogonal to the Hilbert kernel and
prove a uniform bound for this choice. The complementary high-weight estimate implies surjectivity. We obtain
the uniform bound from the equations for the canonical solution pair.

Throughout this section, fix $q>N$.  Choose the low weight $\delta$ by
\begin{equation}
\label{e59}
-\frac2{p-1}<\delta<
\min\left\{0,\,1-\frac2{p-1},\,\operatorname{Re}\gamma_1^-,\,2-n-\mu\right\},
\qquad
\delta_*:=2-n-\delta,
\end{equation}
with $\delta$ and $\delta_*$ away from the indicial spectrum.
By Remark~\ref{rem:weight-compatibility}, $\mu<\delta_*$ and
\[
\delta<\operatorname{Re}\gamma_1^-
\le-\frac{n-2}{2}
\le\operatorname{Re}\gamma_1^+<\delta_*,
\qquad
1-n<\delta_*<1.
\]
Set
\[
\beta:=(n-1)(p-1)>2,
\]
which will be used in the finite-cell bootstrap and in the nonlinear estimates.
Here $(\mu,\nu)$ is a pair for which the chosen cell satisfies H\textup{(i)}
or H\textup{(ii)}, as appropriate. These parameters
remain fixed throughout the linear and nonlinear gluing arguments.

Consider the Dirichlet problem
\begin{equation}\label{e58}
\left\{\begin{array}{l}
L_\ve(v)=f\text{ in }\Omega\\
v=0 \text{ on }\partial\Omega\setminus S
\end{array}
\right..
\end{equation}

Near $S$, write $(y,X)\in S\times\R^n_+$
in the Fermi coordinates of Section~\ref{s6}. Fix a smooth positive function
\[
\rho:\overline\Omega\setminus S\to (0,\infty)
\]
which equals $|X|$ near $S$. Choose its extension so that
$\rho(z)\asymp\rho(x)$ whenever $z\in B_{\rho(x)/2}(x)\cap\Omega$.
We use $\rho$ to define the global weighted norms; the estimates for $u_1$ are those of Theorem~\ref{thm1}.

\subsection{Weighted spaces on the punctured domain}\label{sec:global-weighted-spaces}

For $\gamma\in\R$, let
\[
 L^q_\gamma(\Omega\setminus S)
 :=\{f\in L^q_{\rm loc}(\overline\Omega\setminus S):
                         \|f\|_{L^q_\gamma}<\infty\},
\]
where
\[
\|f\|_{L^q_\gamma(\Omega\setminus S)}
:=
\sup_{x\in\Omega\setminus S}
\rho(x)^{-\gamma-\frac Nq}
\|f\|_{L^q(B_{\rho(x)/2}(x)\cap\Omega)} .
\]
For $\delta\in\R$, the corresponding Sobolev space is
\[
 W^{2,q}_\delta(\Omega\setminus S)
 :=\{\phi\in W^{2,q}_{\rm loc}(\overline\Omega\setminus S):
                        \|\phi\|_{W^{2,q}_\delta}<\infty\},
\]
with norm
\[
\|\phi\|_{W^{2,q}_\delta(\Omega\setminus S)}
:=
\sum_{j=0}^2
\|\nabla^j\phi\|_{L^q_{\delta-j}(\Omega\setminus S)} .
\]

Let $W^{2,q}_{\delta,0}(\Omega\setminus S)$ denote
the closed subspace of $W^{2,q}_\delta(\Omega\setminus S)$ consisting of
functions whose trace vanishes on $\partial\Omega\setminus S$.

We use the following scaled H\"older spaces for $j\ge0$ and
$0<\alpha<1$. Set $E_x=B_{\rho(x)/2}(x)\cap\Omega$ and
\begin{equation}\label{eq:scaled-holder-spaces}
\begin{aligned}
 \|w\|_{C^{j,\alpha}_b}
 :=\sum_{\ell=0}^j\sup_x\rho(x)^{\ell-b}|\nabla^\ell w(x)|
    +\sup_x\rho(x)^{j+\alpha-b}[\nabla^jw]_{\alpha;E_x},\\
 C^{j,\alpha}_b(\Omega\setminus S)
 :=\{w\in C^{j,\alpha}_{\rm loc}(\overline\Omega\setminus S):
                                  \|w\|_{C^{j,\alpha}_b}<\infty\}.
\end{aligned}
\end{equation}
Here $[\cdot]_{\alpha;E}$ is the usual H\"older seminorm on $E$.
The same notation on a homogeneous model uses $\rho=r=|X|$.
Since $q>N$, scaled Sobolev embedding gives
\[
\begin{gathered}
W^{2,q}_\delta(\Omega\setminus S)
\hookrightarrow C^{1,\alpha}_\delta(\Omega\setminus S),
\qquad \alpha=1-\frac Nq,\\
|\phi(x)|\le C\rho(x)^\delta\|\phi\|_{W^{2,q}_\delta}.
\end{gathered}
\]

\subsection{Non-degeneracy of the homogeneous model}\label{sec:homogeneous-liouville}
At scales smaller than $\ve$, the rescaled potential approaches that
of $u_\infty$. We will use the following Liouville property to exclude
this limit: there is no nontrivial solution of
\[
\Delta_Y\Phi+\Delta_X\Phi+pu_\infty(X)^{p-1}\Phi=0
\quad\text{in }\R^k\times(\R^n_+\setminus\{0\})
\]
with zero boundary value, polynomial growth in $Y$, and
$|\Phi(Y,X)|\le C(1+|Y|)^m|X|^{\delta_*}$; that is, its kernel in
$\mathscr X^m_{\delta_*,\delta_*}(\mathcal O_k)$ is zero for every
finite $m$. For $k=0$, this is Lemma~\ref{lem:conic-modes}\,\textup{(i)}, since
$\delta_*$ is nonindicial. For $k>0$, we verify the test-inverse
hypothesis of Lemma~\ref{lem:fourier-testing}. At a nonzero frequency,
a homogeneous solution of the normal equation which is $O(r^\sigma)$ at
zero, with $\sigma>\operatorname{Re}\gamma_1^+$, and polynomially
bounded at infinity vanishes by Lemma~\ref{lem:conic-modes}\,\textup{(ii)}.
Hence the high massive kernel is zero.

Choose $\vartheta>0$ so small that
$s=\delta_*-\vartheta>\operatorname{Re}\gamma_1^+$ and
$[s,\delta_*]$ is root-free. Lemma~\ref{lem:scalar-massive} and
Green duality make the complementary low realization surjective at
$\check s=2-n-s=\delta+\vartheta$. Its Hilbert graph kernel is complemented,
so it has a bounded right inverse $R_{\xi_0}$ at each fixed frequency.
On a neighborhood of $\xi_0$ a smooth right inverse is
\[
 R_\xi=R_{\xi_0}
 \bigl[I+(B_\xi-B_{\xi_0})R_{\xi_0}\bigr]^{-1}.
\]
A frequency partition of unity patches these inverses on each compact
set away from zero. Dyadic estimates and the massive exterior estimate
provide the test bounds~\eqref{eq:fourier-test-estimates} with inner
exponent $\check s$. The boundary exponent is
\[
 \delta_*+\check s+n-2=\vartheta>0.
\]
Thus Lemma~\ref{lem:fourier-testing} implies
$\operatorname{supp}_\xi\widehat\Phi\subset\{0\}$, so $\Phi$ is a
polynomial in $Y$. Each coefficient inherits the normal bound
$r^{\delta_*}$ by finite-dimensional polynomial interpolation. Its
highest-degree coefficients solve the zero-frequency conic equation;
the $k=0$ argument makes them zero. Descending in degree proves
$\Phi=0$.

\subsection{Injectivity and a uniform estimate}
We first apply local elliptic estimates on balls whose radii are
comparable to the distance from $S$. The resulting estimate is valid
at any weight and will be used in both the high- and low-weight arguments.
\begin{lemma}
\label{lem4}
Let $1<s<\infty$, and let $\sigma\in\R$.  There exists a constant
$C>0$, independent of $0<\ve<\ve_0$, such that for every
$\phi\in W^{2,s}_{\sigma,0}(\Omega\setminus S)$,
\[
\|\phi\|_{W^{2,s}_\sigma(\Omega\setminus S)}
\le
C
\left(
\|L_\ve\phi\|_{L^s_{\sigma-2}(\Omega\setminus S)}
+
\|\phi\|_{L^s_\sigma(\Omega\setminus S)}
\right).
\]
\end{lemma}

\begin{proof}
On a bounded-overlap Whitney covering $B_i=B_{r_i}(x_i)\cap\Omega$,
$r_i=\rho(x_i)/8$, rescale $x=x_i+r_i z$ and
$\phi=r_i^\sigma\widetilde\phi_i$. The rescaled potential satisfies
\[
|V_i(z)|=|pr_i^2\bar u_\ve(x_i+r_i z)^{p-1}|
\le C r_i^2\ve^{-2}
\begin{cases}
(r_i/\ve)^{-2},&r_i\le\ve,\\
(r_i/\ve)^{-(n-1)(p-1)},&r_i\ge\ve,
\end{cases}
\le C,
\]
since $(n-1)(p-1)>2$. In a fixed small tubular neighborhood, the
rescaled $O(r_i)$ second-order error is absorbed in the interior and
boundary Calder\'on--Zygmund estimates. Scaling back and taking weighted
suprema proves the claim there, uniformly in $\ve$; a fixed finite
covering handles the complement.
\end{proof}

\begin{proposition}\label{prop4}
Assume H\textup{(i)} when $k=0$ and H\textup{(ii)} when $k>0$. For every $1<s<\infty$ there exist $\ve_0>0$ and $C>0$ such that
\begin{equation}\label{e60}
\|\phi\|_{W^{2,s}_{\delta_*}(\Omega\setminus S)}
\le C
\|L_\ve\phi\|_{L^s_{\delta_*-2}(\Omega\setminus S)}
\end{equation}
for all $0<\ve<\ve_0$ and
$\phi\in W^{2,s}_{\delta_*,0}(\Omega\setminus S)$.  In particular,
\[
\ker\bigl(L_\ve:W^{2,s}_{\delta_*,0}\to L^s_{\delta_*-2}\bigr)=\{0\}.
\]
The same conclusion holds, with constants depending on $\sigma$, when
$\delta_*$ is replaced by any nonindicial weight $\sigma\in[\delta_*,1)$.
\end{proposition}

\begin{proof}
We write the proof for $\delta_*$. It uses only that $\delta_*$ is
nonindicial and satisfies
$\max\{\mu,\operatorname{Re}\gamma_1^+\}<\delta_*<1$; every
$\sigma$ as in the last assertion has the same properties.

\proofstep{1}{Obtain a nonzero limit after rescaling.}
Assume that \eqref{e60} fails.  There are $\ve_j\to0$ and
$L_{\ve_j}\phi_j=f_j$ with
\[
\|\phi_j\|_{W^{2,s}_{\delta_*}}=1,
\qquad
\|f_j\|_{L^s_{\delta_*-2}}\to0.
\]
Lemma~\ref{lem4} provides points $x_j$ and $r_j=\rho(x_j)$ such that
\begin{equation}\label{eq:high-weight-local-mass}
r_j^{-\delta_*-N/s}
\|\phi_j\|_{L^s(B_{r_j/2}(x_j)\cap\Omega)}\ge c_0>0.
\end{equation}
If $r_j\ge r_0>0$, then $\bar u_{\ve_j}\to0$ away from $S$ and local
compactness produces a harmonic Dirichlet limit.  The high-weight bound implies
$O(\rho^{\delta_*})$ near $S$; since $\delta_*>1-n$, the boundary singular
set is removable.  The limit is zero, contradicting
\eqref{eq:high-weight-local-mass}.  Hence $r_j\to0$.

In Fermi coordinates, set
\[
\widetilde\phi_j(Y,Z)=r_j^{-\delta_*}
\phi_j(y_j+r_jY,r_jZ),
\qquad \lambda_j=\frac{r_j}{\ve_j}.
\]
Local elliptic estimates and Rellich compactness ensure, after passing
to a subsequence, weak $W^{2,s}_{\rm loc}$ and strong $L^s_{\rm loc}$
convergence to $\Phi\ne0$, by \eqref{eq:high-weight-local-mass}.
For $r=|Z_0|>0$, the local estimate for the limiting equation gives
\[
|\Phi(Y_0,Z_0)|
\le Cr^{-N/s}\|\Phi\|_{L^s(B_{r/4}(Y_0,Z_0))}
\le Cr^{\delta_*}.
\]
This pointwise bound holds for every $1<s<\infty$ by elliptic
bootstrapping, without requiring $s>N$. Since the rescaled potential is
\[
pr_j^2\bar u_{\ve_j}^{p-1}
=p\eta_R^{p-1}\lambda_j^2u_1(\lambda_jZ)^{p-1},
\]
the ratio $\lambda_j=r_j/\ve_j$ determines the limiting equation.

\proofstep{2}{Exclude the three limiting equations.}
If $\lambda_j\to0$, then
\[
(\Delta_Y+\Delta_Z)\Phi
+p|Z|^{-2}\phi_p^{p-1}\Phi=0.
\]
The weight $\delta_*>\Re\gamma_1^+$ lies above the real part of both
first-mode indicial roots and avoids the indicial spectrum; the Liouville property of Section~\ref{sec:homogeneous-liouville} implies
$\Phi=0$.

If $\lambda_j\to\lambda_0\in(0,\infty)$, a fixed dilation gives
\[
(\Delta_Y+\Delta_Z)\Psi+pu_1^{p-1}\Psi=0.
\]
The inner bound is stronger than $O(|Z|^\mu)$.  At infinity,
$u_1^{p-1}=O(|Z|^{-\beta})$ with $\beta>2$.  Choose
\[
\delta_*=\kappa_0>\kappa_1>\cdots>\kappa_M=\nu
\quad\text{in }(1-n,1),
\qquad
\kappa_{\ell+1}>\kappa_\ell-(\beta-2).
\]
If $\Psi=O(|Z|^{\kappa_\ell})$, then
\[
u_1^{p-1}\Psi=O(|Z|^{\kappa_\ell-\beta})
=o(|Z|^{\kappa_{\ell+1}-2}).
\]
Since the harmonic Dirichlet indicial spectrum has no root in $(1-n,1)$,
the exterior estimate implies
$\Psi=O(|Z|^{\kappa_{\ell+1}})$. Iteration therefore establishes
$\Psi=O(|Z|^\nu)$ at infinity. Rescaling preserves a uniform bound in the
tangential variable, hence the limit has the polynomial tangential growth
allowed in H\textup{(ii)}.  Assumption H\textup{(i)}, or H\textup{(ii)} for $k>0$, implies $\Psi=0$.

If $\lambda_j\to\infty$, then
\[
\lambda_j^2u_1(\lambda_jZ)^{p-1}\to0
\quad\text{locally on }\R^n_+\setminus\{0\},
\]
so $\Phi$ is harmonic with zero flat-boundary data and
$|\Phi|\le C|Z|^{\delta_*}$, $1-n<\delta_*<1$.  Removability of the edge
and odd reflection produce an entire harmonic function of sublinear growth,
hence zero.  Each case contradicts the normalization, proving \eqref{e60}.
\end{proof}

\subsection{A right inverse with a uniform bound}

We next solve \eqref{e58} at the low weight $\delta$. The kernel may be
nontrivial and, when $k>0$, infinite dimensional. We select the solution
orthogonal to it in a weighted Hilbert space, as in
\cite[Section~\textup{8}]{acdfgw}. The high-weight estimate just proved will
exclude the adjoint kernel, while Appendix~\ref{app:fixed-regularity}
places the selected solution in $W^{2,q}_{\delta,0}$ for each fixed
$\ve$. We then prove that its bound is uniform as $\ve\to0$.

For $\sigma\in\R$, the weighted Hilbert space is
\[
 \mathcal H_\sigma:=\left\{h:\|h\|_{\mathcal H_\sigma}^2
 :=\int_{\Omega\setminus S}\rho^{-2\sigma-n}|h|^2\,dx<\infty\right\}.
\]
Define its zero-trace edge Sobolev space by
\begin{equation}\label{eq:global-edge-Sobolev-space}
\begin{aligned}
 \mathscr H^2_{\sigma,0}
 :=\biggl\{v\in H^2_{\rm loc}(\overline\Omega\setminus S):\ &
 v|_{\partial\Omega\setminus S}=0,\\
 &\|v\|_{\mathscr H^2_\sigma}^2
 :=\sum_{j=0}^2\|\nabla^jv\|_{\mathcal H_{\sigma-j}}^2<\infty\biggr\}.
\end{aligned}
\end{equation}
The powers $\rho^j$ in this norm measure derivatives on balls of radius
comparable to $\rho$. The closed Dirichlet realization is
\[
 L_{\ve,\sigma}:\mathcal D_\sigma\subset\mathcal H_\sigma
                  \longrightarrow\mathcal H_{\sigma-2},\qquad
 \mathcal D_\sigma:=\mathscr H^2_{\sigma,0},\quad
 L_{\ve,\sigma}v=L_\ve v.
\]
Equivalently, $\mathcal D_\sigma$ is the space of $v\in\mathcal H_\sigma$
with $L_\ve v\in\mathcal H_{\sigma-2}$ distributionally and zero
physical Dirichlet trace. Local scaled elliptic estimates show that
its edge norm is equivalent to the graph norm
\[
 \|v\|_{\mathcal D_\sigma}
 :=\|v\|_{\mathcal H_\sigma}
       +\|L_\ve v\|_{\mathcal H_{\sigma-2}}.
\]

\begin{lemma}\label{lem:edge-weight-transfer}
Fix $\ve>0$ and choose $a=\delta_-<\delta<-(n-2)/2$ sufficiently
close to $\delta$ that $-a_p<2a-\delta$ and
$I=[2a-\delta,\delta]$ contains no real part of an indicial root. Set
\[
\delta_+=2-n-a,\qquad \widehat\delta=\delta-2a-n+2.
\]
The high realization $L_\ve:\mathcal D_{\delta_+}\to
\mathcal H_{\delta_+-2}$ has closed range and finite-dimensional kernel,
and every kernel element belongs to $W^{2,2}_{\delta_*,0}$.
If $A_\ve:\mathcal D_a\subset\mathcal H_a\to\mathcal H_{a-2}$
is surjective, its canonical pair for $h\in L^q_{\delta-2}$ is
\[
u=A_\ve^*(A_\ve A_\ve^*)^{-1}h,\qquad
z=\rho^{-2a+4-n}(A_\ve A_\ve^*)^{-1}h.
\]
It belongs to $W^{2,q}_{\delta,0}\times W^{2,q}_{\widehat\delta,0}$ and satisfies
\begin{equation}\label{eq:canonical-fixed-regularity}
\|u\|_{W^{2,q}_\delta}+\|z\|_{W^{2,q}_{\widehat\delta}}
\le C_\ve\|h\|_{L^q_{\delta-2}}.
\end{equation}
The constant may depend on the fixed parameter $\ve$.
\end{lemma}

\begin{proof}
Apply Theorem~\ref{thm:fixed} in Appendix~\ref{app:fixed-regularity}
with $V_\ve=p\bar u_\ve^{p-1}$ and $V_0=p\phi_p^{p-1}$.
Theorem~\ref{thm1} ensures uniform convergence of the scaled potential
on closed normal half-annuli, while the Fermi-coordinate errors in
\eqref{e52} tend to zero. The chosen root-free intervals are those of
the appendix. Its fixed-parameter estimates place the high kernel and the canonical
Hilbert pair in the required weighted Sobolev spaces, proving
\eqref{eq:canonical-fixed-regularity}.
\end{proof}

\begin{proposition}\label{prop5}
Under the hypotheses of Proposition~\ref{prop4}, there exist $\ve_0>0$, $C_G>0$, and operators
\[
G_\ve:
L^q_{\delta-2}(\Omega\setminus S)
\longrightarrow
W^{2,q}_{\delta,0}(\Omega\setminus S),
\qquad 0<\ve<\ve_0,
\]
such that
\[
L_\ve G_\ve=I
\]
and
\begin{equation}\label{eq:uniform-right-inverse-low-weight}
\|G_\ve f\|_{W^{2,q}_\delta}
\le C_G\|f\|_{L^q_{\delta-2}}.
\end{equation}
\end{proposition}

\begin{proof}
\proofstep{1}{Construct the canonical Hilbert right inverse.}
Choose $\delta_-<\delta$ sufficiently close to $\delta$ so that
\[
-a_p<\check\delta:=2\delta_- -\delta<\delta_-<\delta,
\qquad
\delta_*<\delta_+:=2-n-\delta_-<\widehat\delta<1,
\]
where
\[
\widehat\delta=\delta-2\delta_- -n+2,
\qquad
\delta_-+\delta_+=\check\delta+\widehat\delta=2-n,
\]
and the intervals $[\check\delta,\delta]$ and
$[\delta_*,\widehat\delta]$ contain no indicial root.

Regard
\[
A_\ve:
\mathcal H_{\delta_-}\supset\mathcal D_{\delta_-}
\longrightarrow\mathcal H_{\delta_- -2},
\qquad A_\ve u=L_\ve u,
\]
as a densely defined closed operator. Put
$\mathcal R_\sigma w=\rho^{-2\sigma-n}w$.
Since $\delta_-+\delta_+=2-n$, the Riesz maps identify the paired
Hilbert spaces isometrically. The adjoint calculation in
Appendix~\ref{app:fixed-regularity} gives
\[
A_\ve^*=\mathcal R_{\delta_-}^{-1}L_{\ve,\delta_+}
\mathcal R_{\delta_- -2},\qquad
D(A_\ve^*)=\mathcal R_{\delta_- -2}^{-1}\mathcal D_{\delta_+}.
\]

By Lemma~\ref{lem:edge-weight-transfer}, the high realization has closed
range and its kernel lies in $W^{2,2}_{\delta_*,0}$; Proposition~\ref{prop4}
makes this kernel zero. Thus $A_\ve^*$ is injective with closed range,
and the closed-range theorem implies, for each fixed sufficiently small $\ve$,
\[
\operatorname{Ran}A_\ve=(\ker A_\ve^*)^\perp=\mathcal H_{\delta_- -2}.
\]
Consequently
$A_\ve A_\ve^*$ is invertible on $\mathcal H_{\delta_- -2}$ and
\begin{equation}\label{eq:canonical-Hilbert-inverse}
G_\ve^H=A_\ve^*(A_\ve A_\ve^*)^{-1}
\end{equation}
is the minimal-norm Hilbert right inverse.  In particular,
$A_\ve G_\ve^Hh=h$.

\proofstep{2}{Establish weighted regularity of the solution pair.}
Each shell $\{r/2<\rho<2r\}$ has $O(r^{-k})$ Whitney balls.
By H\"older's inequality,
\[\|h\|_{L^2}^2\le Cr^{2(\delta-2)+n}\|h\|_{L^q_{\delta-2}}^2.\]
The Hilbert weight leaves $r^{2(\delta-\delta_-)}$, summable over
dyadic shells since $\delta>\delta_-$. Hence
\begin{equation}\label{eq:Lq-Hilbert}
\|h\|_{\mathcal H_{\delta_- -2}}
\le C\|h\|_{L^q_{\delta-2}}.
\end{equation}
For this datum, the canonical Hilbert pair is
\[
u=G_\ve^Hh,\qquad
z=\mathcal R_{\delta_- -2}(A_\ve A_\ve^*)^{-1}h.
\]
Then
\begin{equation}\label{eq:paired-Hilbert-system}
L_\ve u=h,\qquad L_\ve z=\mathcal R_{\delta_-}u,\qquad
u=z=0\quad\text{on }\partial\Omega\setminus S.
\end{equation}
Lemma~\ref{lem:edge-weight-transfer} implies
\[
(u,z)\in W^{2,q}_{\delta,0}\times W^{2,q}_{\widehat\delta,0},
\]
with the bound \eqref{eq:canonical-fixed-regularity} for fixed $\ve$. We use the
second equation to control $z$ uniformly in terms of $u$. Put
$a=\delta_-$. The weight identities are
\[
 \widehat\delta-\delta_*=2(\delta-a)>0,\qquad
 \widehat\delta-2=\delta-2a-n.
\]
Thus multiplication by $\mathcal R_a=\rho^{-2a-n}$ gives
\[
 \|\mathcal R_a u\|_{L^q_{\widehat\delta-2}}
 \le C\|u\|_{L^q_\delta}.
\]
Since $\widehat\delta\in[\delta_*,1)$ is nonindicial, the last
assertion of Proposition~\ref{prop4} applies with
$\sigma=\widehat\delta$. Applying it to the second equation yields
\begin{equation}\label{eq:adjoint-lift}
\|z\|_{W^{2,q}_{\widehat\delta}}
\le C\|u\|_{W^{2,q}_\delta}.
\end{equation}

\proofstep{3}{Rescale the solution pair at a concentration scale.}
To bound $u$ uniformly, we rescale the coupled system. Its second
equation is needed to exclude the kernel allowed by the low weight. Assume that the $W^{2,q}_\delta$ estimate fails. After choosing a sequence
$\ve_j\downarrow0$ and normalizing,
\[
u_j=G_{\ve_j}^Hh_j,
\qquad
\|u_j\|_{W^{2,q}_\delta}=1,
\qquad
\|h_j\|_{L^q_{\delta-2}}\to0.
\]
The rescaling argument in
Proposition~\ref{prop4}, now using the low weight $\delta$ and
Lemma~\ref{lem4}, provides points $x_j$ with $r_j=\rho(x_j)$ such that
\[
 r_j^{-\delta-N/q}
 \|u_j\|_{L^q(B_{r_j/2}(x_j)\cap\Omega)}\ge c_0>0.
\]
Here $r_j\to0$: otherwise the harmonic Dirichlet limit has the bound
$O(\rho^\delta)$, and $1-n<\delta<1$ makes its boundary singularity
removable, contradicting this normalization. Rescale both components of
\eqref{eq:paired-Hilbert-system} in Fermi coordinates:
\[
\widetilde u_j(Y,Z)=r_j^{-\delta}u_j(y_j+r_jY,r_jZ),
\qquad
\widetilde z_j(Y,Z)=r_j^{-\widehat\delta}z_j(y_j+r_jY,r_jZ).
\]
A subsequence converges to $(\Phi,Z_0)$, $\Phi\not\equiv0$, satisfying one
of the three model systems
\begin{equation}\label{eq:paired-model-system}
\mathcal L_{\rm mod}\Phi=0,
\qquad
\mathcal L_{\rm mod}Z_0=|Z|^{-2\delta_- -n}\Phi.
\end{equation}

\proofstep{4}{Exclude a nonzero limit when $k=0$.}
Choose $1-n<\check\delta<\sigma_\infty<\delta_-<\delta$ with
$[\sigma_\infty,\delta]$ root-free. The conic expansion or the finite-cell
exterior bootstrap of Proposition~\ref{prop4} implies
$\Phi=O(r^{\sigma_\infty})$ at infinity. Put
$a=\delta_-$ and $b=\sigma_\infty-2a-n+2$. Then
$\mathcal L_{\rm mod}Z_0=f=r^{-2a-n}\Phi=O(r^{b-2})$, while
$Z_0=O(r^{\widehat\delta})$ and
$\delta_*<b<\delta_+<\widehat\delta<1$.
For the conic and harmonic models, variation of constants improves
this to $Z_0=O(r^b)$: the inherited bound excludes powers above
$\widehat\delta$, the intervening interval is root-free, and the source
is nonresonant.

For the finite fast-cell model write
$\Delta Z_0=f-VZ_0$, where $V=O(r^{-\beta})$ and $\beta>2$.
Choose a finite sequence
\[
 \widehat\delta=b_0>b_1>\cdots>b_M=b,\qquad
 b_{j+1}>b_j-(\beta-2).
\]
If $Z_0=O(r^{b_j})$, then
\[
 f-VZ_0=O(r^{b-2})+O(r^{b_j-\beta})
       =O(r^{b_{j+1}-2}).
\]
The same harmonic exterior estimate, with
$[b_{j+1},b_j]\subset(1-n,1)$, gives $Z_0=O(r^{b_{j+1}})$.
Iteration reaches $b$; local scaled estimates ensure the same bounds
for first derivatives. The inner Green term is
$O(r^{\delta+\widehat\delta+n-2})=O(r^{2(\delta-a)})=o(1)$;
the outer term is $O(r^{2(\sigma_\infty-a)})=o(1)$.

Green's identity on half-annuli then gives
\[
\int_{\R^n_+}|Z|^{-2\delta_- -n}|\Phi|^2
=\int (\mathcal L_{\rm mod}Z_0)\Phi
=\int Z_0\mathcal L_{\rm mod}\Phi=0,
\]
hence $\Phi=0$, a contradiction.

\proofstep{5}{Exclude a nonzero limit when $k>0$.}
For $k>0$, fix $K\Subset\R^k\setminus\{0\}$ and write
\[
 a=\delta_-,\qquad \varpi(r)=r^{-2a-n},\qquad
 B_\xi=\Delta_X-|\xi|^2+pU_{\rm mod}^{p-1},\quad r=|X|.
\]
Use the scalar realizations of Lemma~\ref{lem:scalar-massive},
with the spaces \eqref{eq:massive-two-end-spaces} and Green duality
\eqref{eq:massive-green-dual}. Its hypotheses hold for all three model
potentials. In particular the graph domains are fixed on compact
frequency charts, and the scalar kernels are independent of the
exterior polynomial weights.

Choose $0<\vartheta<\delta-a$ sufficiently small and set
\[
 s=\delta-\vartheta,\qquad \widehat s=\widehat\delta-\vartheta,\qquad
 t=\sigma_\infty<a,\qquad \widehat t=t-2a-n+2.
\]
Here $t$ is fixed and finite; the inner weights remain in the previously
chosen root-free intervals. For smooth $(\psi_1,\psi_2)$ compactly
supported in $K$ and in a normal annulus, we construct the transpose tests
\begin{equation}\label{eq:massive-transpose-test}
 B_\xi\eta-\varpi\zeta=\psi_1,\qquad B_\xi\zeta=\psi_2,
\end{equation}
in the Green-dual domains
\[
 \eta\in\mathcal D_{2-n-s,-n-t},\qquad
 \zeta\in\mathcal D_{2-n-\widehat s,-n-\widehat t}.
\]
The coupling maps the second domain into the first source space, since
\[
 (2-n-\widehat s)-2a-n=(2-n-s)-2,\qquad
 (-n-\widehat t)-2a-n=(-n-t)-2.
\]
Thus it matches the inner source weight exactly and gains a factor
$r^{-2}$ at infinity.

The high scalar realizations are injective. In the cell case their
kernel elements decay exponentially at infinity and have the inner
high weight, hence $e^{iY\cdot\xi}v(X)$ belongs to the class in H\textup{(ii)}.
The homogeneous conic and harmonic cases are excluded by their
high-weight model Liouville results. Green duality therefore makes
the low equation for $\zeta$ surjective. Its kernel $K_\xi$ represents
the range annihilator of the high equation for $\eta$: the low inner weights
\[
 2-n-\widehat s=\check\delta+\vartheta<s=\delta-\vartheta
\]
lie in the same root-free interval, and exponential decay makes the
exterior polynomial weights immaterial for kernels.
On this common kernel, the Gram form
\[
 \mathfrak g_\xi(\kappa,\psi)
 =\int\varpi\kappa\overline\psi\,dX
\]
is finite and positive definite. Indeed, its elements belong to the
inner space of weight $s>a$, so weighted Cauchy--Schwarz ensures
integrability at zero; at infinity they decay exponentially.
Lemma~\ref{lem:gram-correction}, with variables and outputs interchanged,
makes the transpose system an isomorphism. Its inverse is smooth on
compact frequency charts because the graph domains are fixed.

Scaled elliptic estimates and the fixed graph inverse imply, for every
fixed multi-index $\alpha$, uniformly for $\xi\in K$,
\begin{equation}\label{eq:massive-test-ends}
\begin{aligned}
 |\partial_\xi^\alpha\eta|+r|\nabla_X\partial_\xi^\alpha\eta|
 &\le C_\alpha r^{2-n-\delta+\vartheta},\\
 |\partial_\xi^\alpha\zeta|+r|\nabla_X\partial_\xi^\alpha\zeta|
 &\le C_\alpha r^{2a-\delta+\vartheta},\qquad r\le1.
\end{aligned}
\end{equation}
Local elliptic bootstrapping establishes these bounds; differentiation adds
only previously controlled sources in the same spaces. At infinity
the compact data and positive mass ensure exponential decay of the
tests and all finite frequency derivatives. The two inner Green
exponents for $(\Phi,Z_0)$ are
\[
 \delta+(2-n-\delta+\vartheta)+n-2=\vartheta,
 \qquad \widehat\delta+(2a-\delta+\vartheta)+n-2=\vartheta.
\]
The coupling is bounded by $C\int_0^1r^{\vartheta-1}\,dr$.
Lemma~\ref{lem:fourier-testing} establishes zero-frequency support for both
components.

Both components are therefore polynomials in $Y$ as normal distributions. Their
uniform boundedness in $Y$ makes these polynomials constant. The
$k=0$ paired Green argument above then implies $\Phi=0$, contradicting
the normalization. The test functions decay against the inherited
exterior bounds; since the frequency support is separated from zero,
no inverse estimate as $\xi\to0$ is needed.

Thus
\[
\|G_\ve^Hh\|_{W^{2,q}_\delta}
\le C\|h\|_{L^q_{\delta-2}}
\]
uniformly for small $\ve$.  Taking
$G_\ve=G_\ve^H|_{L^q_{\delta-2}}$ proves the proposition.
\end{proof}

\section{Completion of the gluing proofs}\label{s9}

The uniform right inverse reduces \eqref{e56} to a contraction on a
small ball in $W^{2,q}_{\delta,0}$. We estimate the residual and the
nonlinear remainder, retaining the boundary factor to control the
correction near $S$. The resulting solution is positive and satisfies
the very weak boundary identity.

\subsection{Estimate of the error}
Keep $q$, $\delta$, and $\delta_*$ fixed as in Section~\ref{s8}, put
$a_p=2/(p-1)$, and set
\[
\begin{aligned}
\varsigma&:=1-a_p-\delta>0,\qquad
\varsigma_1:=\min\{\varsigma,\,n-1-a_p\}>0.
\end{aligned}
\]
After fixing the geometric radius and the linear bound $C_G$, choose
$R_0>1$ so large that the term $C R_0^{2-\beta}$ in the outer
nonlinear estimate below is less than $1/(4C_G)$.
Denote the resulting residual constant by $C_E=C_E(R_0)$.
We then choose $M>2C_GC_E$ and finally decrease $\ve_0$, in particular
so that $R_0\ve_0<R/2$.

By Proposition~\ref{prop5}, solving \eqref{e56} for
$v\in W^{2,q}_{\delta,0}(\Omega\setminus S)$ reduces to
\begin{equation}
\mathcal K(v):=-G_\ve\bigl(\mathcal E_\ve+\mathcal N(v)\bigr)=v.
\end{equation}
The equation for $u_1$ and the geometric expansion \eqref{e52} give
\begin{align*}
\mathcal E_\ve
&=\Delta \bar u_\ve+\bar u_\ve^p\\
&=(\eta_R^p-\eta_R)u_\ve^p
 +(\Delta\eta_R)u_\ve
 +2\nabla\eta_R\cdot\nabla u_\ve\\
&\quad
 +O(|X|)\nabla^2(\eta_Ru_\ve)
 +O(1)\nabla(\eta_Ru_\ve).
\end{align*}

\begin{lemma}\label{lem:residual-estimate}
For $\ve>0$ sufficiently small,
\begin{equation}\label{eq:error-final}
\|\mathcal E_\ve\|_{L^q_{\delta-2}(\Omega\setminus S)}
\le C\ve^{\varsigma_1}.
\end{equation}
\end{lemma}

\begin{proof}
On $\rho\le R_0\ve<R/2$, the cutoff is one. The bounds
$|\nabla^ju_\ve|\le C_{R_0}\rho^{-a_p-j}$, $j=1,2$, imply
\[
 |\mathcal E_\ve|\le C_{R_0}\rho^{-a_p-1},\qquad
 \|\mathcal E_\ve\|_{L^q_{\delta-2}(\{\rho\le R_0\ve\})}
 \le C_{R_0}\ve^{1-a_p-\delta}.
\]
On $R_0\ve\le\rho\le R/2$, use
$|\nabla^ju_\ve|\le C\ve^{n-1-a_p}\rho^{1-n-j}$ to obtain
\[
 |\mathcal E_\ve|\le C\ve^{n-1-a_p}\rho^{-n},\qquad
 \|\mathcal E_\ve\|_{L^q_{\delta-2}(\{R_0\ve\le\rho\le R/2\})}
 \le C\bigl(\ve^{1-a_p-\delta}+\ve^{n-1-a_p}\bigr),
\]
by taking the maximum of $r^{2-n-\delta}$ at the two endpoints.
On the fixed cutoff annulus, $|\nabla^ju_\ve|\le
C_R\ve^{n-1-a_p}$ for $0\le j\le2$; all linear errors are of this
order, as is the nonlinear error
$O_R(\ve^{p(n-1-a_p)})$, which is of smaller order. Thus its weighted norm is
$O_R(\ve^{n-1-a_p})$. The residual vanishes outside the support.
Comparability of $\rho$ on the norm's balls now gives
\[
 \|\mathcal E_\ve\|_{L^q_{\delta-2}}
 \le C\bigl(\ve^{1-a_p-\delta}+\ve^{n-1-a_p}\bigr)
 \le C\ve^{\varsigma_1}.
\]
\end{proof}

\subsection{The nonlinear estimate and contraction}
The residual estimate implies
$\|G_\ve\mathcal E_\ve\|_{W^{2,q}_\delta}\le C_GC_E\ve^{\varsigma_1}$.
For the chosen $M>2C_GC_E$, set
\[
\mathcal B_{M,\ve}:=
\left\{v\in W^{2,q}_{\delta,0}(\Omega\setminus S):
\|v\|_{W^{2,q}_\delta}\le M\ve^{\varsigma_1}\right\}.
\]
Thus $\|G_\ve\mathcal E_\ve\|_{W^{2,q}_\delta}\le(M/2)\ve^{\varsigma_1}$.
To estimate $\mathcal N$ on this ball, we first use the zero trace to
recover the same angular boundary factor as in $\bar u_\ve$.

\begin{lemma}\label{lem:boundary-factor}
If $w\in W^{2,q}_{\delta,0}(\Omega\setminus S)$ and $\rho\le R$, then
\begin{equation}\label{eq:weighted-boundary-factor}
|w(y,X)|\le C\rho^\delta\theta_n
\|w\|_{W^{2,q}_\delta},
\qquad \theta_n:=\frac{X_n}{|X|}.
\end{equation}
\end{lemma}

\begin{proof}
If $\theta_n\ge1/2$, the estimate follows from
$|w|\le C\rho^\delta\|w\|_{W^{2,q}_\delta}$.  If
$\theta_n<1/2$, then $|X'|\simeq\rho$ and, along the vertical segment
$(y,X',t)$, $0\le t\le X_n$, the distance to $S$ is comparable to $\rho$.
Since the trace of $w$ vanishes at $(y,X',0)$, the weighted gradient estimate implies
\[
|w(y,X',X_n)|
\le X_n\sup_{0<t<X_n}|\partial_{X_n}w(y,X',t)|
\le C X_n\rho^{\delta-1}\|w\|_{W^{2,q}_\delta},
\]
which is \eqref{eq:weighted-boundary-factor}.
\end{proof}

On $0<\rho\le R_0\ve$, the cutoff is one and the cell estimates imply
$\bar u_\ve\asymp\rho^{-a_p}\theta_n$, with constants depending on $R_0$.
Thus every $v\in\mathcal B_{M,\ve}$ satisfies, in the interior,
\begin{equation}\label{eq:inner-correction-small}
 \frac{|v|}{\bar u_\ve}
 \le C M\ve^{\varsigma_1}\rho^{\delta+a_p}
 \le C_{R_0}M\ve^{\varsigma_1+\delta+a_p}=o(1).
\end{equation}
Here \eqref{eq:weighted-boundary-factor} cancels the boundary factor,
and $\varsigma_1>0$, $\delta+a_p>0$. In particular,
$|v|\le\bar u_\ve/4$ for sufficiently small $\ve$, uniformly in the ball.
We may therefore use Taylor's formula in the inner region; in the
outer region we use a difference estimate for the positive-part power.

\begin{lemma}\label{lem6}
There exists $\ve_1\in(0,\ve_0]$ such that, for every $0<\ve<\ve_1$ and $v_1,v_2\in\mathcal B_{M,\ve}$,
\[
\|\mathcal N(v_1)-\mathcal N(v_2)\|_{L^q_{\delta-2}}
\le \frac{1}{2C_G}\|v_1-v_2\|_{W^{2,q}_\delta}.
\]
\end{lemma}
\begin{proof}
Write $w=v_1-v_2$ and $D=\|w\|_{W^{2,q}_\delta}$.

\smallskip
\noindent\emph{The inner region.}
On $\rho\le R_0\ve$, \eqref{eq:inner-correction-small} permits Taylor's
formula along the positive segment between $\bar u_\ve+v_1$ and
$\bar u_\ve+v_2$. The boundary-factor estimate gives
\[
\begin{aligned}
 |\mathcal N(v_1)-\mathcal N(v_2)|
 &\le C\bar u_\ve^{p-2}(|v_1|+|v_2|)|w|\\
 &\le C_{R_0}M\ve^{\varsigma_1}
       \rho^{\delta-2}\rho^{\delta+a_p}\theta_n^p D.
\end{aligned}
\]
The boundary factor cancels even when $1<p<2$:
\[
\theta_n^{p-2}\theta_n^2=\theta_n^p\le1.
\]
Consequently,
\[
 \|\mathcal N(v_1)-\mathcal N(v_2)\|_{L^q_{\delta-2}(\{\rho\le R_0\ve\})}
 \le C_{R_0}M\ve^{\varsigma_1+\delta+a_p}D.
\]

\smallskip
\noindent\emph{The outer region.}
On the entire outer region $\{\rho\ge R_0\ve\}$, use
$\bar u_\ve\le C\ve^{n-1-a_p}\rho^{1-n}$ and the positive-part inequality
\[
 |\mathcal N(v_1)-\mathcal N(v_2)|
 \le C\bigl(\bar u_\ve^{p-1}+|v_1|^{p-1}+|v_2|^{p-1}\bigr)|w|.
\]
This includes the complement of the cutoff support, where
$\mathcal N(v)=(v_+)^p$. With the globally defined $\rho$,
\[
\begin{aligned}
 |\mathcal N(v_1)-\mathcal N(v_2)|
 &\le C\Bigl[\ve^{\beta-2}\rho^{2-\beta}
       +(M\ve^{\varsigma_1})^{p-1}\rho^{2+(p-1)\delta}\Bigr]
       \rho^{\delta-2}D\\
 &\le C\Bigl[R_0^{2-\beta}+(M\ve^{\varsigma_1})^{p-1}\Bigr]
       \rho^{\delta-2}D.
\end{aligned}
\]
Here $\beta>2$, $2+(p-1)\delta>0$, and $\rho$ is bounded above.
Taking weighted local norms yields the same coefficient times $D$.

The fixed choice of $R_0$ controls its term; now decrease $\ve$,
with $M$ fixed, so that each regional bound
is at most $D/(4C_G)$. Adding the two bounds proves the claim, including
balls meeting both regions.
\end{proof}

\smallskip
\noindent\emph{Existence of the fixed point.}
Since $\mathcal N(0)=0$, \eqref{eq:error-final} and Lemma~\ref{lem6} give
\[
\begin{aligned}
 \|\mathcal K(v)\|_{W^{2,q}_\delta}
 &\le C_GC_E\ve^{\varsigma_1}+\tfrac12\|v\|_{W^{2,q}_\delta}
 \le M\ve^{\varsigma_1},\\
 \|\mathcal K(v_1)-\mathcal K(v_2)\|_{W^{2,q}_\delta}
 &\le\tfrac12\|v_1-v_2\|_{W^{2,q}_\delta}.
\end{aligned}
\]
Thus $\mathcal K$ has a unique fixed point in $\mathcal B_{M,\ve}$.

\subsection{Positivity and the very weak formulation}\label{sec:gluing-proof-completion}

The fixed point defines $u=\bar u_\ve+v$ satisfying
\begin{equation}\label{e68}
\begin{cases}
\Delta u+u_+^p=0&\text{in }\Omega,\\
u=0&\text{on }\partial\Omega\setminus S.
\end{cases}
\end{equation}
By \eqref{eq:inner-correction-small}, $u>0$ on
$\{\rho\le R_0\ve\}\cap\Omega$ for small $\ve$.
Any connected component $D$ of $\{u<0\}$ is therefore separated from
$S$. On $D$, $\Delta u=0$ and $u=0$ on $\partial D$, so the maximum
principle excludes $D$. Thus $u\ge0$, and the strong maximum principle implies
\[
u>0\quad\hbox{in }\Omega\setminus S.
\]
Equation~\eqref{e68} reduces to
\eqref{e1}.

For the very weak formulation, the estimate near $S$ is
\[
u(x)\sim \rho(x)^{-\frac2{p-1}}\phi_p(\theta).
\]
Since $a_p<n-1$ and $p a_p=a_p+2$,
\[
 \int_0^1\bigl(r^{n-1-a_p}+r^{n-a_p-2}\bigr)\,dr<\infty.
\]
The boundary distance is $O(r)$ in these coordinates, so these integrals imply
\[
 u\in L^1(\Omega),\qquad
 \operatorname{dist}(x,\partial\Omega)u^p\in L^1(\Omega).
\]

To obtain the identity, remove a tube around $S$: let
$\Omega_r:=\Omega\cap\{\rho>r\}$ and
$\zeta\in C^2(\overline\Omega)$ with $\zeta=0$ on $\partial\Omega$.  Then
\[
\int_{\Omega_r}(u\Delta\zeta+u^p\zeta)\,dx
=
\int_{\{\rho=r\}}
\bigl(u\partial_{\nu_r}\zeta-\zeta\partial_{\nu_r}u\bigr)\,d\sigma_r.
\]
On $\{\rho=r\}$,
\[
|\zeta|\le Cr,
\qquad
|u|\le Cr^{-\frac2{p-1}},
\qquad
|\nabla u|\le Cr^{-\frac2{p-1}-1},
\qquad
d\sigma_r=O(r^{n-1}),
\]
and therefore
\[
\left|
\int_{\{\rho=r\}}
\bigl(u\partial_{\nu_r}\zeta-\zeta\partial_{\nu_r}u\bigr)\,d\sigma_r
\right|
\le Cr^{n-1-\frac2{p-1}}=o(1).
\]
Letting $r\downarrow0$ gives the very weak identity.  Finally,
$v=o(\bar u_\ve)$ non-tangentially along $S$, hence $u(x)\to\infty$ non-tangentially as $x\to S$.  Theorems~\ref{thm2} and~\ref{thm3} follow.

\subsection{Proof of Corollary~\ref{cor:gluing-applications}}\label{sec:completion-main}
The model results of Section~\ref{sec:model-linear} verify the hypotheses
of the two gluing theorems in the stated ranges.

\begin{proof}
For \textup{(i)}, in the stable range, Lemmas~\ref{lem:stable-ground-state}
and~\ref{lem:stable-edge-nondeg} establish H\textup{(i)} and H\textup{(ii)}, respectively.
At an interior stability threshold and on the right neighborhood
specified in Proposition~\ref{prop:near-threshold-H}, that proposition
verifies both hypotheses. Theorems~\ref{thm2} and~\ref{thm3} therefore
apply for $k=0$ and $k>0$, respectively, throughout the stated ranges.

For \textup{(ii)}, apply Proposition~\ref{prop:generic-exponents} with $n=N$
and then Theorem~\ref{thm2} to the resulting cell.
\end{proof}

\appendix
\section{Auxiliary estimates for analytic continuation}\label{app:continuation-estimates}
The analytic continuation in Section~\ref{sec:analytic-continuation}
uses the spectral bounds and the inner and exterior estimates proved below.
We retain its notation, spaces, and boundary data, with normal dimension $n$
and exponent cutoff $\epstr$. Constants are uniform on a fixed truncated
exponent interval, independently of the gluing scale $\ve$.

\subsection{Eigenvalue estimates}\label{app:proof-prop7}
\begin{proof}[Proof of Proposition~\ref{prop7}]
The profile equation and \eqref{e7}, with $\phi_p,\psi_p^1>0$, give
\begin{equation}\label{e96}
 \lambda_p^1\int_{S^{n-1}_+}\phi_p\psi_p^1\,d\sigma
 =-(p-1)\int_{S^{n-1}_+}\phi_p^p\psi_p^1\,d\sigma<0.
\end{equation}
    Since the first eigenspace is one-dimensional and the operator is invariant under rotations of the first $n-1$ coordinates, $\psi_p^1$ is axisymmetric; hence $\lambda_p^1=\lambda_p^{0,1}$.

\smallskip
\noindent\emph{The axisymmetric sector.}
    For the axisymmetric sector $m=0$, extend
    $\phi_p$ evenly across the equator and denote the extension by
    $\overline\phi_p$.  On the full sphere, restrict the operator
    \[
      \mathcal A_p^{(0)}
      :=-\Delta_{S^{n-1}}+\beta_p-p\overline\phi_p^{p-1}
    \]
    to the axisymmetric subspace.  This is a one-dimensional
    Sturm--Liouville operator; let its eigenvalues be
    \[
      \eta_p^1<\eta_p^2<\eta_p^3<\cdots.
    \]
    They are simple. The $\ell$th Dirichlet eigenfunction in the
    $m=0$ sector has $\ell-1$ zeros inside the northern half-meridian.
    Its odd extension has the reflected zeros and one zero at the equator,
    so it has $2(\ell-1)+1=2\ell-1$ interior zeros on the full meridian.
    Sturm oscillation therefore places it at index $j=2\ell$. In particular,
    for $\ell\ge2$,
    \begin{equation}\label{e97}
      \lambda_p^{0,\ell}=\eta_p^{2\ell},\qquad 2\ell\ge4.
    \end{equation}

    On the other hand, the homogeneous solution
    $u_p^*(x)=|x|^{-2/(p-1)}\phi_p(x/|x|)$ satisfies
    $-\Delta(\partial_{x_n}u_p^*)
    =p(u_p^*)^{p-1}\partial_{x_n}u_p^*$.  Since
    $\partial_{x_n}u_p^*$ is homogeneous of degree
    $-(p+1)/(p-1)$, its restriction to the upper hemisphere has
    axisymmetric radial factor
    \[
      G_p(s)
      =
      \cos s\,\Phi_p'(s)-a_p\sin s\,\Phi_p(s).
    \]
    The profile equation at \(s=0\), together with
    \(\Phi_p(0)=0\), implies \(\Phi_p''(0)=0\).  Hence
    \[
      G_p'(0)=0.
    \]
    The even extension of $G_p$ is therefore $C^1$ across the equator
    and solves the Sturm--Liouville equation weakly. Elliptic regularity
    places it in the full-sphere domain, where it is an axisymmetric
    eigenfunction of $\mathcal A_p^{(0)}$ with eigenvalue
    \begin{equation}\label{e98}
      \alpha_p+1.
    \end{equation}
    As $p\downarrow (n+1)/(n-1)$,
    $\overline\phi_p\to0$, $\beta_p\to-(n-1)$, and the axisymmetric
    eigenvalues converge to
    \[
      \ell(\ell+n-2)-(n-1),\qquad \ell=0,1,2,\ldots.
    \]
    The value $\alpha_p+1$ tends to $n+1$, which is the third
    eigenvalue in this list.  Therefore, for $p$ close to the lower
    endpoint, $\alpha_p+1=\eta_p^3$. Simplicity and continuous dependence
    separate each eigenvalue from its neighbors locally in $p$, so the
    index $j(p)$ in $\alpha_p+1=\eta_p^{j(p)}$ is locally constant.
    Connectedness of the exponent interval implies $j(p)\equiv3$. Hence
    \begin{equation}\label{e99}
      \alpha_p+1=\eta_p^3
      \quad\hbox{for all }
      p\in\left(\frac{n+1}{n-1},\frac{n+2}{n-2}\right).
    \end{equation}
    Combining \eqref{e97} and \eqref{e99} yields
    \begin{equation}\label{e100}
      \lambda_p^{0,\ell}>\alpha_p+1,
      \qquad \ell\ge2.
    \end{equation}

\smallskip
\noindent\emph{The nonaxisymmetric sectors.}
    For the degree-one sector, the moving-plane reduction on
the hemisphere \cite[Section~\textup{3}]{bpv1} shows that the positive profile is
strictly decreasing in the geodesic distance from the north pole.
Since \(s\) increases from the equator to the north pole,
\[
\Phi_p'(s)>0,
\qquad
0<s<\frac{\pi}{2}.
\]
For $a=1,\ldots,n-1$, the tangential translation mode
$-\partial_{x_a}u_p^*$ has the form
    \[
      -\partial_{x_a}u_p^*|_{S^{n-1}_+}
      =Y_{1,a}(\omega)\Psi_{p}^{\rm tr}(s).
    \]
    It solves \eqref{e94a} with $m=1$ and
    $\lambda=\alpha_p+1$. Direct differentiation, with
    $\Phi_p(s)=\phi_p(\cos s\,\omega,\sin s)$, gives, up to a positive
    normalization constant,
    \[
      \Psi_p^{\rm tr}(s)
      =a_p\cos s\,\Phi_p(s)+\sin s\,\Phi_p'(s)>0.
    \]
    Hence this is the principal eigenfunction in the degree-one sector, and
    \begin{equation}\label{e101}
      \lambda_p^{1,1}=\alpha_p+1.
    \end{equation}
    All higher eigenvalues in the same sector are strictly larger.

Finally, the degree-$m$ principal eigenfunction for $m\ge2$ belongs
to the degree-one form domain. The angular potential increases by
$(\nu_m-\nu_1)/\cos^2s\ge\nu_m-\nu_1>0$, where $\nu_1=n-2$.
The Rayleigh characterization therefore implies
\begin{equation}\label{e102}
 \lambda_p^{m,1}>\lambda_p^{1,1}=\alpha_p+1,\qquad m\ge2.
\end{equation}
The higher eigenvalues in these sectors satisfy the same bound;
together with \eqref{e96} and \eqref{e100}, this proves \eqref{e95}.
\end{proof}
\subsection{Analyticity of the spherical profile}\label{app:proof-lem:profile-analytic}
The profile vanishes on the equator. We use its positive boundary quotient
to justify analytic dependence of the noninteger power.

\begin{lemma}[Analytic powers]
\label{lem:positive-boundary-power}
Let $(\Theta,\mathfrak m)$ be a measure space and let $d:\Theta\to(0,1]$ be
measurable. Let $E$ be a real Banach space of measurable functions such
that $u\mapsto u/d$ is bounded from $E$ to $L^\infty(\Theta)$. Define
\[
 \mathcal U_d:=\{u\in E:\operatorname*{ess\,inf}_{\Theta}u/d>0\}.
\]
Then $(\zeta,u)\mapsto u^\zeta$ is real analytic from
$(0,\infty)\times\mathcal U_d$ to $L^\infty(\Theta)$. The local bounds are
uniform when $\zeta$ ranges over a compact subinterval of $(0,\infty)$,
the norm of division by $d$ is fixed, and $0<m\le u/d\le M$.
\end{lemma}
\begin{proof}
Write $\mathfrak q=u/d$ and $u^\zeta=d^\zeta\exp(\zeta\log\mathfrak q)$.
The logarithm is analytic in the Banach algebra $L^\infty$ near any
$\mathfrak q_0\ge m>0$, with radius controlled by $m$. Moreover,
\[
 \sup_{0<d\le1}d^{\zeta_0}|\log d|^j\le j!\zeta_0^{-j},
 \qquad \zeta_0>0.
\]
Thus both factors are analytic with the asserted uniform bounds.
\end{proof}

\begin{proof}[Proof of Lemma~\ref{lem:profile-analytic}]
By Proposition~\ref{prop7}, $\lambda_p^1<0<\lambda_p^2$, so
\[
\mathcal A_p:W^{2,q}(S^{n-1}_+)\cap W^{1,q}_0(S^{n-1}_+)
\longrightarrow L^q(S^{n-1}_+)
\]
is an isomorphism.

By the Hopf lemma, $\phi_p/\theta_n\ge c>0$ on the closed
hemisphere. Boundary division maps
$W^{2,q}\cap W^{1,q}_0$ continuously into $L^\infty$ for $q>n$.
Lemma~\ref{lem:positive-boundary-power}, with $d=\theta_n$, therefore
makes
\[
 (p',u)\longmapsto-\Delta_{S^{n-1}}u+\beta_{p'}u-u^{p'}
\]
real analytic into $L^q$ near $(p,\phi_p)$.
Its $u$-derivative is the isomorphism $\mathcal A_p$ above. The analytic
implicit function theorem and uniqueness of the positive profile construct
the local analytic branch $p'\mapsto\phi_{p'}$.

The resulting family $\mathcal A_p$ is analytic on the fixed Dirichlet
domain. Fix $\bar p$ in the exponent interval. Its lowest eigenvalue is simple. Choose a positively oriented contour
$\Gamma$ enclosing only $\lambda_{\bar p}^1$. For $p$ near $\bar p$, analytic
perturbation theory \cite[Chapter~\textup{VII}, Section~\textup{1.3}]{kato} gives
\[
 \mathsf P_1(p)=\frac{1}{2\pi i}\int_\Gamma
                      (z-\mathcal A_p)^{-1}\,dz,\qquad
 \psi_p^1=\frac{\mathsf P_1(p)\psi_{\bar p}^1}
                   {\|\mathsf P_1(p)\psi_{\bar p}^1\|_2}.
\]
The resolvent is analytic from $L^q$ to the fixed $W^{2,q}$ Dirichlet
domain, so these formulas give the stated analytic eigenfunction in
that space. Positivity fixes its sign, and
$\lambda_p^1=\langle\mathcal A_p\psi_p^1,\psi_p^1\rangle$
is analytic as well.
\end{proof}
\subsection{Estimates for the inner linear problem}\label{app:proof-lem9}
\begin{proof}[Proof of Lemma~\ref{lem9}]

\proofstep{1}{Establish mode estimates for both trace spaces.}
Set $f_j(t)=c_p^j[f(\cdot,t)]$, $z_j=c_p^j[z]$, and
$\omega_j(t)=c_p^j[\omega(\cdot,t)]$. Angular decomposition reduces
\eqref{e112} to
    \begin{equation}\label{e116}
        \omega_j''-\alpha_p\omega_j'-\lambda_p^j\omega_j=f_j\quad \text{in}\quad (-\infty,0),\quad \omega_j(0)=z_j, \text{ for }j\geq 1,
   \end{equation}
with the additional condition $\omega'_1(0)=b$.

We first estimate these solutions for both trace spaces needed below:
$H^{3/2}$ for the strong inner inverse and $H^{1/2}_{00}$ for weak
matching to an exterior energy solution. Let
$\mu_j^\pm$ be the roots in Remark~\ref{rem:root-relation}.  Proposition~\ref{prop7}
gives, for $j\ge2$,
\[
 \begin{aligned}
 &\mu_j^-<0<\mu_j^+,\\
 &c_{n,\epstr}\sqrt{1+\lambda_p^j}
 \le \min\{-\mu_j^-,\mu_j^+\}
 \le \max\{-\mu_j^-,\mu_j^+\}
 \le C_{n,\epstr}\sqrt{1+\lambda_p^j}.
 \end{aligned}
\]
The constants are uniform on the fixed truncated interval $(p_-,p_+)$.
For $j\ge2$, the decaying half-line Green kernel and the boundary
term $e^{\mu_j^+t}z_j$ determine the bounded solution. For $j=1$, both roots
have positive real part: the data $\omega_1(0)=z_1$ and $\omega_1'(0)=b$
fix the homogeneous coefficients, and variation of constants produces a
bounded particular solution.  In the repeated case the usual
$t e^{\alpha_pt/2}$ term is used.

The scalar Green formulas are well defined for
$f\in L^2_{\rm ul}(Q_-)$, and Parseval's identity together with the
preceding parameter bounds gives, for $z\in H^{3/2}$ with zero boundary trace,
\begin{equation}\label{e116a}
 \sup_{\tau\le0}
 \|\omega\|_{W^{2,2}(S^{n-1}_+\times(\tau-1,\tau))}
 \le C
 \left(
 \|f\|_{L^2_{\rm ul}(Q_-)}
 +\|z\|_{H^{3/2}(S^{n-1}_+)}+|b|
 \right).
\end{equation}
For the weaker datum $z\in H^{1/2}_{00}(\Sigma_+)$, the same formulas imply
\begin{equation}\label{eq:hilbert-inner-critical-trace}
 \|\omega\|_{W^{1,2}(\Sigma_+\times(-1,0))}
 \le C\bigl(\|z\|_{H^{1/2}_{00}}+\|f\|_{L^2_{\rm ul}}+|b|\bigr).
\end{equation}
Indeed, for $j\ge2$ the Green kernel has amplitude and decay length
$O((1+\lambda_p^j)^{-1/2})$, hence $L^1$ size
$O((1+\lambda_p^j)^{-1})$. This compensates two angular derivatives
of the source term. The boundary layer $e^{\mu_j^+t}z_j$ contributes
$O((1+\lambda_p^j)^{3/2}|z_j|^2)$ to the squared $W^{2,2}$ norm,
and $O((1+\lambda_p^j)^{1/2}|z_j|^2)$ to the squared $W^{1,2}$ norm
near $t=0$. These are the two respective Dirichlet trace norms; the
first mode is controlled by its two data. Approximation in
$H^{1/2}_{00}$ thus constructs a $W^{1,2}_{\rm loc}(Q_-)$ weak solution,
with the prescribed first-mode derivative and projected terminal derivative
in $(H^{1/2}_{00})^*$. For $H^{3/2}$ data the construction proves
\eqref{e116a}.
Uniqueness follows directly from the homogeneous scalar
equations: for $j\ge2$, boundedness at $-\infty$ removes the
$e^{\mu_j^-t}$ mode and the value at $0$ removes the remaining mode;
for $j=1$, the two Cauchy data remove both modes.

\proofstep{2}{Establish uniformly local estimates and analytic dependence.}
Since $L^q_{\rm ul}\hookrightarrow L^2_{\rm ul}$ and
$Z\hookrightarrow H^{3/2}$, \eqref{e116a} supplies the initial
$W^{2,2}_{\rm ul}$ bound. Standard boundary $L^r$ estimates and
finite Sobolev iteration on enlarged unit strips raise it to $q$.
At $t=0$, subtract a $W^{2,q}$ extension of $z$ of norm
$C\|z\|_Z$ and reflect the zero-terminal-data remainder oddly.
The coefficients and strip geometry are uniformly controlled on
$(p_-,p_+)$, so taking strip suprema proves \eqref{e113}.
Analytic dependence follows by inverting the analytic boundary-value
operator on these fixed spaces.
\end{proof}
\subsection{Regularity and asymptotics of the outer solution}\label{app:proof-lem11}
\begin{proof}[Proof of Lemma~\ref{lem11}]

\proofstep{1}{Match the weak solutions at the interface.}
Near $\partial B_1(0)$, the energy regularity of $v$ in $\mathcal D^{1,2}_0$
ensures that $z:=v|_{S^{n-1}_+}$ belongs to $H^{1/2}_{00}(\Sigma_+)$ by
\eqref{eq:energy-critical-trace}. To obtain its $Z$-regularity from the
equation, apply Remark~\ref{rem4} to
$-P_ph$ and choose $f_h\in L^q_{\ul}(\overline{Q_-})$ so that
\[
    \partial_t\omega[p,0,f_h,0](\,\cdot\,,0)=P_ph,
    \qquad \|f_h\|_{0,q}\le C\|h\|_{Z'}.
\]
With
\[
 b:=\frac{2}{p-1}c_p^1[z]+c_p^1[h],
\]
apply \eqref{eq:hilbert-inner-critical-trace} to the data $(z,f_h,b)$.
This constructs a weak inner solution $\omega^-\in W^{1,2}_{\rm loc}(Q_-)$,
whose projected terminal derivative lies in $(H^{1/2}_{00})^*$. Define on the whole cylinder
\[
\overline\omega(\theta,t):=
\begin{dcases}
\omega^-(\theta,t),& t<0,\\
e^{\frac2{p-1}t}v(e^t\theta),& t>0.
\end{dcases}
\]
The two traces agree by construction:
\[
\overline\omega(\cdot,0^-)=\overline\omega(\cdot,0^+)=v|_{S^{n-1}_+},
\]
and \eqref{e124} matches the projected derivatives:
\[
 P_p\partial_t\overline\omega(\cdot,0^-)
 =C_p[v|_{S^{n-1}_+}]+P_ph
 =P_p\partial_\lambda v
 =P_p\partial_t\overline\omega(\cdot,0^+).
\]
For the first mode, the auxiliary condition~\eqref{e123} implies
\[
 c_p^1[\partial_\lambda v]
 =\frac{2}{p-1}c_p^1[v]+c_p^1[h],
\]
which is the parameter chosen on the inner side.  Hence both the
value and the full $t$-derivative match at $t=0$, initially in the
trace space and its dual. There is therefore no interface source, and
$\overline\omega$ is a weak
solution on $Q$ of
\[
\left\{
\begin{aligned}
\partial_t^2\overline\omega-\alpha_p\partial_t\overline\omega-
\beta_p\overline\omega+
\Delta_\theta\overline\omega+p\phi_p^{p-1}\chi_{Q_-}\overline\omega
&=\overline f &&\text{in }Q,\\
\overline\omega&=0 &&\text{on }\partial Q,
\end{aligned}
\right.
\]
where
\[
\overline f(\theta,t):=
\begin{dcases}
f_h(\theta,t),&t<0,\\
-e^{\frac{2p}{p-1}t}g(e^t\theta),&t>0.
\end{dcases}
\]

\proofstep{2}{Establish regularity up to the interface.}
The Hilbert construction bounds the local $W^{1,2}$ norm by
$C(\|\nabla v\|_2+\|h\|_{Z'})$. Since $\overline f\in L^q_{\rm loc}$
and the potential is bounded, the boundary estimates and finite
Sobolev iteration used in Appendix~\ref{app:proof-lem9} give
\[
 \|\overline\omega\|_{W^{2,q}(S^{n-1}_+\times(-1,1))}
 \le C\bigl(|g|_{0,q;\nu}+\|h\|_{Z'}+\|\nabla v\|_2\bigr).
\]
Transforming back and taking the trace now implies $v|_{S^{n-1}_+}\in Z$.

\proofstep{3}{Determine the far-field coefficient and estimate the remainder.}
Choose a smooth cutoff $\chi=0$ on $B_2$ and $\chi=1$ outside $B_3$. Then
$\widetilde v=\chi v$ satisfies
\[
 -\Delta\widetilde v
 =\chi g-2\nabla\chi\cdot\nabla v-v\Delta\chi=:\widetilde g,
 \qquad \widetilde v|_{\partial\R^n_+}=0.
\]
For $D=|g|_{0,q;\nu}+\|h\|_{Z'}+\|\nabla v\|_2$, local regularity
controls the compact commutators, while $\nu<-n-1$ makes the tail moment
summable. Hence
\[
 M(\widetilde g)\le CD,\qquad
 \|\widetilde g\|_{L^q(A_r^+)}\le CD r^{\nu+n/q}\quad(r\ge6).
\]
Since $\widetilde g\in(\mathcal D^{1,2}_0)^*$, its Green potential
belongs to $\mathcal D^{1,2}_0$. The harmonic difference
$\widetilde v-G[\widetilde g]$ has zero Dirichlet trace; testing by
cutoff approximations gives zero Dirichlet energy, hence
$\widetilde v=G[\widetilde g]$. Lemma~\ref{lem7}, together with the
local bound, now proves \eqref{e126} and
\[
\begin{gathered}
 \langle I,v\rangle=\frac{2}{n|B_1|}\int y_n\widetilde g(y)\,dy,\\
 \sup_{|x|=r,\ x_n>0}
 \left|\frac{r^nv(x)}{x_n}-\langle I,v\rangle\right|
 \le CD r^{-\kappa},\qquad r\ge6,
 \quad\kappa=\tfrac12\min\{1,-\nu-n-1\}.
\end{gathered}
\]
For arbitrary $v\in W$, the same cutoff representation with
$g=-\Delta v$ proves $I\in W^*$ and a remainder uniform on bounded
subsets of the entire graph space, since
\[
 M(\widetilde g)+\sup_{r\ge6}r^{-\nu-n/q}
       \|\widetilde g\|_{L^q(A_r^+)}+|\langle I,v\rangle|
 \le C\|v\|_W.
\]

    \end{proof}

\section{Weighted regularity of the canonical inverse}\label{app:fixed-regularity}
To use the canonical Hilbert inverse in Section~\ref{s8}, we need
its weighted $W^{2,q}$ regularity. We prove this for the Hilbert solution
and its paired variable, first on the homogeneous model and then by
localization. All constants may depend on the fixed parameter $\ve$.

Use the geometry $(\Omega,S)$ of Section~\ref{s6}, with $N=n+k$,
the function $\rho$ fixed in Section~\ref{s8}, and the norms of
Section~\ref{sec:global-weighted-spaces}, replacing
$q$ by any $1<s<\infty$. The case $k=0$ is included. As there,
$\mathcal H_b=L^2(\rho^{-2b-n}dx)$. Explicitly,
\[
 C^0_b:=\{w\in C^0(\overline\Omega\setminus S):
       \|w\|_{C^0_b}:=\sup\rho^{-b}|w|<\infty\}.
\]
For $j\ge0$ and $0<\beta<1$, write
$\mathcal C^{j,\beta}_b=C^{j,\beta}_b$ for the scaled H\"older
space in~\eqref{eq:scaled-holder-spaces}, and
$\mathscr H^2_b=\mathscr H^2_{b,0}$ for the zero-trace edge space in
\eqref{eq:global-edge-Sobolev-space}. On the full normal model, these
spaces use $\rho=r=|X|$ and the physical Dirichlet boundary.
Norms on pairs are the sums of the two component norms.

Fix $\varepsilon$ and put $L=\Delta+V_\varepsilon$, where
$V_\varepsilon=p\bar u_\varepsilon^{p-1}$ is the potential of
Section~\ref{s8}. Its scaled Fermi-coordinate limit is
$V_0=p\phi_p^{p-1}$, uniformly on fixed normal half-annuli as
$\rho/\varepsilon\to0$. No rate uniform in $\varepsilon$ is assumed.
The homogeneous model estimates below also include $V_0=0$.

Use the weight $\delta$ in \eqref{e59} and choose $a=\delta_-$
sufficiently close to $\delta$ that
\[
 -a_p<\check\delta:=2a-\delta<a<\delta
 <\operatorname{Re}\gamma_1^-,\qquad
 I=[\check\delta,\delta],\qquad e=2a+n-2.
\]
Proposition~\ref{prop7} gives, for the angular operator
$-\Delta_\theta-p\phi_p^{p-1}$, eigenvalues
$\widetilde\lambda_j=\lambda_p^j-\beta_p$ and
\[
 \gamma_2^-=-a_p-1,\qquad \gamma_2^+=a_p+3-n.
\]
The singular radial branch in angular mode $j$ is admitted by the
inner Hilbert weight $a$ exactly when $\operatorname{Re}\gamma_j^->a$. Here
\[
 \operatorname{Re}\gamma_1^->\delta\ge a,\qquad
 \gamma_j^-\le-a_p-1<\check\delta\quad(j\ge2).
\]
Thus $I$ is root-free, and only the first mode has an admitted singular
branch. For the harmonic model, every lower root is at most
$1-n<\check\delta$, so no singular branch is admitted.
Set
\[
 \delta_+=2-n-a,\qquad \delta_*=2-n-\delta,\qquad
 \widehat\delta=\delta-e,\qquad J_+=[\delta_*,\widehat\delta].
\]
Indicial symmetry makes $J_+$ root-free. These are precisely the
weights used in Proposition~\ref{prop5}.

\begin{theorem}\label{thm:fixed}
For each fixed $\varepsilon$, the following hold.

If $z\in\mathcal D_{\delta_+}$, $Lz=f\in L^s_{\sigma-2}$ and
$\sigma\in J_+$, then
\[
 \|z\|_{W^{2,s}_\sigma}
 \le C_{\varepsilon,s,\sigma}
       \bigl(\|f\|_{L^s_{\sigma-2}}+\|z\|_{\mathcal D_{\delta_+}}\bigr).
\]
In particular this applies to the high Hilbert kernel, with $f=0$.
The high Hilbert realization has closed range and finite-dimensional kernel.

Suppose the low realization
$A:\mathcal D_a\subset\mathcal H_a\to\mathcal H_{a-2}$ is surjective. For $q>N$,
its existing canonical inverse $G^H=A^*(AA^*)^{-1}$ satisfies
\[
 \|G^Hh\|_{W^{2,q}_\delta}\le C_\varepsilon\|h\|_{L^q_{\delta-2}}.
\]
The paired variable
$z=\rho^{-2a+4-n}(AA^*)^{-1}h$ simultaneously belongs to
$W^{2,q}_{\widehat\delta,0}$.
\end{theorem}

\subsection{The Hilbert adjoint and the paired system}
Put $P=\rho^2L$ on $\mathcal H_a$. Multiplication by $\rho^2$ is unitary
from $\mathcal H_{a-2}$ to $\mathcal H_a$, and
\[
 P^*=\rho^{2a+n}L\rho^{-2a-n+2}=\rho^eP\rho^{-e}.
\]
Both domains are $\mathscr H^2_a$: distributional adjoint tests identify
the operator and its zero physical trace, and scaled Dirichlet estimates give
\[
 \|w\|_{\mathscr H^2_a}\le C(\|Tw\|_{\mathcal H_a}+\|w\|_{\mathcal H_a}),
 \qquad T=P,P^*.
\]
Radial cutoff commutators tend to zero with the weighted derivative tails;
Dirichlet approximation on compact annuli then gives a common graph core.
Thus these are closed adjoint realizations with no extra pairing at $S$.
The same argument applies to the models.

For a densely defined closed surjection $P:D(P)\subset H\to H$, put
$G=P^\dagger$ and let $\Pi$ project orthogonally onto $\ker P$. Then
\begin{equation}\label{appB:hilbert}
 \mathcal M=\begin{pmatrix}P&0\\-I&P^*\end{pmatrix},\qquad
 \mathcal M^{-1}=\begin{pmatrix}G&-\Pi\\G^*G&G^*\end{pmatrix}.
\end{equation}
Here $D(\mathcal M)=D(P)\times D(P^*)$. The identities $PG=G^*P^*=I$ and
$GP=P^*G^*=I-\Pi$, with the orthogonality of the ranges, verify both
products on these domains. In particular, data $(F,0)$ give the
original canonical solution $GF$ in the first component.

The local estimates used below are standard Dirichlet estimates on
nested rescaled patches $K_1\Subset K_2$: for the scalar operators
and pairs,
\begin{equation}\label{appB:local}
 \|w\|_{W^{2,s}(K_1)}\le C_s
 \bigl(\|Tw\|_{L^s(K_2)}+\|w\|_{L^2(K_2)}\bigr),\qquad 1<s<\infty.
\end{equation}
Finite Sobolev iteration from $L^2$ extends the estimate to every
finite $s$ because
the potential and coupling are bounded zeroth-order terms. For $s>N$, boundary $C^1$ compactness follows as well.

\subsection{The Fourier and cylinder operators}
Let $\Sigma=S^{n-1}_+$, $\mathcal A=-\Delta_\Sigma-V_0$ with Dirichlet
domain, and write its eigenvalues as $\lambda_j$. On the model,
\[
 P_0=D^2+(n-2)D+r^2\Delta_Y-\mathcal A,\qquad D=r\partial_r,
 \qquad \mathcal M_0=\begin{pmatrix}P_0&0\\-I&r^eP_0r^{-e}\end{pmatrix}.
\]
The scalar roots satisfy
$d_j(\zeta)=\zeta^2+(n-2)\zeta-\lambda_j=0$.
The indicial matrix is obtained by inserting $r^\zeta$ into the
homogeneous radial terms. For the pair it is
\[
 \mathfrak I_j(\zeta)
 :=\begin{pmatrix}d_j(\zeta)&0\\-1&d_j(\zeta-e)\end{pmatrix},\qquad
 \det\mathfrak I_j(\zeta)=d_j(\zeta)d_j(\zeta-e).
\]
Its zeros define the paired indicial spectrum. Since
$e=2a+n-2$ and $\gamma_j^++\gamma_j^-=2-n$,
\[
 e+\gamma_j^\pm=2a-\gamma_j^\mp,\qquad
 \operatorname{Spec}_{\rm ind}\mathcal M_0
 =\bigcup_j\{\gamma_j^\pm,\,2a-\gamma_j^\pm\}.
\]
Consequently, with
\begin{equation}\label{appB:gap}
 \iota=\delta-a,\qquad \tau_*=\min_{j,\pm}|\operatorname{Re}\gamma_j^\pm-a|,
 \qquad \iota<\tau_*,
\end{equation}
the interval $I=[a-\iota,a+\iota]$ avoids every paired root.
Moreover $2\iota<\delta+a_p<1$ by the choice of weights.

\smallskip
\noindent\emph{Graph domains at nonzero frequencies.}
For $\kappa=|\xi|>0$ put $B_\xi=\Delta_X-\kappa^2+r^{-2}V_0$ and
$P_\xi=r^2B_\xi$ on $H_a^X=L^2(r^{-2a-n}dX)$.
Its graph domain is precisely the realization $\mathcal D_{a,a-2}(\xi)$
of Lemma~\ref{lem:scalar-massive}, with equivalent graph norm after multiplying the
source by $r^2$. To check the exterior part of this assertion, for
$\xi$ in a fixed compact set $K\Subset\mathbb R^k\setminus0$ one has
\begin{equation}\label{appB:mass}
 \begin{aligned}
 &\|w\|_{H_{a-2}^X(r>2R)}+
 \|\nabla w\|_{H_{a-2}^X(r>2R)}+
 \|\nabla^2w\|_{H_{a-2}^X(r>2R)}
 \\
 &\le C_K\bigl(\|B_\xi w\|_{H_{a-2}^X(r>R)}+
                  \|w\|_{H_a^X(R<r<3R)}\bigr).
 \end{aligned}
\end{equation}
This estimate initially requires only $w\in H_a^X$ and the displayed
source condition. Conjugate by $r^{-a-n/2}$ and multiply the conjugated
function by $m_T(r)=r^2/(1+r^2/T^2)$. Its logarithmic derivative is
bounded by $2/r$, uniformly in $T$. Outside a sufficiently large fixed
$R$, the induced first-order terms $O(r^{-1})$ and zeroth-order terms
$O(r^{-2})$ are absorbed by $\kappa^2\ge\min_K|\xi|^2>0$ in the
Dirichlet energy estimate. The source is bounded by its $r^2$-weighted
norm and the inner cutoff terms are on a fixed annulus. Let $T\to\infty$;
unit-scale Dirichlet estimates imply the two derivative bounds in
\eqref{appB:mass}. This also proves the reverse inclusion of graph domains.
The singular-end inclusion follows from scaled local estimates.

In particular $r^2:D(P_\xi)\to H_a^X$ is bounded, uniformly on $K$.
Conjugating by $r^e$ identifies the adjoint with the scalar realization
$\mathcal D_{a-e,a-e-2}(\xi)$ and preserves the same normalized domain.
Thus both domains are fixed on $K$, and
\[
 P_\xi-P_{\xi_0}=-(|\xi|^2-|\xi_0|^2)r^2
\]
is a bounded map from that domain to $H_a^X$.

Lemma~\ref{lem:scalar-massive} identifies the range annihilator of the
low scalar realization with the high kernel at inner weight $2-n-a$. This kernel is zero.
Indeed, dyadic local estimates turn the graph bounds of a kernel
element into $O(r^{2-n-a})$ at zero and a polynomial bound at
infinity, and $2-n-a>\operatorname{Re}\gamma_1^+$;
Lemma~\ref{lem:conic-modes}\,\textup{(ii)} applies.
Thus $P_\xi$ is surjective. Formula~\eqref{appB:hilbert} provides its paired
inverse. Smooth dependence on $\xi$ and bounds for every finite-order
derivative on $K$ follow from the fixed-domain inverse identity.
The unitary radial dilation $w(r)\mapsto\kappa^{-a}w(\kappa r)$
identifies $P_\xi$ with $P_{|\xi|=1}$. Direct integration in $\xi$ therefore
constructs the Hilbert inverse $\mathcal B_0=\mathcal M_0^{-1}$.
Local $L^2$ estimates identify its domain with the full model edge graph
space.

\smallskip
\noindent\emph{The cylinder inverse on different weighted spaces.}
To improve the weight of this Hilbert solution at $r=0$, put $t=\log r$
and conjugate both components by $r^a$. The limiting cylinder operator is
\begin{equation}\label{appB:Q}
 \mathcal Q=\begin{pmatrix}\mathcal P_+&0\\-I&\mathcal P_-\end{pmatrix},
 \qquad \mathcal P_\pm=\partial_t^2\pm e\partial_t+a(a+n-2)-\mathcal A.
\end{equation}
The exponential weights define cylinder spaces; that is,
\[
\begin{aligned}
 L^2_\sigma&:=\{f:e^{-\sigma t}f\in L^2(\R\times\Sigma)\},\\
 H^2_\sigma&:=e^{\sigma t}\{v\in H^2(\R\times\Sigma):
                                  v|_{\R\times\partial\Sigma}=0\}.
\end{aligned}
\]
Their norms are $\|f\|_{L^2_\sigma}:=\|e^{-\sigma t}f\|_2$ and
$\|w\|_{H^2_\sigma}:=\|e^{-\sigma t}w\|_{H^2}$.
For $|\sigma|<\tau_*$, Fourier transformation in $t$ and angular
decomposition construct an inverse $R_\sigma:(L^2_\sigma)^2\to(H^2_\sigma)^2$.
Indeed, in mode $j$ its symbol is the inverse of
\[
 \begin{pmatrix}q_{+,j}(z)&0\\-1&q_{-,j}(z)\end{pmatrix},\qquad
 q_{\pm,j}(z)=z^2\pm ez+a(a+n-2)-\lambda_j,\quad z=\sigma+i\eta.
\]
The inverse entries are $q_{+,j}^{-1}$, $q_{-,j}^{-1}$, and
$(q_{+,j}q_{-,j})^{-1}$. There are no poles in $|\operatorname{Re}z|<\tau_*$.
On compact substrips the diagonal entries are bounded by
$C(1+\eta^2+\lambda_j+C_0)^{-1}$; the product has the corresponding
squared bound. High modes follow from the negative real quadratic part,
and only finitely many remaining modes need the root-free condition.
The Dirichlet graph norm of $\mathcal A$ is equivalent to $H^2$ since
$V_0$ is bounded. These bounds prove the asserted cylinder estimate.
For compactly supported smooth data, shifting the Fourier--Laplace
contour between two such lines crosses no pole. Approximation by common
compact truncations in the two weighted $L^2$ norms shows
\begin{equation}\label{appB:consistent}
 R_{\sigma_1}g=R_{\sigma_2}g
 \quad\text{if }g\in (L^2_{\sigma_1}\cap L^2_{\sigma_2})^2.
\end{equation}
For $k=0$, $\mathcal Q$ is the entire transformed model, so $R_0$
provides its Hilbert inverse.

\subsection{Improved decay at nonzero frequencies}
\begin{lemma}\label{appB:F}
For $k\ge1$, suppose
$\widehat F\in C_c^\infty((\R^k\setminus\{0\})\times\R^n_+;\C^2)$.
For $k=0$, suppose $F\in C_c^\infty(\R^n_+;\C^2)$.
Then the Hilbert solution $w=\mathcal B_0F$ belongs to $C^0_b$ for every
$b\in I$. In the first case no compact support in $Y$ is required.
\end{lemma}
\begin{proof}
\proofstep{1}{Improve decay at the singular end.}
Choose $\iota<\sigma<\min\{\tau_*,2\}$, which is possible by
\eqref{appB:gap}. For each frequency, $W=r^{-a}\widehat w$ satisfies
$(\mathcal Q-|\xi|^2e^{2t})W=r^{-a}\widehat F$.
Choose $\chi_-(t)$ supported below the data, equal to one sufficiently
far to the left. The Hilbert graph regularity implies $W_-:=\chi_-W\in H^2_0$,
and
\begin{equation}\label{appB:boost}
 \mathcal QW_-=[\mathcal Q,\chi_-]W+|\xi|^2e^{2t}W_-=:g.
\end{equation}
The commutator is supported on a fixed strip. The other term belongs
to $L^2_0\cap L^2_\sigma$, since $W_-\in L^2_0$ and
$e^{(2-\sigma)t}$ is bounded on its support. By
\eqref{appB:consistent},
\[
W_-=R_0g=R_\sigma g\in H^2_\sigma.
\]
Since $R_\sigma$ inverts the full matrix $\mathcal Q$, both components
gain the weight $\sigma$.

\proofstep{2}{Estimate frequency derivatives and decay at infinity.}
All frequency derivatives first belong to the fixed Hilbert graph
domain by the inverse identity and the bound for $r^2$. Differentiating
\eqref{appB:boost} adds only terms $e^{2t}\partial_\xi^\alpha W_-$
of lower frequency order, so the same argument establishes $H^2_\sigma$
for every finite derivative, uniformly on $K$. On unit cylinder strips,
\eqref{appB:local} and finite bootstrapping give
\begin{equation}\label{appB:ends}
 |\partial_\xi^\alpha\widehat w|+
 r|\nabla_X\partial_\xi^\alpha\widehat w|
 \le C_\alpha r^{a+\sigma}\quad(0<r\le1).
\end{equation}
Only boundedness of the angular potential is needed for this step.

For decay at infinity, apply the energy estimate in
Lemma~\ref{lem:scalar-massive} first to $u$, then to $z=r^{-e}v$,
whose equation has source $f=r^{-2a-n}u$. Use a smaller exponential
rate for $z$: the additional term
$-\operatorname{Re}\int q_T^2 f\overline z$ is absorbed by
Cauchy--Schwarz because $q_Tf\in L^2$. Frequency differentiation adds
only lower-order derivatives already estimated, so induction gives
exponential decay, with first normal derivatives, at every finite
frequency order on $K$.

\proofstep{3}{Recover the physical-space solution.}
As in the proof of Lemma~\ref{lem:fourier-testing}, integration by parts
in frequency gives rapid decay in $Y$. The two normal-end estimates
then imply $w\in C^0_b$, since $a+\sigma>\delta\ge b$.
For $k=0$, use $R_0=R_\sigma=R_{-\sigma}$ on annular data;
the bounds $r^{a+\sigma}$ at zero and $r^{a-\sigma}$ at infinity
give the same conclusion for every $b\in I$.
\end{proof}

\subsection{Liouville theorem and the supremum estimate}
\begin{lemma}\label{appB:L}
Let $b\in I$ and $q>N$. If
$w\in W^{2,q}_{\rm loc}(\mathcal O_k;\C^2)$ has zero physical trace,
$\mathcal M_0w=0$ distributionally, and $\|w\|_{C^0_b}<\infty$, then
$w=0$.
\end{lemma}
\begin{proof}
\proofstep{1}{Construct test functions for the transposed system.}
Write $w=(u,v)$, $z=r^{-e}v$ and $\varpi=r^{-2a-n}$. The equations are
$\mathcal Lu=0$, $\mathcal Lz=\varpi u$, where
$\mathcal L=\Delta_{X,Y}+r^{-2}V_0$. Fix a compact nonzero frequency
set. Choose $\vartheta>0$ so that
$s=b-\vartheta$ and $\ell=2a-b+\vartheta$ lie in the same slightly
larger root-free interval containing $I$, still below $-(n-2)/2$.
For arbitrary annular test data smooth in frequency, solve
\begin{equation}\label{appB:transpose}
 B_\xi\eta-\varpi\zeta=\psi_1,\qquad B_\xi\zeta=\psi_2,
 \qquad
 \eta\in\mathcal D_{2-n-s,t_1},\quad
 \zeta\in\mathcal D_{\ell,t_2},
\end{equation}
where $t_2=t_1+2a+n-2$. The source spaces are defined in Lemma~\ref{lem:scalar-massive}.
Multiplication by $\varpi$ matches the inner source weight exactly,
since $\ell-2a-n=-n-s=(2-n-s)-2$, and gains a factor $r^{-2}$ at infinity.

Lemma~\ref{lem:conic-modes}\,\textup{(ii)} gives high injectivity,
and Lemma~\ref{lem:scalar-massive} gives low surjectivity and identifies
its kernel with the high range annihilator: $s,\ell$ lie in the same
root-free gap, and exterior polynomial weights do not change the kernel.
The root bounds preceding Theorem~\ref{thm:fixed}, together with
$\operatorname{Re}\gamma_j^+\ge\operatorname{Re}\gamma_1^-$, place every admitted
inner exponent strictly above $\delta>a$, including repeated roots
with their logarithmic factors. Hence
\[
 \mathfrak g(\kappa,\psi)=\int_{\mathbb R^n_+}
             \varpi\kappa\overline\psi\,dX
\]
is finite and positive definite on this kernel, even for $b\le a$.
Lemma~\ref{lem:gram-correction} therefore supplies the transpose inverse
and all finite-order frequency derivatives on compact frequency sets.

\proofstep{2}{Use Green identities to prove that the solution vanishes.}
Local regularity and the graph weights give, for every finite frequency
derivative, the inner bounds
\[
 \eta=O(r^{2-n-b+\vartheta}),\qquad
 \zeta=O(r^{2a-b+\vartheta}),
\]
also for first scaled normal derivatives. The compact data and positive
mass ensure exponential outer decay, as in Lemma~\ref{appB:F}.
The two inner Green exponents are
\[
 b+(2-n-b+\vartheta)+n-2=\vartheta,
 \qquad (b-e)+(2a-b+\vartheta)+n-2=\vartheta,
\]
and the coupling is bounded by $C\int_0^1r^{\vartheta-1}\,dr$.
Lemma~\ref{lem:fourier-testing} puts both Fourier supports at zero;
their uniform $Y$ bounds make the components independent of $Y$.

In each angular mode, $t=\log r$ reduces the remaining equation to a
constant-coefficient system with exponents
$\gamma_j^\pm$ and $2a-\gamma_j^\pm$. None has real part $b$, so
no nonzero exponential polynomial is $O(e^{bt})$ at both ends.
This proves the result and also treats $k=0$ directly.
\end{proof}

Fix $q>N$. Uniformly for $b\in I$, every zero-trace solution
$w\in C^0_b$ with $\mathcal M_0w\in L^q_b$ satisfies
\begin{equation}\label{appB:A}
 \|w\|_{C^0_b}\le C_q\|\mathcal M_0w\|_{L^q_b}.
\end{equation}
Otherwise normalize the left side to one, let $b_j\to b_\infty$, and
choose $(r_j\theta_j,Y_j)$ with weighted value at least $1/2$. The rescaling
\[
 \widetilde w_j(X,Y)=r_j^{-b_j}w_j(r_jX,Y_j+r_jY)
\]
preserves the equations and the scaled $L^q_{b_j}$ source norm, which
tends to zero. Estimate~\eqref{appB:local} gives local $W^{2,q}$ bounds
and $C^1$ compactness up to the physical boundary. The zero trace thus
keeps $\theta_j$ away from the equator. The limit is nonzero, solves
the homogeneous system, and is bounded by $r^{b_\infty}$, contrary to
Lemma~\ref{appB:L}. The source weight here is $b$, since the equations
have been multiplied by $r^2$.

\subsection{Weighted bounds for the Hilbert inverse}
For smooth data supported in a normal annulus and compactly supported
in $Y$, cut off the tangential frequencies by
\[
 \widehat F_\eta(X,\xi)=
 (1-\chi(\xi/\eta))\chi(\eta\xi)\widehat F(X,\xi),\qquad 0<\eta<1,
\]
where $\chi\in C_c^\infty$ equals one near zero. The convolution
kernels have bounded $L^1$ norms. Minkowski's inequality on the scaled
balls therefore gives $\|F_\eta\|_{L^q_b}\le C\|F\|_{L^q_b}$;
Plancherel gives convergence in the Hilbert space. By
Lemma~\ref{appB:F}, $\mathcal B_0F_\eta\in C^0_b$, so
\eqref{appB:A} and the scaled local estimate~\eqref{appB:local} yield
\begin{equation}\label{appB:C0}
 \|\mathcal B_0F\|_{W^{2,q}_b}\le C_q\|F\|_{L^q_b},
 \qquad q>N,\quad b\in I.
\end{equation}
The limit is the original Hilbert solution by Hilbert convergence.
For $k=0$, Lemma~\ref{appB:F} applies without a frequency cutoff.
For general Hilbert data in $L^q_b$, compact cutoffs and smoothing
preserve bounded $L^q_b$ norms and converge in Hilbert norm. For
arbitrary $L^q_b$ data, the same bounded approximations converge locally
in $L^q$; local compactness and Lemma~\ref{appB:L} give a unique
extension, agreeing with the Hilbert inverse on the intersection.
Formula~\eqref{appB:hilbert} also gives this bound for
$S_0=G_0^*$, the Hilbert left inverse of $P_0^*$.

For the full range $1<s<\infty$ in Theorem~\ref{thm:fixed}, a finite
local regularization suffices. For $T_0=P_0^*$ or $\mathcal M_0$, patch
Dirichlet inverses of its diagonal part minus a fixed large constant
on Whitney balls, retaining the bounded potential. The resulting $Q$
and $E_0=I-T_0Q$ satisfy
\[
 Q:L^s_b\to W^{2,s}_{b,0},\qquad E_0:L^s_b\to L^t_b,
 \quad t>s,\quad \frac1t>\max\{\tfrac1s-\tfrac1N,0\}.
\]
The defect consists only of first-order cutoff commutators and bounded
zeroth-order terms. For the Hilbert left inverse $B_0=S_0$ or
$\mathcal B_0$, iterate $B_0=Q+B_0E_0$ until the remaining data lie
in some $L^q_b$, $q>\max\{N,s\}$, and apply~\eqref{appB:C0} and
local regularity. Thus
\begin{equation}\label{appB:Ls}
 B_0:L^s_b\longrightarrow W^{2,s}_{b,0},
 \qquad 1<s<\infty,\quad b\in I.
\end{equation}
Bounded overlap gives the same finite identity on the Hilbert spaces,
so this extension still agrees with the specified Hilbert operator.

\subsection{Localization and regularity of the canonical pair}

\proofstep{1}{Establish weighted regularity of the Hilbert solution.}
Fix $\varepsilon$ and let $w$ be the Hilbert graph solution of
$Tw=F\in L^s_b$, where $T=P^*$ or $\mathcal M$ and $b\in I$.
We transfer \eqref{appB:Ls} to this solution using a finite tangential atlas of diameter $\ell$, followed
by a radial cutoff at $r_0$. For each frozen homogeneous model, use
the left inverse \eqref{appB:Ls}. The coefficient differences are
bounded by
\[
 C\ell+Cr_0+\vartheta(r_0/\varepsilon),\qquad \vartheta(t)\to0\ (t\to0),
\]
and tangential cutoff commutators have coefficients
$O(r_0/\ell)+O((r_0/\ell)^2)$. Radial commutators are supported
away from $S$. The exact localized Hilbert identities imply
\begin{equation}\label{appB:localization}
 w=\mathcal BF+Ew+Kw,\qquad
 \|E\|_{\mathscr H^2_a\to\mathscr H^2_a}<\tfrac12,
 \quad \|E\|_{W^{2,s}_b\to W^{2,s}_b}<\tfrac12,
\end{equation}
by choosing $\ell$ and then $r_0$ sufficiently small.
The term $K$ contains the exterior part and model inverses applied
to radial commutators. On their supporting annuli, local regularity
uses the known Hilbert norm and $F\in L^s_b$. Consequently
\[
 R:=\mathcal BF+Kw\in\mathscr H^2_a\cap W^{2,s}_{b,0},
 \qquad \|R\|_{W^{2,s}_b}\le
 C(\|F\|_{L^s_b}+\|w\|_{\mathscr H^2_a}).
\]
By \eqref{appB:localization},
\[
w=R+Ew,\qquad
\sum_{j\ge0}E^jR\ \hbox{converges in both }\mathscr H^2_a
\hbox{ and }W^{2,s}_{b,0}.
\]
The two limits agree as distributions on every compact subset away from
$S$. Convergence in the Hilbert space identifies their common limit
with $w$, so
\[
w=\sum_{j\ge0}E^jR\in W^{2,s}_{b,0}.
\]
This proves membership and the estimate
\begin{equation}\label{appB:fixed}
 \|w\|_{W^{2,s}_b}\le
 C(\|Tw\|_{L^s_b}+\|w\|_{\mathscr H^2_a}).
\end{equation}

\proofstep{2}{Establish closed range and finite-dimensional kernel at the higher weight.}
For the scalar adjoint, the same Hilbert localization and local
interpolation on the commutator annuli give
\begin{equation}\label{appB:compact}
 \|w\|_{\mathscr H^2_a}\le
 C(\|P^*w\|_{\mathcal H_a}+\|w\|_{L^2(K_0)}),
\end{equation}
where $K_0$ is compact and separated from $S$.
Restriction to $K_0$ is compact, so $P^*$ has closed range and
finite-dimensional kernel.

For the high realization set $v=\rho^ez$. Then
$P^*v=\rho^{e+2}f$, and $b=e+\sigma\in I$ for $\sigma\in J_+$.
The powers identify the Hilbert graphs and the target norms.
Equations \eqref{appB:fixed}--\eqref{appB:compact} prove the first part
of Theorem~\ref{thm:fixed}.

\proofstep{3}{Establish weighted estimates for the canonical pair.}
For the canonical pair, \eqref{eq:Lq-Hilbert} with $a=\delta_-$ implies
\begin{equation}\label{appB:embedding}
 \|h\|_{\mathcal H_{a-2}}\le C\|h\|_{L^q_{\delta-2}}.
\end{equation}
Assume the stated surjectivity of $A$. With $U h=\rho^2h=F$, take
the canonical Hilbert pair
\[
 u=P^\dagger F=G^Hh,\qquad v=(PP^*)^{-1}F,
 \qquad Pu=F,\quad P^*v=u.
\]
Its Hilbert graph norm is bounded by $C_\varepsilon\|h\|_{L^q_{\delta-2}}$.
Apply \eqref{appB:fixed} at $b=\delta$, $s=q$ to obtain
\[
 \|u\|_{W^{2,q}_\delta}+\|v\|_{W^{2,q}_\delta}
 \le C_\varepsilon\|h\|_{L^q_{\delta-2}}.
\]
Finally $PP^*=UAA^*U^{-1}$, and hence
\[
 \rho^{-e}v=\rho^{-2a+4-n}(AA^*)^{-1}h=z
 \in W^{2,q}_{\delta-e,0}=W^{2,q}_{\widehat\delta,0}.
\]
These estimates prove the second part of Theorem~\ref{thm:fixed}
for the canonical pair chosen above, with $\varepsilon$ fixed.

\section*{Acknowledgements}
The research of W. Ao is supported by the National Natural Science Foundation of China (No.~12471111) and the National Key Research and Development Program of China (No.~2022YFA1006800). The research of S. Katayama is supported by JSPS KAKENHI Grant Number 26KJ0296. The research of J. Wei is partially supported by the Hong Kong General Research Fund of the Research Grants Council of Hong Kong under the project ``New frontiers in singularity formations in nonlinear partial differential equations.''

\begin{samepage}
\section*{Declaration on the use of generative AI}
Generative AI tools were used to assist with language editing, consistency
checking, and the
organization of the manuscript. The authors are responsible for the mathematical content
and conclusions of the article.
\par\end{samepage}

\par\medskip\noindent\begin{minipage}{\textwidth}
\section*{Author information}
\noindent
Weiwei Ao, School of Mathematics and Statistics, Wuhan University,
Wuhan 430072, China.  \texttt{wwao@whu.edu.cn}

\medskip
\noindent
Sho Katayama, Department of Mathematics, Faculty of Science and Technology,
Keio University, 3-14-1 Hiyoshi, Kohoku-ku, Yokohama 223-8522, Japan.
\texttt{sho-katayama@keio.jp}

\medskip
\noindent
Juncheng Wei, Department of Mathematics, The Chinese University of Hong Kong,
Shatin, N.T., Hong Kong.  \texttt{wei@math.cuhk.edu.hk}

\end{minipage}


\begin{thebibliography}{99}

\bibitem{acdfgw}
W. Ao, H. Chan, A. DelaTorre, M. A. Fontelos, M. del Mar Gonz\'alez, and J. Wei, On higher-dimensional singularities for the fractional Yamabe problem: a nonlocal Mazzeo--Pacard program, Duke Math. J. \textbf{168} (2019), no. 17, 3297--3411.

\bibitem{acdfgw-ode}
W. Ao, H. Chan, A. DelaTorre, M. A. Fontelos, M. del Mar Gonz\'alez, and J. Wei, ODE methods in non-local equations, J. Math. Study \textbf{53} (2020), no. 4, 370--401.

\bibitem{bdn}
H. Berestycki, I. Capuzzo Dolcetta, and L. Nirenberg, Superlinear indefinite elliptic problems and nonlinear Liouville theorems, Topol. Methods Nonlinear Anal. \textbf{4} (1994), no. 1, 59--78.

\bibitem{bmn}
M. Bhakta, M. Marcus, and P.-T. Nguyen, Boundary value problems for semilinear Schr\"odinger equations with singular potentials and measure data, Math. Ann. \textbf{390} (2024), no. 1, 351--379.

\bibitem{bvhnv}
M.-F. Bidaut-V\'eron, G. Hoang, Q.-H. Nguyen, and L. V\'eron, An elliptic semilinear equation with source term and boundary measure data: the supercritical case, J. Funct. Anal. \textbf{269} (2015), no. 7, 1995--2017.

\bibitem{bpv1}
M.-F. Bidaut-V\'eron, A. C. Ponce, and L. V\'eron, Isolated boundary singularities of semilinear elliptic equations, Calc. Var. Partial Differential Equations \textbf{40} (2011), no. 1--2, 183--221.

\bibitem{bvv}
M.-F. Bidaut-V\'eron and L. Vivier, An elliptic semilinear equation with source term involving boundary measures: the subcritical case, Rev. Mat. Iberoamericana \textbf{16} (2000), no. 3, 477--513.

\bibitem{cl}
C.-C. Chen and C.-S. Lin, Existence of positive weak solutions with a prescribed singular set of semilinear elliptic equations, J. Geom. Anal. \textbf{9} (1999), no. 2, 221--246.

\bibitem{cww-bif}
R. M. Chen, S. Walsh, and M. H. Wheeler, Existence and qualitative theory
for stratified solitary water waves, Ann. Inst. H. Poincar\'e Anal.
Non Lin\'eaire \textbf{35} (2018), no.~2, 517--576.

\bibitem{jwx}
B. Choi, X. Jiang, H.-Y. Wang, and J. Xiong, Optimal weighted smoothing and asymptotics of ancient solutions for fast diffusion equations, arXiv:2605.14622v2, 2026.

\bibitem{dmp}
M. del Pino, M. Musso, and F. Pacard, Boundary singularities for weak solutions of semilinear elliptic problems, J. Funct. Anal. \textbf{253} (2007), no. 1, 241--272.

\bibitem{fik}
M. Fila, K. Ishige, and T. Kawakami, Positive solutions of a semilinear elliptic equation with singular Dirichlet boundary data, J. Elliptic Parabol. Equ. \textbf{1} (2015), no. 2, 335--362.

\bibitem{gh-bif}
C. Garc\'ia and S. V. Haziot, Global bifurcation for corotating and
counter-rotating vortex pairs, Comm. Math. Phys. \textbf{402} (2023),
1167--1204.

\bibitem{gs}
B. Gidas and J. Spruck, Global and local behavior of positive solutions of nonlinear elliptic equations, Comm. Pure Appl. Math. \textbf{34} (1981), no. 4, 525--598.

\bibitem{gk}
K. T. Gkikas, Boundary singularities on a wedge-like domain of a semilinear elliptic equation, Proc. Roy. Soc. Edinburgh Sect. A \textbf{145} (2015), no. 5, 979--1006.

\bibitem{gkn}
K. T. Gkikas and P.-T. Nguyen, Semilinear elliptic equations involving power nonlinearities and Hardy potentials with boundary singularities, Proc. Roy. Soc. Edinburgh Sect. A \textbf{155} (2025), no. 3, 993--1050.

\bibitem{hmr}
J. Hor\'ak, P. J. McKenna, and W. Reichel, Very weak solutions with boundary singularities for semilinear elliptic Dirichlet problems in domains with conical corners, J. Math. Anal. Appl. \textbf{352} (2009), no. 1, 496--514.

\bibitem{hormander}
L. H\"ormander, The Analysis of Linear Partial Differential Operators I.
Distribution Theory and Fourier Analysis, Classics in Mathematics,
Springer-Verlag, Berlin, 2003, reprint of the second (1990) edition.

\bibitem{ios}
K. Ishige, S. Okabe, and T. Sato, A supercritical scalar field equation with a forcing term, J. Math. Pures Appl. (9) \textbf{128} (2019), 183--212.

\bibitem{kat-cone}
S. Katayama, Supercritical Lane--Emden equation on a cone with an inhomogeneous Dirichlet boundary condition, Nonlinear Differ. Equ. Appl. \textbf{32} (2025), article no.~120.

\bibitem{kato}
T. Kato, Perturbation Theory for Linear Operators, Classics in Mathematics,
Springer-Verlag, Berlin, 1995, reprint of the 1980 edition.

\bibitem{lb}
Y. Li and J. G. Bao, Semilinear elliptic system with boundary singularity, Discrete Contin. Dyn. Syst. \textbf{40} (2020), no. 4, 2189--2212.

\bibitem{lm}
J.-L. Lions and E. Magenes, Non-Homogeneous Boundary Value Problems and
Applications I, Grundlehren der mathematischen Wissenschaften 181,
Springer-Verlag, Berlin, 1972.

\bibitem{mp}
R. Mazzeo and F. Pacard, A construction of singular solutions for a semilinear elliptic equation using asymptotic analysis, J. Differential Geom. \textbf{44} (1996), no. 2, 331--370.

\bibitem{mr}
P. J. McKenna and W. Reichel, A priori bounds for semilinear equations and a new class of critical exponents for Lipschitz domains, J. Funct. Anal. \textbf{244} (2007), no. 1, 220--246.

\bibitem{p}
F. Pacard, Existence and convergence of positive weak solutions of $-\Delta u=u^{\frac{n}{n-2}}$ in bounded domains of $\R^n$, $n\geq 3$, Calc. Var. Partial Differential Equations \textbf{1} (1993), no. 3, 243--265.

\bibitem{s}
P. Souplet, Optimal regularity conditions for elliptic problems via $L^p_\delta$-spaces, Duke Math. J. \textbf{127} (2005), no. 1, 175--192.

\bibitem{wx}
H.-Y. Wang and J. Xiong, The critical semilinear elliptic equation with isolated boundary singularities II, arXiv:2608.16106v2, 2026.

\bibitem{x}
J. Xiong, The critical semilinear elliptic equation with isolated boundary singularities, J. Differential Equations \textbf{263} (2017), no. 3, 1907--1930.

\end{thebibliography}
\end{document}